\documentclass[10pt]{article}

\usepackage{a4wide}
\usepackage{amssymb}
\usepackage{amsfonts}
\usepackage{amsmath}
\input xy
\xyoption{arrow} \xyoption{matrix}

\date{}

\newtheorem{proposition}{Proposition}[section]
\newtheorem{theorem}[proposition]{Theorem}
\newtheorem{lemma}[proposition]{Lemma}

\newtheorem{corollary}[proposition]{Corollary}

\def\GK{{\rm  GK}\,}

\def\der{\partial }

\def\nFM0{{\nu }_{F,M_0}}
\def\nFN0{{\nu }_{F,N_0}}
\def\nGN0{{\nu }_{G,N_0}}

\def\N0{ {\bf N}_0 }

\def\t{\otimes}
\def\g{\gamma}
\def\v{\varphi}
\def\ra{\rightarrow}

\def\lra{\leftrightarrow}
\def\Xpm{X^{\pm }}

\def\s{\sigma}
\def\Z{\mathbb{Z}}

\def\l1{{\lambda}_1}

\def\a{\alpha}
\def\a0{ {\alpha }_0}
\def\a1{ {\alpha }_1}

\def\l{\lambda}
\def\o{\omega}

\def\nFGM0{{\nu }_{F,G,M_0}}

\def\nFN0{{\nu}_{F,N_0}}

\def\sm{{\sigma}^m}

\def\sm1{{\sigma}^{-1}}

\def\smtp1{{\sigma}^{-t+1}}

\def\o{\omega }
\def\S1{S^{-1}}

\def\Xpm1{X^{\pm 1}_1}

\def\sPM1{{\sigma }^{\pm 1}}
\def\sMP1{{\sigma }^{\mp 1 }}

\def\d{\delta}

\def\di{{\rm d.ind}}

\def\L{\Lambda}

\def\CD{{\cal D}}

\def\Ytm1{Y^{t-1}}
\def\Yim1{Y^{i-1}}

\def\CL{{\cal L}}
\def\CM{{\cal M}}
\def\CN{{\cal N}}

\def\CF{{\cal F}}

\def\CH{{\cal H}}

\def\Aut{{\rm Aut}}

\def\ad{{\rm ad }}
\def\dim{{\rm dim }}

\def\ker{ {\rm ker } }

\def\CJ{ {\cal J}}

\def\SL2Z{ {\rm SL}_2({\bf Z}) }

\def\CL{{\cal L}}

\def\Gp1{ G^{1 , 1 } }
\def\P11{ P^{-1 , 1 } }
\def\Pp1{ P^{1 , 1 } }

\def\Supp{{\rm Supp}}

\def\vol{{\rm vol}}

\def\CV{{\cal V}}

\def\nCLsr{{}^\nu\kern-2pt {\cal L}^{\sigma , \rho  }}
\def\nP{{}^\nu \kern-2pt P}
\def\nL{{}^\nu\kern-2pt L}
\def\nLL{{}^\nu\kern-2pt \Lambda}
\def\nPsr{{}^\nu\kern-2pt P^{\sigma , \rho  }}
\def\nLsr{{}^\nu\kern-2pt L^{\sigma , \rho  }}
\def\nuCL{{}^\nu\kern-2pt  {\cal L}}
\def\nCLsr{{}^\nu\kern-2pt {\cal L}^{\sigma , \rho  }}
\def\nCL1m{{}^\nu\kern-2pt {\cal L}^{-1 , 1  }}
\def\x1nu{x^\frac{1}{\nu}}
\def\xm1nu{x^{-\frac{1}{\nu}}}

\def\ob{\overline{b}}

\def\CN{{\cal N}}
\def\ra{\rightarrow }

\def\CB{{\cal B}}

\def\coker{{\rm coker}}

\def\CC{ {\cal C}}

\def\CH{ {\cal H}}

\def\nAM0{{\nu }_{{\cal A},M_0}}
\def\nAN0{{\nu }_{{\cal A},N_0}}

\def\End{ {\rm End }}

\def\CJ{ {\cal J }}

\def\det{ {\rm det }}
\def\ad{ {\rm ad }}

\def\spec{{\rm spec}}

\def\ga{\mathfrak{a}}
\def\gb{\mathfrak{b}}

\def\gn{\mathfrak{n}}

\def\gp{\mathfrak{p}}
\def\gq{\mathfrak{q}}

\def\GL{{\rm GL}}
\def\SL{{\rm SL}}

\def\Spec{{\rm Spec}}

\def\di!{\frac{\der^i}{i!}}
\def\dik!{\frac{\der^k_i}{k!}}

\def\id{{\rm id}}

\def\Max{{\rm Max}}

\def\gldim{{\rm gldim}}

\def\Fun{{\rm Fun}}

\def\N{\mathbb{N}}

\def\0{\overline{0}}
\def\1{\overline{1}}

\def\Ln1{\L_{n,\overline{1}}}

\def\oa{\overline{a}}

\def\a1{a_{\overline{1}}}

\def\St{{\rm St}}
\def\S{\Sigma}

\def\vn1{\overrightarrow{n-1}}

\def\hx{\widehat{x}}

\def\im{{\rm im}}

\def\mA{\mathbb{A}}

\def\Sub{{\rm Sub}}
\def\SSub{{\rm SSub}}
\def\Inc{{\rm Inc}}
\def\Min{{\rm Min}}

\def\bu{\overline{u}}

\def\F{\mathbb{F}}

\def\Inn{{\rm Inn}}

\def\mS{\mathbb{S}}

\def\lann{{\rm l.ann}}
\def\rann{{\rm r.ann}}
\def\Cen{{\rm Cen}}
\def\clKdim{{\rm cl.Kdim}}

\def\hht{{\rm ht}}

\def\heta{\widehat{\eta}}
\def\mF{\mathbb{F}}
\def\mM{\mathbb{M}}
\def\mT{\mathbb{T}}
\def\ind{{\rm ind}}
\def\ud{\underline{d}}
\def\stCH{{\rm st}_{G_n}(\CH_1)}

\def\mU{\mathbb{U}}

\begin{document}

\author{V. V. \  Bavula 
}

\title{The group of automorphisms of the 
algebra of one-sided inverses of a polynomial algebra }

\maketitle

$$Dedicated\; to\; T.\; Lenagan \;on \; the\; occasion\; of\; his \;
60'th\; birthday $$

\begin{abstract}
The 
algebra $\mS_n$ in the title  is  obtained from a polynomial
algebra $P_n$ in $n$ variables by adding commuting,  {\em left}
(but not two-sided) inverses of the canonical generators of $P_n$.
Ignoring non-Noetherian property, the algebra $\mS_n$ belongs to a
family of algebras like the Weyl algebra $A_n$ and the polynomial
algebra $P_{2n}$. The group of automorphisms $G_n$ of the algebra
$\mS_n$ is found:
$$ G_n=S_n\ltimes \mT^n\ltimes \Inn (\mS_n) \supseteq S_n\ltimes \mT^n\ltimes \underbrace{\GL_\infty (K)\ltimes\cdots \ltimes
\GL_\infty (K)}_{2^n-1 \;\; {\rm times}}=:G_n' $$ where  $S_n$ is
the symmetric group, $\mT^n$ is the $n$-dimensional torus,  $\Inn
(\mS_n)$ is the group of inner automorphisms of $\mS_n$ (which is
huge), and $\GL_\infty (K)$ is  the group of invertible infinite
dimensional matrices. This result may help in understanding of the
structure of the groups of automorphisms of the Weyl algebra $A_n$
and the polynomial algebra $P_{2n}$. An analog of the {\em
Jacobian homomorphism}: $\Aut_{K-{\rm alg}}(P_{2n} )\ra K^*$, the
so-called  {\em global determinant} is introduced for the group
$G_n'$ (notice that the algebra $\mS_n$ is {\em noncommutative}
and neither left nor right Noetherian).



$${\bf Contents}$$
\begin{enumerate}
\item Introduction. \item Preliminaries on the 
algebras $\mS_n$. \item Certain subgroups of $\Aut_{K-{\rm
alg}}(\mS_n )$. \item The groups $\Aut_{K-{\rm alg}}(\mS_1)$ and
$\mS_1^*$. \item The group of automorphisms of the algebra
$\mS_n$. \item A membership criterion for elements of the algebra
$\mS_n$. \item The groups $\mM_n^*$ and $G_n'$.  \item An analogue
of the Jacobian map - the global determinant. \item Stabilizers in
$\Aut_{K-{\rm alg}}(\mS_n )$ of the prime or idempotent ideals of
$\mS_n$. \item Endomorphisms of the algebra $\mS_n$.
\end{enumerate}
\end{abstract}


\section{Introduction}
Throughout, ring means an associative ring with $1$; module means
a left module;
 $\N :=\{0, 1, \ldots \}$ is the set of natural numbers; $K$ is a
field of characteristic zero and  $K^*$ is its group of units;
$P_n:= K[x_1, \ldots , x_n]$ is a polynomial algebra over $K$;
$\der_1:=\frac{\der}{\der x_1}, \ldots , \der_n:=\frac{\der}{\der
x_n}$ are the partial derivatives ($K$-linear derivations) of
$P_n$.

$\noindent $

{\it Definition}, \cite{shrekalg}. The 
{\em algebra} $\mathbb{S}_n$ {\em of one-sided inverses} of $P_n$
is an algebra generated over a field $K$ of characteristic zero by
$2n$ elements $x_1, \ldots , x_n, y_n, \ldots , y_n$ that satisfy
the defining relations:
$$ y_1x_1=\cdots = y_nx_n=1 , \;\; [x_i, y_j]=[x_i, x_j]= [y_i,y_j]=0
\;\; {\rm for\; all}\; i\neq j,$$ where $[a,b]:= ab-ba$, the
commutator of elements $a$ and $b$.

$\noindent $

By the very definition, the algebra $\mS_n$ is obtained from the
polynomial algebra $P_n$ by adding commuting, left (or right)
inverses of its canonical generators. The algebra $\mS_1$ is a
well-known primitive algebra \cite{Jacobson-StrRing}, p. 35,
Example 2. Over the field
 $\mathbb{C}$ of complex numbers, the completion of the algebra
 $\mS_1$ is the {\em Toeplitz algebra} which is the
 $\mathbb{C}^*$-algebra generated by a unilateral shift on the
 Hilbert space $l^2(\N )$ (note that $y_1=x_1^*$). The Toeplitz
 algebra is the universal $\mathbb{C}^*$-algebra generated by a
 proper isometry.

$\noindent $

{\it Example}, \cite{shrekalg}. Consider a vector space $V=
\bigoplus_{i\in \N}Ke_i$ and two shift operators on $V$, $X:
e_i\mapsto e_{i+1}$ and $Y:e_i\mapsto e_{i-1}$ for all $i\geq 0$
where $e_{-1}:=0$. The subalgebra of $\End_K(V)$ generated by the
operators $X$ and $Y$ is isomorphic to the algebra $\mS_1$
$(X\mapsto x$, $Y\mapsto y)$. By taking the $n$'th tensor power
$V^{\t n }=\bigoplus_{\alpha \in \N^n}Ke_\alpha$ of $V$ we see
that the algebra $\mS_n$ is isomorphic to the subalgebra of
$\End_K(V^{\t n })$ generated by the $2n$ shifts $X_1, Y_1, \ldots
, X_n, Y_n$ that act in different directions.

$\noindent $

It is an experimental fact \cite{shrekalg} that the algebra
$\mS_1$ has properties that are a mixture of the properties of the
polynomial algebra $P_2$ in two variable and the {\em first Weyl }
algebra $A_1$, which is not surprising when we look at their
defining relations:
\begin{eqnarray*}
 P_2:& yx-xy=0; \\
 A_1:&yx-xy=1;\\
 \mS_1:& yx=1.
\end{eqnarray*}
The same is true for their higher analogues: $P_{2n}=P_2^{\t n}$,
$A_n:=A_1^{\t n }$ (the $n$'th {\em Weyl } algebra), and $\mS_n =
\mS_1^{\t n }$. For example,
\begin{eqnarray*}
 \clKdim (\mS_n )&\stackrel{\cite{shrekalg}}{=} &2n=\clKdim (P_{2n}), \\
\gldim (\mS_n )&\stackrel{\cite{shrekalg}}{=} &n=\gldim (A_n), \\
\GK (\mS_n )&\stackrel{\cite{shrekalg}}{=} &2n=\GK (A_n) = \GK(P_{2n}), \\
\end{eqnarray*}
where $\clKdim$, $\gldim$, and $\GK$ stand for the classical Krull
dimension, the global homological dimension, and the
Gelfand-Kirillov dimension respectively. The big difference
between the algebra $\mS_n$ and the algebras $P_{2n}$ and $A_n$ is
that $\mS_n$ is neither left nor right Noetherian and is not a
domain either.

The algebras $\mS_n$ are fundamental non-Noetherian algebras, they
are universal non-Noetherian algebras of their own kind in a
similar way as the polynomial algebras are universal in the class
of all the commutative algebras and the Weyl algebras are
universal in the class of algebras  of differential operators.

The algebra $\mS_n$ often appears as a subalgebra or a factor
algebra  of many non-Noetherian algebras. For example, $\mS_1$ is
a factor algebra of certain non-Noetherian down-up algebras as was
shown by Jordan  \cite{Jordan-Down-up}  (see also Benkart and Roby
\cite{Benk-Roby}; Kirkman, Musson, and Passman
\cite{Kirk-Mus-Kuz}; Kirkman and Kuzmanovich \cite{Kirk-Kuz}); and
$\mS_n$ is a subalgebra of the Jacobian algebra $\mA_n$ (see
below) \cite{Bav-Jacalg}.

The aim of this paper is to find the group $G_n:=\Aut_{K-{\rm
alg}}(\mS_n)$ of automorphisms of the algebra $\mS_n$.
\begin{itemize}
\item (Theorem \ref{5Feb9}) $ G_n=S_n\ltimes \mT^n\ltimes \Inn
(\mS_n)$. \item (Lemma \ref{x1Apr9}) $G_n\supseteq G_n':=
S_n\ltimes \mT^n\ltimes \underbrace{\GL_\infty (K)\ltimes\cdots
\ltimes \GL_\infty (K)}_{2^n-1 \;\; {\rm times}}$,
\end{itemize}
where  $S_n$ is the symmetric group, $\mT^n$ is the
$n$-dimensional torus, $\Inn (\mS_n)$ is the group of inner
automorphisms of the algebra $\mS_n$, and $\GL_\infty (K)$ is the
group of all the invertible infinite dimensional matrices of the
type $1+M_\infty (K)$ where the algebra (without 1) of infinite
dimensional matrices $M_\infty (K) :=\varinjlim
M_d(K)=\bigcup_{d\geq 1}M_d(K)$ is the injective limit of matrix
algebras. A semi-direct product $H_1\ltimes H_2\ltimes \cdots
\ltimes H_m$ of several groups means that  $H_1\ltimes (H_2\ltimes
( \cdots \ltimes (H_{m-1}\ltimes  H_m)\cdots )$.

The proof of Theorem \ref{5Feb9} is rather long (and non-trivial)
and based upon several results proved in this paper (and in
\cite{shrekalg}) which are interesting on their own. Let me
explain briefly the logical structure of the proof. There are two
cases to consider when $n=1$ and $n>1$. The proofs of both cases
are based on  different ideas. The case $n=1$ is a kind of a
degeneration of the second case and is much more easier. The key
point in finding the group $G_1$ is to use the {\em index} of
linear maps in infinite dimensional vector spaces and the fact
that each automorphism of the algebra $\mS_n$ is determined by its
action on the  set $\{ x_1, \ldots , x_n\}$ (or  $\{ y_1, \ldots ,
y_n\}$):

\begin{itemize}
\item (Theorem \ref{6Feb9}) (Rigidity of the group $G_n$)  {\em
Let $\s , \tau \in G_n$. Then the following statements are
equivalent.}
\begin{enumerate}
 \item $\s = \tau$. \item $\s (x_1) = \tau (x_1), \ldots , \s (x_n) = \tau
 (x_n)$.
 \item$\s (y_1) = \tau (y_1), \ldots , \s (y_n) = \tau (y_n)$.
\end{enumerate}
\end{itemize}

For $n>1$, one of the key ideas in finding the group $G_n$ is  to
use the action of the group $G_n$ on the set $\CH_1$ of all the
height 1 prime ideals of the algebra $\mS_n$. The set $\CH_1=\{
\gp_1, \ldots , \gp_n\}$ is finite and is found in
\cite{shrekalg}. It follows that the group
$$G_n=S_n\ltimes \St_{G_n}(\CH_1)$$
is the semi-direct product of the symmetric group $S_n$ and the
stabilizer of the set $\CH_1$ in $G_n$,
$$ \St_{G_n}(\CH_1):=\{ \s \in G_n\, | \, \s (\gp_1) = \gp_1,
\ldots , \s (\gp_n) = \gp_n\}.$$ The group $\St_{G_n}(\CH_1)$
contains the $n$-dimensional torus $\mT^n$. Using a Membership
Criterion (Theorem \ref{10Feb9}) for elements of the algebra
$\mS_n$, it follows that
$$ \St_{G_n}(\CH_1) = \mT^n\ltimes \stCH$$
where 
\begin{equation}\label{stCH}
\stCH =\{ \s \in \St_{G_n}(\CH_1)\, | \, \s (x_i) \equiv x_i \mod
\gp_i, \s (y_i) \equiv y_i\mod \gp_i, i=1, \ldots , n\}.
\end{equation}
Moreover,
\begin{itemize}
\item (Corollary \ref{a21Feb9}) $\stCH =\Inn (\mS_n)$.
\end{itemize}
One of the key points of the proof of Theorem \ref{5Feb9} and
Corollary  \ref{a21Feb9}  is the fact that
\begin{itemize}
\item (\cite{shrekalg}, Corollary 3.3):  $P_n$ {\em is the only
simple, faithful $\mS_n$-module (up to isomorphism),}
\end{itemize}
and so the algebra $\mS_n$ can be seen as the subalgebra of the
endomorphism algebra $E_n:=\End_K(P_n)$ of all the linear maps
from the vector space $P_n$ to itself and we can visualize the
group $G_n$ via the group $\Aut_K(P_n)$ of units of the algebra
$E_n$ as follows:
\begin{itemize}
\item (Theorem \ref{a5Feb9}) $G_n=\{ \s_\v \, | \, \v \in
\Aut_K(P_n)$ {\em such that} $\v\mS_n \v^{-1}= \mS_n\}$ {\em
where} $\s_\v (a):= \v a \v^{-1}$, $a\in \mS_n$.
\end{itemize}
To represent the group $G_n$ via linear maps in an infinite
dimensional space helps not much unless we have a criterion of
when a linear map represents an element of the group $G_n$ (or an
element of  the algebra $\mS_n$). Several membership criteria are
proved in Section \ref{MEMB} which are used at the final stage of
the proof of Theorem \ref{5Feb9}:

\begin{itemize}
\item (Theorem \ref{10Feb9}) {\em  Let $\v \in \End_K(P_n)$. Then
} $\v \in \mS_n$ iff $[x_1, \v]\in \v\gp_1 +\gp_1, \ldots , [x_n,
\v ] \in \v \gp_n +\gp_n.$\item (Corollary \ref{a15Feb9}) {\em Let
$F_n:=\gp_1\cdots \gp_n$. Then} $$\{ \v \in \End_K(P_n)\, | \,
[x_i, \v ] \in F_n, [y_i, \v ] \in F_n, i=1, \ldots , n\}
=\begin{cases}
\mS_1& \text{if }n=1,\\
K+F_n& \text{if }n>1.\\
\end{cases} $$
\end{itemize}
The structure of the group $G_1=\mT^1\ltimes \GL_\infty (K)$ is
yet another confirmation of `similarity' of the algebras $P_2$,
$A_1$, and $\mS_1$. The groups of automorphisms of the polynomial
algebra $P_2$ and the Weyl algebra $A_1$ were found by Jung
\cite{jung}, Van der Kulk \cite{kulk}, and Dixmier \cite{Dix}
respectively. These two groups have almost identical structure,
they are `infinite $GL$-groups' in the sense that they are
generated by the torus $\mT^1$ and by the obvious automorphisms:
$x\mapsto x+\l y^i$, $y\mapsto y$; $x\mapsto x$, $y\mapsto y+ \l
x^i$, where $i\in \N$ and $\l \in K$; which are sort of
`elementary infinite dimensional matrices' (i.e. `infinite
dimensional transvections`). The same picture as for the group
$G_1$. In prime characteristic, the group of automorphism of the
Weyl algebra $A_1$ was found by Makar-Limanov \cite{Mak-LimBSMF84}
(see also Bavula \cite{A1rescen}  for a different approach and for
further developments). More on polynomial automorphisms the reader
can find in the book of Van den Essen \cite{vdEssen-book}.

$\noindent $

There is an important homomorphism from the group $\Aut_{K-{\rm
alg}}(P_{2n})$ of automorphisms of the polynomial algebra $P_{2n}$
to the group $K^*$,  the so-called {\em Jacobian} (map or
homomorphism):
$$ \CJ : \Aut_{K-{\rm alg}}(P_{2n})\ra K^*, \;\; \s \mapsto \det (\frac{\der \s (x_i)}{\der
x_j}).$$ Note that the Jacobian homomorphism is a determinant. In
this paper (Section \ref{JACDET}), its analogue is introduced for
the algebra $\mS_n$ which is called the {\em global determinant}:
$$\det : G_n'\ra K^*, \;\; \s \mapsto \det (\s ).$$ It is  a group
homomorphism (Corollary  \ref{e7Mar9})  which is defined as
follows. By Lemma \ref{x1Apr9}, each element $\s$ of $G_n'$ is a
unique product $\s = \tau t_\l \s_1\cdots \s_{2^n-1}$ where $\tau
\in S_n$, $t_\l \in \mT^n$,  $\l = (\l_1, \ldots , \l_n)\in
K^{*n}$, and $\s_i\in \GL_\infty (K)$. Then 
\begin{equation}\label{detgen}
\det (\s ) := {\rm sgn} (\tau ) \cdot \prod_{i=1}^n \l_i\cdot
\prod_{j=1}^{2^n-1}\det (\s_j)
\end{equation}
where ${\rm sgn} (\tau )$ is the parity of the permutation $\tau$
and $\det (\s_j)$ is the `usual' determinant of the element $\s_j$
 of the group $\GL_\infty (K)$. It is an interesting question of
 whether it is possible to extend the global determinant to the
 group $G_n$.

$\noindent $

The paper is organized as follows. In Section \ref{PREL}, some
useful results from \cite{shrekalg} are collected which are used
later.

In Section \ref{CSOA}, several subgroups of the group $G_n$ are
introduced, a useful description (Theorem \ref{a5Feb9}) of the
group $G_n$  is given, and a criterion of equality of two elements
of the group $G_n$ is proved (Theorem \ref{6Feb9}).

In Section \ref{TG1S}, the group $G_1$ is found (Theorem
\ref{A5Feb9}).

In Section \ref{GAS}, the group $G_n$ is found (Theorem
\ref{5Feb9}). Several corollaries are obtained.  It is proved that
the groups $G_n$ and $\Inn (\mS_n)$ have trivial centre (Corollary
\ref{b24Feb9}).

In Section \ref{MEMB}, several Membership Criteria are proved for
the algebras $\mS_n$, $P_n+F_n$ and $K+F_n$ (Theorem \ref{10Feb9},
Corollaries \ref{a13Feb9} and \ref{a15Feb9}).

In Section \ref{JACDET}, the global determinant is extended to a
certain monoid $S_n\ltimes \mT^n\ltimes \mM_n$, the group of units
of which is isomorphic to the group $G_n'$ (Corollary
\ref{b11Mar9}.(1)). Moreover,
\begin{itemize}
\item (Corollary \ref{b11Mar9}.(2)) $G_n'\simeq \{ a\in S_n\ltimes
\mT^n\ltimes \mM_n\, | \, \det (a) \neq 0\}$.
\end{itemize}
Intuitively, the pair $(S_n\ltimes \mT^n\ltimes \mM_n, G_n')$, a
monoid and its group of units,  is an infinite dimensional
analogue of the pair $(M_d(K), \GL_d(K))$. Theorem
\ref{d7Mar9}.(3) shows that   the global determinant can be
computed effectively (in finitely many steps).

In Section \ref{STABSN}, the stabilizers in the group $G_n$ of
several classes of ideals of the algebra $\mS_n$ are computed. In
particular, the stabilizers of all the prime ideals of $\mS_n$ are
found (Corollary \ref{b10Feb9}.(2) and Corollary \ref{a8Mar9}).

The ideal $\ga_n:= \gp_1+\cdots + \gp_n$ is a prime, idempotent
ideal of the algebra $\mS_n$ of height $n$, \cite{shrekalg}.

\begin{itemize}
\item (Theorem \ref{B10Feb9}) {\em The ideal $\ga_n$ is the only
nonzero, prime, $G_n$-invariant ideal of the algebra $\mS_n$.}
\item  (Theorem \ref{10Mar9}) {\em Let $\gp$ be a prime ideal of
$\mS_n$. Then its stabilizer $\St_{G_n}(\gp )$ is a maximal
subgroup of the group $G_n$ iff $n>1$ and $\gp$ is of height 1,
and, in this case, $[G_n:\St_{G_n}(\gp )]=n$. }\item (Corollary
\ref{a10Mar9}) {\em Let $\ga$ be a proper ideal of $\mS_n$. Then
its stabilizer $\St_{G_n}(\ga )$ has finite index in the group
$G_n$ iff $\ga^2= \ga$.} \item (Corollary \ref{aA10Feb9}) {\em If
$\ga$ is a generic idempotent ideal of $\mS_n$ then its stabilizer
is written via the wreath products of symmetric groups:}
$$  \St_{G_n}(\ga )= (S_m\times \prod_{i=1}^t(S_{h_i}\wr
S_{n_i}))\ltimes \mT^n\ltimes \Inn (\mS_n). $$
\end{itemize}

In Section \ref{ENDSN}, we classify all the algebra endomorphisms
of $\mS_n$ that stabilize the elements $x_1, \ldots , x_n$ and
show that each such endomorphism is a {\em monomorphism} but {\em
not} an isomorphism provided it is not the identity map (Corollary
\ref{a7Feb9}). Therefore, an analogous question to the Question of
Dixmier, namely, {\em is a monomorphism of the algebra $\mS_n$ is
an automorphism?} has  a negative answer. The original
Question/Problem of Dixmier states \cite{Dix}: {\em is every
homomorphism of the Weyl algebra $A_n$ an automorphism?} The Weyl
algebra $A_n$ is a simple algebra, so any homomorphism is
automatically a monomorphism. In \cite{Dix}, Dixmier poses this
question only for the first Weyl algebra $A_1$.


\section{Preliminaries on the 
algebras $\mS_n$}\label{PREL}

In this section, we collect some results without proofs on the
algebras $\mS_n$ from \cite{shrekalg} that will be used in this
paper, their proofs can be found in \cite{shrekalg}.

  Clearly,
$\mathbb{S}_n=\mS_1(1)\t \cdots \t\mS_1(n)\simeq \mathbb{S}_1^{\t
n}$ where $\mS_1(i):=K\langle x_i,y_i \, | \, y_ix_i=1\rangle
\simeq \mS_1$ and $$\mS_n=\bigoplus_{\alpha , \beta \in \N^n}
Kx^\alpha y^\beta$$ where $x^\alpha := x_1^{\alpha_1} \cdots
x_n^{\alpha_n}$, $\alpha = (\alpha_1, \ldots , \alpha_n)$,
$y^\beta := y_1^{\beta_1} \cdots y_n^{\beta_n}$, $\beta =
(\beta_1,\ldots , \beta_n)$. In particular, the algebra $\mS_n$
contains two polynomial subalgebras $P_n$ and $Q_n:=K[y_1, \ldots
, y_n]$ and is equal,  as a vector space,  to their tensor product
$P_n\t Q_n$. Note that also the Weyl algebra $A_n$ is a tensor
product (as a vector space) $P_n\t K[\der_1, \ldots , \der_n]$ of
its two polynomial subalgebras.

When $n=1$, we usually drop the subscript `1' if this does not
lead to confusion.  So, $\mS_1= K\langle x,y\, | \,
yx=1\rangle=\bigoplus_{i,j\geq 0}Kx^iy^j$. For each natural number
$d\geq 1$, let $M_d(K):=\bigoplus_{i,j=0}^{d-1}KE_{ij}$ be the
algebra of $d$-dimensional matrices where $\{ E_{ij}\}$ are the
matrix units, and
$$M_\infty (K) :=
\varinjlim M_d(K)=\bigoplus_{i,j\in \N}KE_{ij}$$ be the algebra
(without 1) of infinite dimensional matrices. The algebra $\mS_1$
contains the ideal $F:=\bigoplus_{i,j\in \N}KE_{ij}$, where
\begin{equation}\label{Eijc}
E_{ij}:= x^iy^j-x^{i+1}y^{j+1}, \;\; i,j\geq 0.
\end{equation}
For all natural numbers $i$, $j$, $k$, and $l$,
$E_{ij}E_{kl}=\d_{jk}E_{il}$ where $\d_{jk}$ is the Kronecker
delta function.  The ideal $F$ is an algebra (without 1)
isomorphic to the algebra $M_\infty (K)$ via $E_{ij}\mapsto
E_{ij}$. For all $i,j\geq 0$, 
\begin{equation}\label{xyEij}
xE_{ij}=E_{i+1, j}, \;\; yE_{ij} = E_{i-1, j}\;\;\; (E_{-1,j}:=0),
\end{equation}
\begin{equation}\label{xyEij1}
E_{ij}x=E_{i, j-1}, \;\; E_{ij}y = E_{i, j+1} \;\;\;
(E_{i,-1}:=0).
\end{equation}
\begin{equation}\label{mS1d}
\mS_1= K\oplus xK[x]\oplus yK[y]\oplus F,
\end{equation}
the direct sum of vector spaces. Then 
\begin{equation}\label{mS1d1}
\mS_1/F\simeq K[x,x^{-1}]=:L_1, \;\; x\mapsto x, \;\; y \mapsto
x^{-1},
\end{equation}
since $yx=1$, $xy=1-E_{00}$ and $E_{00}\in F$.

$\noindent $

The algebra $\mS_n = \bigotimes_{i=1}^n \mS_1(i)$ contains the
ideal
$$F_n:= F^{\t n }=\bigoplus_{\alpha , \beta \in
\N^n}KE_{\alpha \beta}, \;\; {\rm where}\;\; E_{\alpha
\beta}:=\prod_{i=1}^n E_{\alpha_i \beta_i}(i).$$ Note that
$E_{\alpha \beta}E_{\g \rho}=\d_{\beta \g }E_{\alpha  \rho}$ for
all elements $\alpha, \beta , \g , \rho \in \N^n$ where $\d_{\beta
 \g }$ is the Kronecker delta function.
\begin{itemize}
 \item
{\em $F_n a \neq 0$ and $aF_n\neq 0$ for all nonzero elements $a$
of the algebra $\mS_n$.} \item  {\em $F_n$ is the smallest (with
respect to inclusion) nonzero ideal of the algebra $\mS_n$ (i.e.
$F_n$ is contained in all nonzero ideals of $\mS_n$); $F_n^2=
F_n$; $F_n$ is an essential left and right submodule of $\mS_n$;
$F_n$ is the socle of the left and right $\mS_n$-module $\mS_n$;
$F_n$ is the socle of the $\mS_n$-bimodule $\mS_n$ and $F_n$ is a
simple $\mS_n$-bimodule.} \end{itemize}

{\bf The involution $\eta$ on $\mS_n$}. The algebra $\mS_n$ admits
the {\em involution}
$$ \eta : \mS_n\ra \mS_n, \;\; x_i\mapsto y_i, \;\; y_i\mapsto
x_i, \;\; i=1, \ldots , n,$$ i.e. it is a $K$-algebra
anti-isomorphism ($\eta (ab) = \eta (b) \eta (a)$ for all $a,b\in
\mS_n$) such that $\eta^2 = \id_{\mS_n}$, the identity map on
$\mS_n$. So, the algebra $\mS_n$ is {\em self-dual} (i.e. it is
isomorphic to its opposite algebra, $\eta : \mS_n\simeq
\mS_n^{op}$). The involution $\eta$ acts on the `matrix' ring
$F_n$ as the transposition,  
\begin{equation}\label{eEij1}
\eta (E_{\alpha \beta} )=E_{\beta \alpha}.
\end{equation}
 The canonical generators $x_i$,
$y_j$ $(1\leq i,j\leq n)$ determine the ascending filtration $\{
\mS_{n, \leq i}\}_{i\in \N}$ on the algebra $\mS_n$ in the obvious
way (i.e. by the total degree of the generators): $\mS_{n, \leq
i}:= \bigoplus_{|\alpha |+|\beta |\leq i} Kx^\alpha y^\beta$ where
$|\alpha | \; = \alpha_1+\cdots + \alpha_n$ ($\mS_{n, \leq
i}\mS_{n, \leq j}\subseteq \mS_{n, \leq i+j}$ for all $i,j\geq
0$). Then $\dim (\mS_{n,\leq i})={i+2n \choose 2n}$ for $i\geq 0$,
and so the Gelfand-Kirillov dimension $\GK (\mS_n )$ of the
algebra $\mS_n$ is equal to $2n$. It is not difficult to show
 that the algebra $\mS_n$ is neither left nor
right Noetherian. Moreover, it contains infinite direct sums of
left and right ideals (see \cite{shrekalg}).

\begin{itemize}
 \item   {\em The
 algebra $\mS_n$ is central, prime, and catenary. Every nonzero
 ideal  of $\mS_n$ is an essential left and right submodule of}
 $\mS_n$.
 \item  {\em The ideals of
 $\mS_n$ commute ($IJ=JI$);  and the set of ideals of $\mS_n$ satisfy the a.c.c..}
 \item  {\em The classical Krull dimension $\clKdim (\mS_n)$ of $\mS_n$ is $2n$.}
  \item  {\em Let $I$
be an ideal of $\mS_n$. Then the factor algebra $\mS_n / I$ is
left (or right) Noetherian iff the ideal $I$ contains all the
height one primes of $\mS_n$.}
\end{itemize}

{\bf The set of height 1 primes of $\mS_n$}.   Consider the ideals
of the algebra $\mS_n$:
$$\gp_1:=F\t \mS_{n-1},\; \gp_2:= \mS_1\t F\t \mS_{n-2}, \ldots ,
 \gp_n:= \mS_{n-1} \t F.$$ Then $\mS_n/\gp_i\simeq
\mS_{n-1}\t (\mS_1/F) \simeq  \mS_{n-1}\t K[x_i, x_i^{-1}]$ and
$\bigcap_{i=1}^n \gp_i = \prod_{i=1}^n \gp_i =F^{\t n }$. Clearly,
$\gp_i\not\subseteq \gp_j$ for all $i\neq j$.

\begin{itemize}
 \item
{\em The set $\CH_1$ of height 1 prime ideals of the algebra
$\mS_n$ is} $\{ \gp_1, \ldots , \gp_n\}$.
\end{itemize}

 Let
$\ga_n:= \gp_1+\cdots +\gp_n$. Then the factor algebra
\begin{equation}\label{SnSn}
\mS_n/ \ga_n\simeq (\mS_1/F)^{\t n } \simeq \bigotimes_{i=1}^n
K[x_i, x_i^{-1}]= K[x_1, x_1^{-1}, \ldots , x_n, x_n^{-1}]=:L_n
\end{equation}
is a skew Laurent polynomial algebra in $n$ variables,  and so
$\ga_n$ is a prime ideal of height and co-height $n$ of the
algebra $\mS_n$. The algebra $L_n$ is commutative, and so
\begin{equation}\label{comab}
 [a,b]\in \ga_n\;\; {\rm for \; all}\;\; a,b\in \mS_n .
\end{equation}
That is $[\mS_n , \mS_n]\subseteq \ga_n$. In particular, $[\mS_1,
\mS_1]\subseteq F$. Since $\eta (\ga_n) = \ga_n$, the involution
of the algebra $\mS_n$ induces the {\em automorphism}
$\overline{\eta}$ of the factor algebra $\mS_n / \ga_n$ by the
rule: 
\begin{equation}\label{1SnSn}
\overline{\eta}: L_n\ra L_n, \;\; x_i\mapsto x_i^{-1}, \;\; i=1,
\ldots , n.
\end{equation}
It follows that $\eta (ab) - \eta (a) \eta (b)\in \ga_n$ for all
elements $a, b\in \mS_n$.


\section{Certain subgroups of $\Aut_{K-{\rm alg}}(\mS_n
)$}\label{CSOA}

Recall that  $G_n:=\Aut_{K-{\rm alg}}(\mS_n)$ is the group of
automorphisms
of the 
algebra $\mS_n$. In this section, a useful description of the
group $G_n$ is given (Theorem \ref{a5Feb9}), an important (rather
peculiar) criterion of the equality of two elements of $G_n$
(Theorem \ref{6Feb9}) is found, and several subgroups of $G_n$ are
introduced that are building blocks of the group $G_n$. These
results are important in finding the group $G_n$.

\begin{proposition}\label{a19Dec8}
\cite{shrekalg} The polynomial algebra $P_n$
 is the only faithful, simple $\mS_n$-module.
\end{proposition}

In more detail, ${}_{\mS_n}P_n\simeq \mS_n / (\sum_{i=0}^n \mS_n
y_i) =\bigoplus_{\alpha \in \N^n} Kx^\alpha \overline{1}$,
$\overline{1}:= 1+\sum_{i=1}^n \mS_ny_i$; and the action of the
canonical generators of the algebra $\mS_n$ on the polynomial
algebra $P_n$ is given by the rule:
$$ x_i*x^\alpha = x^{\alpha + e_i}, \;\; y_i*x^\alpha = \begin{cases}
x^{\alpha - e_i}& \text{if } \; \alpha_i>0,\\
0& \text{if }\; \alpha_i=0,\\
\end{cases}  \;\; {\rm and }\;\; E_{\beta \g}*x^\alpha = \d_{\g
\alpha} x^\beta,
$$
where $e_1:= (1,0,\ldots , 0),  \ldots , e_n:=(0, \ldots , 0,1)$
is the canonical basis for the free $\Z$-module
$\Z^n=\bigoplus_{i=1}^n \Z e_i$.  We identify the algebra $\mS_n$
with its image in the algebra $\End_K(P_n)$ of all the $K$-linear
maps from the vector space $P_n$ to itself, i.e. $\mS_n \subset
\End_K(P_n)$. Let $\Aut_K(P_n)$ be the group of units of the
algebra $\End_K(P_n)$. $\Aut_K(P_n)$ is the group of all the {\em
invertible} $K$-linear maps from $P_n$ to itself. Each element $\v
\in \Aut_K(P_n)$ yields an inner automorphism $\o_\v : f\mapsto \v
f\v^{-1}$ of the algebra $\End_K(P_n)$. Suppose that the
automorphism $\o_\v$ respects the subalgebra $\mS_n$, that is
$\o_\v (\mS_n ) = \mS_n$, then its restriction $\s_\v : \o_\v
|_{\mS_n} : a\mapsto \v a\v^{-1}$ is an automorphism of the
algebra $\mS_n$.

$\noindent $

The next result shows that all the automorphisms of the algebra
$\mS_n$ can be obtained in this way.

\begin{theorem}\label{a5Feb9}
$G_n=\{ \s_\v \, | \, \v \in \Aut_K(P_n)$ such that $\v\mS_n
\v^{-1}= \mS_n\}$ where $\s_\v (a):= \v a \v^{-1}$, $a\in \mS_n$.
\end{theorem}

{\it Proof}. Let $\s \in G_n$. The twisted by the automorphism
$\s$ module $P_n$, denoted ${}^{\s}P_n$, is simple and faithful.
Recall that as a vector space the module ${}^{\s}P_n$ coincides
with the module $P_n$ but the action of the algebra $\mS_n$ is
given by the rule: $a\cdot p :=\s (a) *p$ where $a\in \mS_n$ and
$p\in P_n$. By Proposition \ref{a19Dec8}, the $\mS_n$-modules
$P_n$ and ${}^{\s}P_n$ are isomorphic. So, there  exists an
element $\v \in \Aut_K(P_n)$ such that $\v a = \s (a)\v$ for all
$a\in \mS_n$, and so $\s (a) = \v a \v^{-1}$, as required. $\Box $

\begin{theorem}\label{24Dec8}
\cite{shrekalg} The ideal $\ga_n$ is the smallest ideal of the
algebra $\mS_n$ such that the factor algebra $\mS_n / \ga_n$ is
commutative.
\end{theorem}

\begin{lemma}\label{b5Feb9}
$\s (\ga_n) = \ga_n$ for all $\s \in G_n$.
\end{lemma}

{\it Remark}. We will see that the ideal $\ga_n$ is the {\em only}
nonzero, prime, $G_n$-invariant ideal of the algebra $\mS_n$
(Theorem \ref{B10Feb9}).

$\noindent $

 {\it Proof}. For each element $\s \in
G_n$, the map
$$ \mS_n / \ga_n \ra \s (\mS_n)/ \s (\ga_n),\;\; a+\ga_n \mapsto \s
(a) +\s (\ga_n), $$ is an isomorphism of  algebras. By Theorem
\ref{24Dec8}, $\s (\ga_n) = \ga_n$ for all $\s \in G_n$ since
$\mS_n / \ga_n$ is a commutative algebra. $\Box $

$\noindent $

{\bf The automorphism $\heta\in \Aut (G_n)$}. The involution
$\eta$ of the algebra $\mS_n$ yields the automorphism $\heta\in
\Aut (G_n)$ of the group $G_n$: 
\begin{equation}\label{etah}
\heta : G_n\ra G_n, \;\; \s \mapsto \eta \s \eta^{-1}.
\end{equation}
Clearly, $\heta^2= e$ and $\heta (\s ) = \eta \s \eta$ since
$\eta^2= e$. By Lemma \ref{b5Feb9}, we have the group homomorphism
(recall that $L_n= \mS_n/ \ga_n$): 
\begin{equation}\label{xiaut}
\xi : G_n \ra \Aut_{K-{\rm alg}}(L_n), \;\; \s\mapsto
(\overline{\s} : a+\ga_n \mapsto \s (a) +\ga_n).
\end{equation}
The homomorphisms $\heta $ and $\xi$ will be used often in the
study of the group $G_n$. We can easily find the group
$\Aut_{K-{\rm alg}}(L_n)$ of algebra automorphisms of the Laurent
polynomial algebra $L_n$. We are interested in finding the image
and the kernel of the homomorphism $\xi$ (Corollary
\ref{a10Feb9}). We will see that the image of $\xi$ is small (and
the homomorphism $\xi$ is far from being surjective).

Next, several important subgroups of $G_n$ are introduced, they
are building blocks of the group $G_n$ (Theorem \ref{5Feb9}).

$\noindent $

{\bf The group $\Inn (\mS_n)$ of inner automorphism of $\mS_n$}.
Let $\mS_n^*$ be the group of units of the algebra $\mS_n$. The
centre $Z(\mS_n)$ of the algebra $\mS_n$ is $K$, \cite{shrekalg}.
For each element $u\in \mS_n^*$, let $\o_u:\mS_n\ra \mS_n$, $
a\mapsto uau^{-1}$, be the inner automorphism associated with the
element $u$. Then the group of inner automorphisms of the algebra
$\mS_n$,
$$\Inn (\mS_n)= \{ \o_u\, | \, u\in \mS_n^*\}\simeq \mS_n^*/
K^*,$$ is a normal subgroup of $G_n$. It follows from
\begin{equation}\label{etaH1}
\heta (\o_u) = \o_{\eta (u)^{-1}}, \;\; u\in \mS_n^*,
\end{equation}
that $\heta (\Inn (\mS_n)) = \Inn (\mS_n)$. The factor algebra
$\mS_n/\ga_n= L_n$ is commutative, and so $\xi (\Inn (\mS_n))= \{
e\}$.

$\noindent $

{\bf The torus $\mT^n$}. The $n$-dimensional torus $\mT^n:= \{
t_\l \, | \, \l = (\l_1, \ldots , \l_n) \in K^{*n}\}$ is a
subgroup of $G_n$ where
$$t_\l (x_i) = \l_ix_i, \;\; t_\l (y_i)= \l_i^{-1}y_i, \;\; i=1,
\ldots , n.$$ The torus $\mT^n:= \{ t_\l \, | \, \l \in K^{*n}\}$
is also a subgroup of the group $\Aut_{K-{\rm alg}}(L_n)$ where
$$t_\l
(x_i) = \l_ix_i, \;\;  i=1, \ldots , n.$$ Then $\heta (\mT^n) =
\mT^n$ and $\heta (t_\l ) = t_\l^{-1} = t_{\l^{-1}}$ where
$\l^{-1} := (\l_1^{-1}, \ldots , \l_n^{-1})$; $\xi (\mT^n) =
\mT^n$ and $\xi (t_\l ) = t_\l$. So, the maps $\heta : \mT^n\ra
\mT^n$ and $\xi : \mT^n\ra \mT^n$ are group isomorphisms. Note
that 
\begin{equation}\label{tlEab}
t_\l (E_{\alpha \beta})=\l^{\alpha - \beta}E_{\alpha, \beta }
\end{equation}
where $\l^{\alpha - \beta}:= \prod_{i=1}^n\l_i^{\alpha_i -
\beta_i}$.

$\noindent $

{\bf The symmetric group $S_n$}. The group $G_n$ contains the
symmetric group $S_n$ where each elements $\tau $ of $S_n$ is
identified with the automorphism of the algebra $\mS_n$ given by
the rule:
$$ \tau (x_i) = x_{\tau (i)}, \;\; \tau (y_i) = y_{\tau (i)}, \;\;
i=1, \ldots , n.$$ The group $S_n$ is also a subgroup of the group
$\Aut_{K-{\rm alg}}(L_n)$ where
$$ \tau (x_i) = x_{\tau (i)},\;\;
i=1, \ldots , n.$$ Clearly, $\heta (S_n) = S_n$ and $\heta (\tau )
= \tau$ for all $\tau \in S_n$;  $\xi (S_n) = S_n$ and $\xi (\tau
) = \tau$ for all $\tau \in S_n$. Note that 
\begin{equation}\label{tlEab1}
\tau (E_{\alpha \beta })= E_{\tau (\alpha ) \tau (\beta )}
\end{equation}
where $\tau (\alpha ):= (\alpha_{\tau^{-1}(1)}, \ldots ,
\alpha_{\tau^{-1}(n)})$.

$\noindent $

The groups $G_n$ and $\Aut_{K-{\rm alg}}(L_n)$ contain the
semi-direct product $S_n\ltimes \mT^n$ since $\mT^n\cap S_n = \{
e\}$ and 
\begin{equation}\label{tt1}
 \tau t_\l \tau^{-1}= t_{\tau (\l )}\;\; {\rm where}\;\; \tau (\l
) := (\l_{\tau^{-1}(1)},  \ldots , \l_{\tau^{-1}(n)}),
\end{equation}
 for
all $\tau \in S_n$ and $t_\l \in \mT^n$. Clearly, the maps
\begin{eqnarray*}
\heta : S_n\ltimes \mT^n\ra S_n\ltimes \mT^n, &  \tau t_\l\mapsto \tau t_\l^{-1},  \\
 \xi : S_n\ltimes \mT^n\ra S_n\ltimes \mT^n, &  \tau t_\l\mapsto \tau t_\l,
\end{eqnarray*}
are group isomorphisms.

\begin{lemma}\label{e5Feb9}
$S_n\ltimes \mT^n \ltimes \Inn (\mS_n)\subseteq G_n $.
\end{lemma}

{\it Proof}. We know already that $\Inn (\mS_n)$ and $ S_n\ltimes
\mT^n$ are subgroups of $G_n$. Since $\Inn (\mS_n) \subseteq \ker
(\xi) $ and $\xi : S_n\ltimes \mT^n \simeq S_n\ltimes \mT^n$, we
see that $\Inn (\mS_n) \cap (S_n\ltimes \mT^n)=\{ e\}$, and  the
result follows. $\Box $

$\noindent $

Let $r$ be an element of a ring $R$. The
 element $r$ is called {\em regular} if $ \lann_R(r)=0$ and
 $\rann_r(r)=0$ where $ \lann_R(r):= \{ s\in R\, | \, sr=0\}$ is
 the {\em left annihilator} of $r$ and $ \rann_R(r):= \{ s\in R\, | \, rs=0\}$ is
 the {\em right annihilator} of $r$.

The next lemma shows that the elements $x$ and $y$ of the algebra
$\mS_1$ are not regular.
\begin{lemma}\cite{shrekalg}\label{a7Dec8}
\begin{enumerate}
\item $\lann_{\mS_1}(x) = \mS_1E_{00} = \bigoplus_{i\geq 0}
KE_{i,0}=\bigoplus_{i\geq 0} Kx^i(1-xy)$ and $\rann_{\mS_1}(x)=0$.
\item $\rann_{\mS_1}(y) = E_{00}\mS_1 =  \bigoplus_{i\geq 0}
KE_{0, i}=\bigoplus_{i\geq 0} K (1-xy)y^i$ and
$\lann_{\mS_1}(y)=0$.
\end{enumerate}
\end{lemma}

It follows from Lemma \ref{a7Dec8} that, for each $i=1, \ldots ,
n$, 
\begin{equation}\label{lanxin}
\lann_{\mS_n}(x_i) = \mS_{n-1}\t \lann_{\mS_1(i)}(x_i)=
 \bigoplus_{j\geq 0} \mS_{n-1}E_{j,0}(i)=\bigoplus_{j\geq 0}
 \mS_{n-1}x_i^jE_{00}(i),
\end{equation}
\begin{equation}\label{lanyin}
\rann_{\mS_n}(y_i) = \mS_{n-1}\t \rann_{\mS_1(i)}(y_i)=
 \bigoplus_{j\geq 0} E_{0,j}(i)\mS_{n-1}=\bigoplus_{j\geq
 0}E_{00}(i)y_i^j\mS_{n-1},
\end{equation}
where $\mS_{n-1}$ stand for $\bigotimes_{k\neq i}\mS_1(k)$.

For an algebra $A$ and a subset $S\subseteq A$, $\Cen_A(S):=\{
a\in A\, | \, as=sa$ for all $s\in S\}$ is the {\em centralizer}
of the set $S$ in $A$. It is a subalgebra of $A$. It was proved in
\cite{shrekalg} that 
\begin{equation}\label{Cenxy}
\Cen_{\mS_n}(x_1, \ldots , x_n)=K[x_1, \ldots , x_n], \;\;
\Cen_{\mS_n}(y_1, \ldots , y_n)=K[y_1, \ldots , y_n].
\end{equation}

Let $E_n:= \End_{K-{\rm alg}}(\mS_n)$ be the monoid of all the
$K$-algebra endomorphisms of $\mS_n$. The group of units of this
monoid is $G_n$. The automorphism $\heta \in \Aut (G_n)$ can be
extended to an automorphism $\heta \in \Aut (E_n)$ of the monoid
$E_n$: 
\begin{equation}\label{1etah}
\heta : E_n\ra E_n, \;\; \s \mapsto \eta \s \eta^{-1}.
\end{equation}
The next (curious) result is instrumental in finding the group of
automorphisms of the algebra $\mS_n$.
\begin{theorem}\label{6Feb9}
Let $\s , \tau \in G_n$. Then the following statements are
equivalent.
\begin{enumerate}
 \item $\s = \tau$. \item $\s (x_1) = \tau (x_1), \ldots , \s (x_n) = \tau
 (x_n)$.
 \item$\s (y_1) = \tau (y_1), \ldots , \s (y_n) = \tau (y_n)$.
\end{enumerate}
\end{theorem}

{\it Proof}. Without loss of generality we may assume that $\tau
=e$, the identity automorphism. Consider the following two
subgroup of $G_n$, the stabilizers of the sets $\{ x_1, \ldots ,
x_n\}$ and $\{ y_1, \ldots , y_n\}$:
\begin{eqnarray*}
\St  (x_1, \ldots , x_n):= & \{ g\in G_n\, | \, g(x_1) = x_1,
\ldots , g(x_n)=x_n\}, \\
 \St  (y_1, \ldots , y_n):= & \{ g\in G_n\,
| \, g(y_1) = y_1, \ldots , g(y_n)=y_n\}.
\end{eqnarray*}
Then $$\heta (\St  (x_1, \ldots , x_n))=\St  (y_1, \ldots , y_n),
\;\; \heta (\St  (y_1, \ldots , y_n))= \St  (x_1, \ldots , x_n).$$
Therefore, the theorem (where $\tau = e$) is equivalent to the
single statement that $\St  (x_1, \ldots , x_n)=\{ e\}$. So, let
$\s \in \St  (x_1, \ldots , x_n)$. We have to show that $\s = e$.
For each  $i=1, \ldots , n$, $1=\s (y_ix_i) = \s (y_i) x_i$ and
$1=y_ix_i$. By taking the difference of these equalities we see
that $a_i:= \s (y_i) -y_i\in \lann_{\mS_n}(x_i)$. By
(\ref{lanxin}), $a_i=\sum_{j\geq 0} \l_{ij}E_{j0}(i)$ for some
elements $\l_{ij}\in \bigotimes_{k\neq i}\mS_1(i)$, and so $$ \s
(y_i) = y_i+\sum_{j\geq 0}\l_{ij} E_{j0}(i).$$ The element $\s
(y_i)$ commutes with the elements $\s (x_k) = x_k$, $k\neq i$,
hence all $\l_{ij}\in K[x_1, \ldots , \hx_i, \ldots , x_n]$, by
(\ref{Cenxy}). Since $E_{j0}(i) = x_i^j E_{00}(i)$, we can write
$$ \s (y_i) =y_i+ p_iE_{00}(i) \;\; {\rm for \; some}\;\; p_i\in
P_n.$$ We have to show that all $p_i=0$. Suppose that this is not
the case. Then $p_i\neq 0$ for some $i$. We seek a contradiction.
Note that $\s^{-1}\in \St (x_1, \ldots , x_n)$, and so $\s (y_i)
=y_i+ q_iE_{00}(i)$ for some  $q_i\in P_n$. Recall that $E_{00}(i)
= 1-x_iy_i$. Then $\s^{-1} (E_{00}(i))= 1-x_i(y_i+q_iE_{00}(i))=
(1-x_iq_i)E_{00}(i)$, and
$$ y_i=\s^{-1}\s (y_i) = \s^{-1} (y_i+p_iE_{00}(i))=
y_i+(q_i+p_i(1-x_iq_i))E_{00}(i), $$ and so $q_i+p_i= x_ip_iq_i$
since the map $P_n\ra P_nE_{00}$, $p\mapsto pE_{00}$, is an
isomorphism of $P_n$-modules as it follows from (\ref{xyEij}).
This is impossible by comparing the degrees of the polynomials on
both sides of the equality. $\Box $

$\noindent $

Theorem \ref{6Feb9} states that each automorphism of the
non-commutative, finitely generated, non-Noetherian algebra
$\mS_n$ is uniquely determined by its action on its commutative,
finitely generated subalgebra $P_n$. A similar result is true for
the ring $\CD (P_n)$ of differential operators on the polynomial
algebra $P_n$ over a field of {\em prime} characteristic. The
algebra $\CD (P_n)$ is a non-commutative, {\em not finitely
generated}, non-Noetherian algebra.

\begin{theorem}
{\rm \cite{autlaur}   (Rigidity of the group $\Aut_{K-{\rm
alg}}(\CD (P_n)$))} Let $K$ be a field of prime characteristic,
and   $\s , \tau \in \Aut_{K-{\rm alg}}(\CD (P_n)$. Then $\s =
\tau$ iff  $\s (x_1) = \tau (x_1), \ldots , \s (x_n) = \tau
 (x_n)$.
\end{theorem}
The above theorem doest not
hold in characteristic zero and doest not hold in prime
characteristic for the ring of differential operators on a Laurent
polynomial algebra \cite{autlaur}.


\section{The groups $\Aut_{K-{\rm alg}}(\mS_1)$ and
$\mS_1^*$}\label{TG1S}

In this section, the groups $\Aut_{K-{\rm alg}}(\mS_1)$ and
$\mS_1^*$ are found (Theorems \ref{A5Feb9} and \ref{a13Dec8}). The
case $n=1$ is rather special and much more simpler than the
general case. It is a sort of a degeneration of the general case.
Briefly, the key idea in finding the group of automorphisms of the
algebra $\mS_1$ is to use Theorem \ref{6Feb9} and some properties
of the index of linear maps in the vector space $P_1=K[x]$. We
start this section with a sketch of the proof of Theorem
\ref{A5Feb9}.  Then we prove necessary results about the index of
certain elements of the algebra $\mS_1$, and  using them we find
the group $\mS_1^*$ of units of the algebra $\mS_1$ and the group
$\Inn (\mS_1)$ of inner automorphisms of $\mS_1$;  and finally we
give the proof of Theorem \ref{A5Feb9}. The proof is constructive
in the sense that for each automorphism $\s$ of the algebra
$\mS_1$ it gives explicitly the presentation $\s = t_\l \o_\v$ of
 $\s$ as the product of an inner automorphism $\o_\v$ and and
 element $t_\l$ of the torus $\mT^1$ (Corollary \ref{d8Feb9}).

\begin{theorem}\label{A5Feb9}
$\Aut_{K-{\rm alg}}(\mS_1 )=\mT^1\ltimes \Inn (\mS_1)\simeq
\mT^1\ltimes \GL_\infty (K)$.
\end{theorem}

{\it Sketch of the Proof}. {\it Step 1}. Let $\s \in G_1$. By
Lemma \ref{e5Feb9}, we have to show that $\s \in \mT^1\ltimes \Inn
(\mS_1)$. Using some properties of the index of linear maps from
$\End_K(P_1)$ that have finite dimensional kernel and cokernel, we
show that
\begin{eqnarray*}
 \s (x)&\equiv & \l x \mod F,  \\
 \s (y)&\equiv & \l^{-1}y \mod F,
\end{eqnarray*}
for some element $\l \in K^*$.

{\it Step 2}. Changing $\s $ for $t_{\l^{-1}}\s$ we may assume
that $\l =1$.

{\it Step 3}. Changing $\s$  for $\o_\v\s$ for a suitable choice
of a unit $\v$ of the algebra $\mS_1$ we may assume that $\s (y) =
y$.

{\it Step 4}. Then, by Theorem \ref{6Feb9},  $\s =e$.  $\Box $

$\noindent $

{\it Remark}. The multiplication in the skew product $\mT^1\ltimes
\GL_\infty (K)$ is given by the rule: 
\begin{equation}\label{tvtp}
\v t_\l  \cdot \psi t_\mu  = \v t_\l (\psi ) t_{\l \mu }
\end{equation}
where $t_\l , t_\mu \in \mT^1$; $\v , \psi \in \GL_\infty (K)$;
and
 $t_\l (\psi ) $ is defined in (\ref{tlEab}).

$\noindent $

{\bf The index $\ind$ of linear maps and its properties}. Let $\CC
= \CC (K)$ be the family of all $K$-linear maps with finite
dimensional kernel and cokernel.

$\noindent $

{\it Definition}. For a linear map $\v \in \CC$, the integer
$$ \ind (\v ) := \dim \, \ker (\v ) - \dim \, \coker (\v )$$
is called the {\em index} of the map $\v$.

$\noindent $

{\it Example}. Note that $\mS_1\subset \End_K(P_1)$. One can
easily prove that 
\begin{equation}\label{indxy}
\ind (x^i)= -i\;\; {\rm and }\;\; \ind (y^i)= i, \;\; i\geq 1.
\end{equation}

Lemma \ref{b8Feb9} shows that $\CC$ is a multiplicative semigroup
with zero element (if the composition of two elements of $C$ is
undefined we set their product to be zero).

\begin{lemma}\label{b8Feb9}
Let $\psi : M\ra N $ and $\v : N\ra L$ be $K$-linear maps. If two
of the following three maps: $\psi$, $\v$,  and $\v \psi$, belong
to the set $\CC$ then so does the third; and in this case, $$ \ind
(\v\psi ) = \ind (\v ) + \ind (\psi ).$$
\end{lemma}

{\it Proof}.  For an arbitrary $K$-linear map $f: V\ra U$, we use
the following notation: ${}_fV:= \ker (f)$ and $U_f:=\coker (f)$.
The result follows from the long exact sequence of $K$-linear maps
(where all the maps are natural): 
\begin{equation}\label{MNLex}
0\ra{}_\psi M \ra {}_{\v\psi}M\stackrel{\psi}{\ra} {}_\v N\ra
N_\psi \stackrel{\v}{\ra} L_{\v \psi} \ra L_\v \ra 0.
\end{equation}
In particular, taking the Euler characteristic of the long exact
sequence (\ref{MNLex}) gives the identity $\ind (\psi ) -\ind
(\v\psi ) + \ind (\v ) =0$. $\Box $
\begin{lemma}\label{a8Feb9}
Let
$$
\xymatrix{0\ar[r] & V_1\ar[r]\ar[d]^{\v_1}  & V_2 \ar[r]\ar[d]^{\v_2} & V_3 \ar[r]\ar[d]^{\v_3 } & 0 \\
0\ar[r] & U_1\ar[r]  & U_2\ar[r] & U_3 \ar[r] & 0 }
$$
be a commutative diagram of $K$-linear maps with exact rows.
Suppose that $\v_1, \v_2, \v_2\in \CC$. Then
$$ \ind (\v_2) = \ind (\v_1)+\ind (\v_3).$$
\end{lemma}

{\it Proof}. The Snake Lemma yields the long exact sequence:
$$ 0\ra \ker (\v_3) \ra \ker (\v_2) \ra \ker (\v_1)\ra  \coker (\v_3) \ra \coker (\v_2) \ra \coker
(\v_1)\ra 0 $$ Taking its Euler characteristic gives $ \ind
(\v_1)-   \ind (\v_2)+\ind (\v_3)=0$. $\Box $

$\noindent $

Each nonzero element $u$ of the Laurent polynomial algebra $L_1=
K[x,x^{-1}]$ is a unique sum $u=\l_sx^s+\l_{s+1}x^{s+1}+\cdots
+\l_dx^d$ where all $\l_i\in K$, $\l_d\neq 0$, and $\l_dx^d$ is
the {\em leading term} of the element $u$. The integer
$\deg_x(u)=d$ is called the {\em degree} of the element $u$. It is
an extension to $L_1$ of the usual degree of polynomials in
$K[x]$. The next lemma explains how to compute the indices of the
elements $\mS_1\backslash F$ using the degree function $\deg_x$
and shows that the index is a $G_1$-invariant concept. Note that
$F\cap \CC = \emptyset$.
\begin{lemma}\label{c8Feb9}
\begin{enumerate}
\item $\mS_1\backslash F\subseteq \CC $ (recall that $\mS_1\subset
\End_K(P_1)$) and, for each element $a\in \mS_1\backslash F$,
$$ \ind (a) = -\deg_x(\oa )$$
where $\oa = a+F\in \mS_1/F= L_1$. \item $\ind (\s (a)) = \ind
(a)$ for all $\s \in G_1$ and $a\in \mS_1\backslash F$.
\end{enumerate}
\end{lemma}

{\it Proof}. 1. Let $a\in \mS_1\backslash F$ and $d:= \deg_x(\oa
)$. The element of the algebra $\mS_1$,
$$b:= \begin{cases}
y^da& \text{if }\; d\geq 0,\\
ax^{-d}& \text{if }\; d<0, \\
\end{cases}
$$
does not belong to the ideal $F$ (since $\ob = x^{-d}\oa \neq 0$),
and $\deg_x(\ob )=0$. By Lemma \ref{b8Feb9} and (\ref{indxy}), it
suffices to prove that $\ind (b)=0$ since then
$$ 0=\ind (b) =
d+\ind (a),$$   that is $\ind (a) = -\deg_x(\oa )$. The element
$b$ can be written as a sum $b=\l +\sum_{i\geq 1}\l_iy^i+f$ for
some elements $\l\in K^*$, $\l_i\in K$, and $f\in F$. Fix a
natural number $m$ such that $f\in M_{m+1}(K)$ (recall that
$F=\cup_{i\geq 1}M_i(K)$. Abusing notation, let $K[b]$ be the
$K$-subalgebra of $\End_K(P_1)$ generated by the element $b$. Then
$V:= \bigoplus_{i=0}^m Kx^i$ is a $K[b]$-submodule of $P_1$, and
$U:= P_1/V$ is the factor module. Let $b_1$ and $b_2$ be the
linear maps that are determined by the action of the element $b$
on the vector spaces $V$ and $U$ respectively. Then $\ind (b_1)=0$
since $\dim (V)<\infty$; and $\ind (b_2)=0$ since $b_2=\l
+\sum_{i\geq 1}\l_iy^i$ is a bijection. Applying Lemma
\ref{a8Feb9} to the commutative diagram
$$
\xymatrix{0\ar[r] & V\ar[r]\ar[d]^{b_1}  & P_1 \ar[r]\ar[d]^{b} & U \ar[r]\ar[d]^{b_2 } & 0 \\
0\ar[r] & V\ar[r]  & P_1\ar[r] & U \ar[r] & 0 }
$$
we have the result: $\ind (b) = \ind (b_1)+\ind (b_2)=0$.

 2. By Theorem \ref{a5Feb9}, $\ind (\s (a))= \ind (\v a \v^{-1}) =
 \ind (a)$ where $\s = \s_\v$. $\Box $

$\noindent $

{\bf The group of units $(1+F)^*$ and $\mS_1^*$}. Recall that the
 algebra (without 1) $F=\bigoplus_{i,j\in \N} KE_{ij}$ is the union
$M_\infty (K) := \bigcup_{d\geq 1}M_d(K)= \varinjlim M_d(K)$ of
the matrix algebras $M_d(K):= \bigoplus_{1\leq i,j\leq
d-1}KE_{ij}$, i.e. $F= M_\infty (K)$.

For each $d\geq 1$, consider the (usual) determinant $\det_d=\det
: 1+M_d(K)\ra K$, $u\mapsto \det (u)$. These determinants
determine the (global) {\em determinant}, 
\begin{equation}\label{gldet}
\det : 1+M_\infty (K)= 1+F\ra K, \;\; u\mapsto \det (u),
\end{equation}
where $\det (u)$ is the common value of all determinants
$\det_d(u)$, $d\gg 1$. The (global) determinant has usual
properties of the determinant. In particular, for all $u,v\in
1+M_\infty (K)$,
$$\det (uv) = \det (u) \cdot \det (v).$$ It follows
from this equality and the  Cramer's formula for the inverse of a
matrix that the group $\GL_\infty (K):= (1+M_\infty (K))^*$ of
units of the monoid $1+M_\infty (K)$ is equal to 
\begin{equation}\label{GLiK}
\GL_\infty (K) = \{ u\in 1+M_\infty (K) \, | \, \det (u) \neq 0\}.
\end{equation}
Therefore, 
\begin{equation}\label{1GLiK}
(1+F)^* = \{ u\in 1+F \, | \, \det (u) \neq 0\}=\GL_\infty (K).
\end{equation}
 The kernel
$$\SL_\infty (K):= \{ u\in
\GL_\infty (K)\, | \, \det (u) =1\}$$
 of the group epimorphism
$\det : \GL_\infty (K)\ra K^*$ is a {\em normal} subgroup of
$\GL_\infty (K)$.

Let $V$ be an infinite dimensional vector space that has countable
basis. A sequence $\CV$ of finite dimensional vector spaces in
$V$, $V_0\subseteq V_1\subseteq \cdots \subseteq V_i\subseteq
\cdots$, such that $V=\bigcup_{i\geq 0}V_i$ is called  a {\em
finite dimensional vector space filtration} on $V$. The next
result reveals an invariant nature of the (global) determinant.

\begin{lemma}\label{a6Mar9}
Let $\CV = \{ V_i\}_{i\geq 0}$ be a finite dimensional vector
space filtration on the polynomial algebra $P_1=K[x]$ and $a\in
\mM_1:=1+M_\infty (K)$. Then $a(V_i) \subseteq V_i$ for all $i\gg
0$, and $\det (a|_{V_i})=\det (a|_{V_j})$ for all $i,j\gg 0$ where
$\det (a|_{V_i})$ is the determinant of the linear map
$a|_{V_i}:V_i\ra V_i$. Moreover, this common value of the
determinants, $\det (a) =  det_\CV (a)$, does not depend on the
filtration $\CV$ and, therefore, coincides with the determinant
defined in (\ref{gldet}).
\end{lemma}

{\it Proof}. Let $a\in \mM_1$. Then
$a=1+\sum_{i,j=0}^d\l_{ij}E_{ij}$ for some scalars $\l_{ij}\in K$
and $d\in \N$. Note that the global determinant $\det (a)$, as
defined in (\ref{gldet}), is equal to the usual determinant $\det
(a|_{P_{1, \leq i}})$ for all $i\geq d$, where $\{ P_{1, \leq i}:=
\sum_{j=0}^iKx^i\}_{i\in \N}$ is the {\em degree filtration} on
$P_1$. Then $\im (a-1)\subseteq P_{1, \leq d}\subseteq V_e$ for
some $e\in \N$. Since $a=1+(a-1)$, we have $a(V_i) \subseteq V_i$
and $\det (a|_{V_i})= \det (a|_{V_e})$ for all $i\geq e$. Note
that this is true for an arbitrary  finite dimensional vector
space filtration $\CV$. Consider the following finite dimensional
vector space filtration $$\CV':= \{ V_i':= P_{1, \leq d}, \; i=0,
\ldots , e-1; \; V_j':=V_j, \; j\geq e\}.$$ Then $$ \det (a) =
\det(a|_{P_{1, \leq d}}) = \det (a|_{V_{e-1}'}) = \det (a|_{V_j'})
\det ( a|_{V_j}), \;\; j\geq e.$$
 This completes
the proof of the lemma. $\Box $

$\noindent $

The {\em centre} of a group $G$ is denoted $Z(G)$.

\begin{theorem}\label{a13Dec8}
\begin{enumerate}
\item $\mS_1^*=K^*(1+F)^*\simeq K^*\times \GL_\infty (K)$.\item
$Z(\mS_1^*) = K^*$ and $Z((1+F)^*)=\{ 1 \}$.  \item $\Inn
(\mS_1)\simeq \GL_\infty (K)$, $\o_u\lra u$.
\end{enumerate}
\end{theorem}

{\it Proof}. 1. Note that $\mS_1^*\supseteq K(1+F)^*\simeq
K^*\times (1+F)^*\simeq K^*\times \GL_\infty (K)$ since $K^*\cap
(1+F)^*=\{ 1\}$. It remains to prove the reverse inclusion. If an
element $u$ is a unit of the algebra $\mS_1$ then the element $\bu
=u+F$ is a unit of the factor algebra $L_1=\mS_1/F$, and so $\bu =
\l x^i$ for some $\l\in K^*$ and $i\in \Z$. Therefore, either $u =
\l x^i+f$ or $u=\l y^i+f$ for some $\l \in K^*$ and $i\in \N$. The
element $u\in \mS_1\backslash F$ is a unit, hence $u\in
\End_K(P_1)$ is an invertible linear map (recall that
$\mS_1\subset \End_K(P_1)$), and so $\ind (u) =0$. By Lemma
\ref{c8Feb9}.(1) and (\ref{indxy}), $i=0$, and so $u\in
K^*(1+F)^*$.

2. Note that $Z(\mS_1^*)= K^* Z((1+F)^*)$. It suffices to show
that $Z((1+F)^*)=\{ 1\}$. Let $z= 1+\sum \l_{ij}E_{ij}\in
Z((1+F)^*)$ where $\l_{ij}\in K$. For all $k\neq l$, $1+E_{kl}\in
(1+F)^*$ since $\det (1+E_{kl})=1$. Now, $z(1+E_{kl})=(1+E_{kl})z$
for all $k\neq l$ iff $\sum_i \l_{ik}E_{il}=\sum_j\l_{lj}E_{kj}$
for all $k\neq l$ iff all $\l_{ij}=0$ iff $z=1$.

3. $\Inn (\mS_1) \simeq \mS_1^*/Z(\mS_1^*)\simeq (K^*\times
\GL_\infty (K))/K^*\simeq \GL_\infty (K)$.  $\Box $

$\noindent $

{\bf Proof of Theorem \ref{A5Feb9}}. By Theorem \ref{a13Dec8}.(3),
$\mT^1\ltimes \Inn (\mS_1) = \mT^1\ltimes \GL_\infty (K)$.

Let $\s \in G_1$. By Lemma \ref{e5Feb9}, in order to finish the
proof of the theorem  we have to show that $\s \in \mT^1\ltimes
\Inn (\mS_1)$. By Lemma \ref{b5Feb9}, $\s (F) = F$, and so the map
$$ \overline{\s}: L_1=\mS_1/F\ra  L_1=\mS_1/F, \;\; \oa = a+F\mapsto \s (a)+F,$$
is an isomorphism of the Laurent polynomial algebra $L_1=
K[x,x^{-1}]$. Therefore, either $ \overline{\s}(y) = \l x^{-1}$
or, otherwise, $ \overline{\s} (y) = \l x$ for some scalar $\l \in
K^*$. Equivalently, either $\s (y) = \l y +f$ or $\s (y) = \l x
+f$ for some element $f\in F$. By Lemma \ref{c8Feb9}, the second
case is impossible since, by (\ref{indxy}),
$$ 1=\ind (y) = \ind (\s (y))= \ind (\l x+f) = -\deg_x(\l x) =
-1.$$ Therefore, $\s (y) = \l y +f$. Then, $t_\l\s (y) = y+g$
where $g:= t_\l (f) \in F$ since $t_\l (F) = F$ (Lemma
\ref{e5Feb9}). Fix a natural number $m$ such that $g\in
M_{m+1}(K)$. Then the finite dimensional vector spaces
$$ V:= \bigoplus_{i=0}^m Kx^i\subset V':=
\bigoplus_{i=0}^{m+1}Kx^i$$ are $y'$-invariant where $y':=t_\l \s
(y)= y+g$. Note that $y'*x^{m+1} = y*x^{m+1}= x^m$ since
$g*x^{m+1}=0$.
 Note that $P_1=\bigcup_{i\geq 1}\ker (y^i)$ and $\dim \, \ker_{P_1} (y)
 =1$. Since the $\mS_1$-modules $P_1$ and ${}^{t_\l \s}P_1$ are
 isomorphic, $P_1=\bigcup_{i\geq 1}\ker (y'^i)$ and $\dim \, \ker_{P_1} (y')
 =1$. This implies that the elements $x_0', x_1', \ldots , x_m',
 x^{m+1}$ are a $K$-basis for the vector space $V'$ where
 $$ x_i':= y'^{m+1-i}*x^{m+1}, \;\; i=0,1, \ldots , m; $$
 and the elements  $x_0', x_1', \ldots , x_m'$ are a $K$-basis for the vector space
 $V$. Then the elements
 $$ x_0', x_1', \ldots , x_m', x^{m+1}, x^{m+2},  \ldots  $$
are a $K$-basis for the vector space $P_1$. The $K$-linear map
\begin{equation}\label{chm}
\v : P_1\ra P_1, \;\; x^i\mapsto x_i'\; (i=0,1,\ldots , m), \;
x^j\mapsto x^j\; (j>m),
\end{equation}
belongs to the group $(1+F)^*= \GL_\infty (K)$ and satisfies the
property that
$$ y'\v = \v y, $$
the equality is in $\End_K(P_1)$. This equality can be rewritten
as follows:
$$ \o_{\v^{-1}}t_\l \s (y) = y \;\; {\rm where}\;\;
\o_{\v^{-1}}\in \Inn (\mS_1).$$ By Theorem \ref{6Feb9}, $\s =
t_{\l^{-1}}\o_\v \in \mT^1\ltimes \Inn (\mS_1)$, as required.
$\Box $

\begin{corollary}\label{d8Feb9}
Each automorphism $\s$ of the algebra $\mS_1$ is a unique product
$\s = t_{\l^{-1}}\o_\v$ where $\s (y) \equiv \l y \mod F$ and
$\v\in (1+F)^*= \GL_\infty (K)$ is defined as in (\ref{chm}).
\end{corollary}

{\it Proof}. The result was established in the proof of Theorem
\ref{A5Feb9} apart from the uniqueness of $\v$ which follows from
the fact that the centre of the group $(1+F)^* = \GL_\infty (K)$
is $\{ 1\}$ (Theorem \ref{a13Dec8}.(3)). $\Box $

\begin{proposition}\label{b7Dec8}
Each algebra endomorphism of $\mS_1$ is either a monomorphism or,
otherwise, its image is  a commutative finite dimensional algebra.
In the second case, all positive integers occur as the dimension
of the image.
\end{proposition}

{\it Proof}. Recall that $F$ is the  smallest nonzero ideal of the
algebra $\mS_1$, and $\mS_1/F\simeq K[x,x^{-1}]$ (see
(\ref{mS1d1})). If an algebra endomorphism $\s$ of $\mS_1$ is not
a monomorphism then $F\subseteq \ker (\s )$, and so $\s (x) \in
\mS_1^* = K^* (1+F)^*$ (Theorem \ref{a13Dec8}.(1)) since the
equalities $yx=1$ and $xy=1-E_{00}$ imply the equalities $\s (y)
\s (x) =1$ and $\s (x) \s (y) =1$; and $\im
 ( \s ) =K\langle \s (x),\s (x^{-1})\rangle $. Therefore, the image
of $\s$ is a commutative  finite dimensional algebra since the
algebra  $K\langle \s (x),\s (x^{-1})\rangle$ can be seen as a
subalgebra of the matrix algebra $ M_d(K)$ for some $d$. The image
of the endomorphism $\mS_1\ra \mS_1$, $x\mapsto 1$, $y\mapsto 1$,
is $K$, hence one-dimensional. For each natural number $n\geq 2$,
the image of the endomorphism
$$\s_n : \mS_1\ra \mS_1, \;\; x\mapsto 1+\gn, \;\; y\mapsto (1+\gn
)^{-1}, \;\; (\gn := \sum_{i=0}^{n-2}E_{i,i+1})$$ has dimension
$n$ since the set $1, \gn, \gn^2, \ldots , \gn^{n-1}$ is a
$K$-basis of the image of $\s_n$. $\Box $


\section{The group of automorphisms of the algebra $\mS_n$}\label{GAS}

In this section, the group $G_n$ is found (Theorem \ref{5Feb9}).
It is shown that the groups $G_n$ and $\Inn (\mS_n)$ have trivial
centre (Corollary \ref{b24Feb9}).

By the very definition, the subset $\stCH$ of $\St_{G_n}(\CH_1)$
(see (\ref{stCH})) is a subgroup of $\St_{G_n}(\CH_1)$.

\begin{theorem}\label{5Feb9}
 $G_n=S_n\ltimes \mT^n\ltimes \Inn (\mS_n)$.
\end{theorem}

{\it Proof}. The group $G_n$ acts in the obvious way, $(\s , \gp_i
) \mapsto \s (\gp_i)$, on the set $\CH_1:= \{ \gp_1, \ldots ,
\gp_n\}$ of all the height 1 prime ideals of the algebra $\mS_n$.
In particular, the symmetric group $S_n$, which is a subgroup of
$G_n$, permutes the ideals $\gp_1, \ldots , \gp_n$, i.e. $\tau
(\gp_i) = \gp_{\tau (i)}$ for $\tau \in S_n$. The stabilizer
$$\St_{G_n}(\CH_1 )= \{ \s \in G_n\, | \, \s (\gp_1) = \gp_1,
\ldots , \s (\gp_n) = \gp_n\}$$ is a {\em normal} subgroup of
$G_n$ such that $G_n = S_n \St_{G_n}(\CH_1)$ and $S_n\cap
\St_{G_n}(\CH_1)=\{ e\}$, and so 
\begin{equation}\label{GSSn}
G_n=S_n \ltimes \St_{G_n}(\CH_1).
\end{equation}
Clearly, $\mT^n\ltimes \Inn (\mS_n)\subseteq \St_{G_n}(\CH_1)$.
So, in order to finish the proof of the theorem we have to prove
that the inverse inclusion holds.

Let $\s \in \St_{G_n}(\CH_1)$. We have to show that $\s
\in\mT^n\ltimes \Inn (\mS_n)$. Since $\s (\gp_n) = \gp_n$, the
automorphism $\s$ induces the automorphism $$\s_n: \mS_n / \gp_n =
\mS_{n-1}\t L_1\ra \mS_n / \gp_n = \mS_{n-1}\t L_1, \;\;
a+\gp_n\mapsto \s (a) +\gp_n.$$ The restriction of the
automorphism $\s_n$ to the centre $Z(\mS_{n-1}\t L_1)= K[x_n,
x_n^{-1}]$ of the algebra $\mS_n / \gp_n$ yields its automorphism,
and so either $\s_n (x_n) = \l x_n$ or $\s_n(x_n )= \l x_n^{-1}$
for some scalar $\l \in K^*$. Therefore, there are two options:
\begin{eqnarray*}
 (i) & \s (x_n) =\l_n x_n+p_n, \;\; \s (y_n) = \l_n^{-1}y_n+q_n; \\
(ii)  &\s (x_n) =\l_n y_n+p_n, \;\; \s (y_n) = \l_n^{-1}x_n+q_n;
\end{eqnarray*}
for some $\l_n  \in K^*$ and $p_n, q_n\in \gp_n$. We aim to show
that the second case is impossible. This is true for $n=1$, by
Theorem \ref{A5Feb9}. So, let $n>1$.  Suppose that $\s (x_n) =\l_n
y_n+p_n$, wee seek a contradiction. By symmetry of the indices,
for each $i=1, \ldots , n$, there are two options:
\begin{eqnarray*}
 (i) & \s (x_i) =\l_i x_i+p_i, \;\; \s (y_i) = \l_i^{-1}y_i+q_i; \\
(ii)  &\s (x_i) =\l_i y_i+p_i, \;\; \s (y_i) = \l_i^{-1}x_i+q_i;
\end{eqnarray*}
for some $\l_i  \in K^*$ and $p_i, q_i\in \gp_n$. Since $\s
(\gp_1+\cdots +\gp_{n-1})=\gp_1+\cdots +\gp_{n-1}$ and $\mS_n/
(\gp_1+\cdots +\gp_{n-1})\simeq L_{n-1}\t \mS_1(n)$ where
$L_{n-1}= K[x_1^{\pm 1}, \ldots , x_n^{\pm 1}]$, the automorphism
$\s$ of the algebra $\mS_n$ induces an automorphism, say
$\overline{\s}$, of the algebra $L_{n-1}\t \mS_1(n)$ such that
either $\overline{\s}(x_i) = \l_ix_i$ or $\overline{\s}(x_i) =
\l_ix_i^{-1}$ for all $i=1, \ldots , n$. We see that
$\overline{\s}(L_{n-1}) = L_{n-1}$. Let $\g$ be the restriction of
the automorphism $\overline{\s}$ to the algebra $L_{n-1}$. Then
$\g \t \id_{\mS_1(n)}$ is the automorphism of the algebra $L_{n-1}
\t \mS_1(n)$. Then $\widetilde{\s}:= (\g \t
\id_{\mS_1(n)})^{-1}\overline{\s}$ is the $L_{n-1}$-algebra
automorphism of the algebra $L_{n-1}\t \mS_1(n)$ which can be
uniquely extended to a ${\rm Frac} (L_{n-1})$-automorphism of the
algebra ${\rm Frac} (L_{n-1})\t \mS_1(n)$ over the field of
fractions ${\rm Frac} (L_{n-1})= K(x_1, \ldots , x_{n-1})$ of the
algebra $L_{n-1}$. By Theorem \ref{A5Feb9} (or Corollary
\ref{d8Feb9}), we must have the case $(i)$ for $x_n$ and $y_n$.

By symmetry of the indices, it follows from the case (i) that
\begin{equation}\label{sxiyi}
\s (x_i) = \l_ix_i+p_i, \;\; \s (y_i) = \l_i^{-1}y_i+q_i, \;\;
i=1, \ldots ,n,
\end{equation}
for some scalars $\l_i\in K^*$ and some elements $p_i, q_i\in
\gp_i$.

Changing $\s$ for $t_{\l^{-1}}\s$, where $\l = (\l_1, \ldots ,
\l_n)$, we may assume that $\l_1=\cdots = \l_n=1$, that is, $\s
\in \stCH $. It follows that $G_n=S_n\mT^n\stCH $. To finish the
proof of the theorem it suffices to show that $\stCH\subseteq \Inn
(\mS_n)$ since then, by Lemma \ref{e5Feb9}, $G_n =S_n\ltimes \mT^n
\ltimes \Inn (\mS_n )$ and also 
\begin{equation}\label{stISn}
\stCH = \Inn (\mS_n).
\end{equation}
Let $\s \in \stCH$. Then $\s^{-1} \in \stCH$ since $\stCH$ is a
group. By Theorem \ref{a5Feb9}, $\s = \s_\v$ for some element $\v
\in \Aut_K(P_n)$ such that $\v \mS_n \v^{-1} = \mS_n$. For each
number $i=1, \ldots , n$, $p_i:= \s (x_i) - x_i\in \gp_i$ since
$\s \in \stCH$. By multiplying this equality on the left by
$\v^{-1}$, we obtain the equality $x_i\v^{-1} = \v^{-1} (x_i+p_i)$
for each $i=1, \ldots , n$. By Theorem \ref{10Feb9}, $\v^{-1} \in
\mS_n$. Repeating the same arguments for the automorphism
$\s^{-1}= \s_{\v^{-1}}\in \stCH$, we have $\v \in \mS_n$, that is
$\v \in \mS_n^*$, and so $\s $ is an inner automorphism of the
algebra $\mS_n$. $\Box $

\begin{corollary}\label{1a21Feb9}
The group ${\rm Out} (\mS_n):=G_n/\Inn (\mS_n)$ of outer
automorphisms of the algebra $\mS_n$ is isomorphic to the group
$S_n\ltimes \mT^n$.
\end{corollary}

{\it Proof}. By Theorem \ref{5Feb9}, ${\rm Out} (\mS_n)=S_n\ltimes
\mT^n\ltimes \Inn (\mS_n)/\Inn (\mS_n)\simeq S_n\ltimes \mT^n$.
$\Box$

$\noindent $

The next corollary describes the image and the kernel of the group
homomorphism $\xi : G_n\ra \Aut_{K-{\rm alg}}(L_n)$, see
(\ref{xiaut}).

\begin{corollary}\label{a10Feb9}
\begin{enumerate}
\item $\im (\xi ) = S_n\ltimes \mT^n$. \item  $\ker (\xi ) = \Inn
(\mS_n)$.
\end{enumerate}
\end{corollary}

{\it Proof}.  By Theorem \ref{5Feb9}, $G_n= S_n\ltimes
\mT^n\ltimes \Inn (\mS_n)$; $\Inn (\mS_n) \subseteq \ker (\xi )$
since $L_n$ is a commutative algebra. Now, the results follow from
the fact that the homomorphism $\xi$ maps isomorphically the
subgroup $S_n\ltimes \mT^n$ of $G_n$ to the subgroup $S_n\ltimes
\mT^n$ of
 $\Aut_{K-{\rm alg}}(L_n)$.  $\Box $

\begin{corollary}\label{a24Feb9}
The group $G_n$ contains an isomorphic copy of each linear
algebraic group over $K$. In particular, $G_n$ contains an
isomorphic copy of each finite group.
\end{corollary}

{\it Proof}.  The result is obvious since the group $G_n$ contains
the group $\GL_\infty (K)$ and any linear algebraic group can be
embedded in $\GL_\infty (K)$.  $\Box $

\begin{corollary}\label{a21Feb9}
\begin{enumerate}
\item $\stCH = \Inn (\mS_n)$. \item {\rm (Characterization of the
inner automorphisms $\Inn (\mS_n)$ via the height 1 primes of
$\mS_n$)} An automorphism $\s \in G_n$ is an inner automorphism
iff $\s (\gp_1) = \gp_1, \ldots ,\s (\gp_n) = \gp_n$ and
$$ \s (x_1) \equiv x_i\mod \gp_i, \;\; \s (y_i) \equiv y_i\mod\gp_i,
\;\; i=1, \ldots , n.$$ \item If $\s \in \Inn (\mS_n)$ then $\s =
\o_\v$ for a unique element $\v \in \mS_n^*/K^*$ and $\s (x_i)  =
x_i+p_i$, $\s (y_i) =y_i+q_i$ where $p_i=[\v, x_i]\v^{-1}$ and $
q_i= [ \v , y_i]\v^{-1}$ for $i=1, \ldots , n$.
\end{enumerate}
\end{corollary}

{\it Proof}. 1. See (\ref{stISn}).

2. Statement 2 is equivalent to  statement 1.

3. \begin{eqnarray*}
 \v x_i\v^{-1}&=& \s (x_i) = x_i+p_i \Leftrightarrow p_i=[\v, x_i]\v^{-1}, \\
  \v y_i\v^{-1}&=& \s (y_i) = y_i+q_i \Leftrightarrow q_i= [ \v ,
  y_i]\v^{-1}. \;\;\; \Box
\end{eqnarray*}
The inner automorphism $\s \in \Inn (\mS_n)$ can be defined in two
different ways:

(i) $\s = \o_\v$ for a unique element $\v \in \mS_n^*/K^*$; or

(ii) by the elements $p_i:= \s (x_i) - x_i$, $q_i:=\s (y_i) -
y_i$, $i=1, \ldots , n$.

Corollary \ref{a21Feb9}.(3) explains how to pass from (i) to (ii).
 The reverse passage, i.e. from (ii) to (i), is more subtle.
 Suppose that the elements $\{ p_i, q_i \, | \, i=1, \ldots , n\}$
 are given. Below, it is explained how to construct the element
 $\v\in \mS_n^*\subseteq E_n$ which is unique up to $K^*$. By
 Theorem \ref{a5Feb9}, the map  $\v : P_n\ra {}^\s P_n$ is an
 isomorphism of the $\mS_n$-modules $P_n$ and ${}^\s P_n$ (which
 is unique up to $K^*$ since $\End_{\mS_n}(P_n) \simeq K$, \cite{shrekalg}). The isomorphism $\v$  is
 determined by the polynomial $v:= \v (1)\in P_n$ which is unique
 up to $K^*$:
 $$ Kv = \bigcap_{i=1}^n \ker_{P_n} (\s (y_i))= \bigcap_{i=1}^n \ker_{P_n}
 (y_i+q_i).$$
 Then $\v$ is the change-of-the-basis matrix
 $$ x^\alpha \mapsto \prod_{i=1}^n (x_i+p_i)^{\alpha_i}*v.$$
 Note that $\{ x^\alpha \}_{\alpha \in \N^n}$ and $\{ \s (x^{\alpha})*v=\prod_{i=1}^n
 (x_i+p_i)^{\alpha_i}*v\}_{\alpha \in \N^n}$ are two bases for the
 vector space $P_n$.

The next corollary shows that the groups $G_n$ and $\Inn (\mS_n)$
have trivial centre as well as some of the subgroups of $G_n$.
\begin{corollary}\label{b24Feb9}
\begin{enumerate}
\item $Z(G_n) = \{ e\}$.  \item $Z(\mT^n\ltimes \Inn (\mS_n)) = \{
e\}$.\item $Z(\Inn (\mS_n)) = \{ e\}$. \item $Z(S_n\ltimes \mT^n)
= \{ t_{(\l , \ldots , \l )}\, | \, \l \in K^*\}\simeq
\mT^1$.\item $Z(S_n\ltimes \Inn (\mS_n)) = \{ e\}$.
\end{enumerate}
\end{corollary}

{\it Proof}. 3. To prove statement 3 we use induction on $n$. The
 case $n=1$ is true (Theorem \ref{a13Dec8}). So, let $n>1$ and we
 assume that the statement holds for all $n'<n$. Since $\Inn
 (\mS_n)\simeq \mS_n^* / K^*$, we have show that $Z(\mS_n^*)=K^*$.
 Let $z\in \Z(\mS_n^*)$. For each $i=1, \ldots , n$, let
 $\mS_{n-1, i}:=\bigotimes_{j\neq i} \mS_1(j)$ and consider the
 obvious algebra homomorphisms:
 $$ \mS_n\ra \mS_n / \gp_i \simeq K[x_i, x_i^{-1}]\t \mS_{n-1,
 i}\ra K(x_i) \t \mS_{n-1, i}.$$
 By induction, the centre of the group of units of the algebra $ K(x_i) \t \mS_{n-1,
 i}$ is $K(x_i)^*$, hence the image of the element $z$ under the
 first map ($a\mapsto a+\gp_i$) belongs to the Laurent polynomial
 algebra $K[x_i, x_i^{-1}]$. This implies that $z\in \CL_1(i)
 +\gp_i$ where $\CL_1(i) := (\bigoplus_{j\geq 1} Ky_i^j) \bigoplus
 K \bigoplus (\bigoplus_{j\geq 1} Kx_i^j)$, and so
 $$ z\in \bigcap_{i=1}^n (\CL_1(i) +\gp_i) \subseteq \bigcap_{i=1}^n
 (K+\gp_i) \subseteq K+F_n.$$ In particular, $z\in
 Z((K+F_n)^*)=K^*$ since $K+F_n\simeq K+M_\infty (K)$ and
 $Z((K+M_\infty (K))^*)=K$ (see Theorem \ref{a13Dec8}).

4. This is obvious.

2. Let $z=t_\l \o_u\in Z(\mT^n\ltimes \Inn (\mS))$ where $t_\l \in
\mT^n$ and $\o_u\in \Inn (\mS_n)$. For all sufficiently large
natural numbers $k$ and $l$, the elements of the group $\mM_n^*$,
$u$ and $v(k,l,i):= 1+E_{kl}(i)$, $i=1, \ldots , n$ commute.
Therefore, the elements $t_\l$ and $\o_{v(k,l,i)}$ commute. By
(\ref{tlEab}), $t_\l = e$, and so $z=\o_u\in Z(\mT^n\ltimes \Inn
(\mS))\cap \Inn (\mS_n ) \subseteq Z(\Inn (G_n))= \{ e \}$ (by
statement 3), hence $z=e$.

1. Let $z\in Z(G_n)$. Then $z=\tau t_\l \o_u$ for some elements
$\tau \in S_n$, $t_\l \in \mT^n$, and $\o_u\in \Inn (G_n)$. The
element $\tau$ is the image of the element $z$ under the group
epimorphism $G_n\ra G_n/\mT^n\ltimes \Inn (\mS_n) \simeq  S_n$.
The element $\tau$ belongs to the centre of the group $S_n$ which
is equal to $Z(S_n) = \begin{cases}
S_2& \text{if }n=2,\\
e & \text{if }n\neq 2.\\
\end{cases}$ Therefore, $\tau = e$ if $n\neq 2$. If $n=2$ then the
element $\tau t_\l$ is the image of the element $z$ under the
group epimorphism $G_2\ra G_2/\Inn (\mS_2)\simeq S_2\ltimes
\mT^2$, and so it belongs to the centre of the group $S_2\ltimes
\Inn (\mS_2)$, and so $\tau = e$,  by statement 4. Therefore, in
general, $\tau =e$, and so  $z\in Z(G_n) \cap \mT^n\ltimes \Inn
(\mS_n) \subseteq Z(\mT^n\ltimes \Inn (\mS_n))=\{ e \}$ (by
statement 2), hence $z=e$.

5. Let $z=\tau \o_u\in Z(S_n\ltimes \Inn (\mS_n))$. Using the same
arguments as in the proof of statement 2, the elements $\tau $ and
$\o_{v(k,l,i)}$ commute for all large natural numbers $k$ and $l$,
and all $i=1, \ldots , n$. Then $\tau = e$, by (\ref{tlEab1}), and
so $z=\o_u\in Z(S_n\ltimes \Inn (\mS_n))\cap \Inn (\mS_n)
\subseteq Z(\Inn (\mS_n)) = \{ e \}$ (by statement 3), hence
$z=e$. $\Box $


\section{A membership criterion for elements of the algebra $\mS_n$}\label{MEMB}
This section is independent of Section \ref{GAS}.  In this
section, membership criteria for the algebras $\mS_n$, $P_n+F_n$,
and $K+F_n$ are found in terms of commutators (Theorem
\ref{10Feb9},  Corollaries \ref{a13Feb9} and \ref{a15Feb9}).  The
most difficult result of this section is Theorem \ref{10Feb9}
which is used in the proof of Theorem \ref{5Feb9}. Corollary
\ref{a15Feb9} is used in the proof of Theorem \ref{20Feb9}. A
general result of constructing algebras using commutators is
proved (Theorem \ref{a11Feb9}) which shows that the obtained
criteria are rather special (and tight).

For each $i=1, \ldots , n$, the equality (\ref{mS1d}) can be
written as follows 
\begin{equation}\label{SLiF}
\mS_1(i) = \CL_1(i)\bigoplus F(i) \;\ {\rm where}\;\;
\CL_1(i):=(\bigoplus_{j\geq 1}Ky_i^j)\bigoplus K\bigoplus
(\sum_{j\geq 1}Kx_i^j)=\bigoplus_{j\in \Z}Kv_j(i),
\end{equation}
where $$ v_j(i):=\begin{cases}
x_i^j& \text{if } j\geq 0,\\
y_i^{-j}& \text{if } j<0.\\
\end{cases}
$$
 So, each element $a\in \mS_1(i)$ can be uniquely written as a
sum
$$ a= \sum_{j\geq 1}\l_{-j}y_i^j+\l_0+\sum_{j\geq
1}\l_jx_i^j+\sum_{k,l\in \N}\l_{kl}E_{kl}(i)=\sum_{j\in
\Z}\l_jv_j(i)+\sum_{k,l\in \N}\l_{kl}E_{kl}(i)$$ where the
coefficients are scalars. On the other hand, each element $a\in
\mS_1(i)$ is a unique sum $a= \sum_{k,l\in \N} \mu_{kl}x_i^ky_i^l$
where $\mu_{kl}\in K$. Using the formula (\ref{Eijc}) the second
presentation of the element $a$ can be easily obtained from the
first one; and the other way round can be done  using the formula
(\ref{razxy}) below.

For all $i,j\in \N$, 
\begin{equation}\label{razxy}
x^iy^j=\begin{cases}
x^{i-j}-\sum_{k=0}^{j-1}E_{i-j+k,k}& \text{if }\; i\geq j ,\\
y^{j-i}-\sum_{k=0}^{i-1}E_{k, j-i+k}& \text{if }\; i<j.\\
\end{cases}
\end{equation}
It suffices to prove the equality (\ref{razxy}) in the case when
$i\geq j$ since then the second case can be obtained from the
first case: indeed,  for $i<j$,
$$ x^iy^j=x^iy^iy^{j-i}=(1-\sum_{k=0}^{i-1}E_{kk})y^{j-i}=
y^{j-i}-\sum_{k=0}^{i-1}E_{k,j-i+k}.$$ To prove the first case we
use induction on $j$. The result is obvious for $j=0$. So, let
$j>0$ and we assume that the formula (\ref{razxy}) holds for all
$j'<j$. Using induction and the equality $xy= 1-E_{00}$, we have
the result:
\begin{eqnarray*}
 x^iy^j&=& x^iy^{j-1}y= (x^{i-j+1}-\sum_{k=0}^{j-2}E_{i-(j-1)+k,k})y\\
 &=&
 x^{i-j}(1-E_{00})-\sum_{k=0}^{j-2}E_{i-j+k+1,k+1}=x^{i-j}-\sum_{k=0}^{j-1}E_{i-j+k,k}.
\end{eqnarray*}

Let $\CB_n$ be the set of all functions $f:\{ 1, 2, \ldots , n\}
\ra \F_2:= \{ 0,1\}$ where $\F_2:= \Z / 2\Z$ is the  field that
contains two elements. $\CB_n$ is a commutative ring with respect
to addition and multiplication of functions.  For $f,g\in \CB_n$,
we write $f\geq g$ iff $f(i) \geq g(i)$ for all $i=1, \ldots , n$
where $1>0$. Then $(\CB_n, \geq )$ is a partially ordered set. For
each function $f\in \CB_n$, let $|f|:= \sum_{i=1}^nf_i=\#\{i\, |
\, f_i=1\}$ and $\mS_{n,f}:=\bigotimes_{i=1}^n \mS_{1, f_i}(i)$
where
$$
\mS_{1, f_i}(i):= \begin{cases}
\CL_1(i)& \text{if }\; f_i=1,\\
F(i)& \text{if }\; f_i=0.\\
\end{cases}
$$
By (\ref{SLiF}) and $\mS_n = \bigotimes_{i=1}^n \mS_1(i)$, we have
the direct sum 
\begin{equation}\label{SnSnf}
\mS_n = \bigoplus_{f\in \CB_n}\mS_{n,f},
\end{equation}
and so each element $a\in \mS_n$ is a unique sum
\begin{equation}\label{1SnSnf}
a = \sum_{f\in \CB_n}a_f,
\end{equation}
where $a_f\in \mS_{n,f}$. The vector space
$\CL_n:=\bigotimes_{i=1}^n \CL_1(i)= \bigoplus_{\alpha \in
\Z^n}Kv_\alpha$, where $v_\alpha := \prod_{i=1}^nv_{\alpha_i}(i)$,
is not an algebra but it is an algebra modulo the ideal $\ga_n$
which is canonically isomorphic to the Laurent polynomial algebra
$L_n$ (via $v_\alpha +\ga_n\lra x^\alpha$): $(\CL_n+\ga_n)/\ga_n =
\mS_n/ \ga_n=L_n$. The elements $\{ v_\alpha\}_{\alpha \in \Z^n}$
have remarkable properties which are used in the proof of the
Membership Criterion for the elements of the algebra $\mS_n$
(Theorem \ref{10Feb9}). 
\begin{equation}\label{vax1}
v_\alpha *x^\beta = \begin{cases}
x^{\alpha +\beta}& \text{if }  \alpha + \beta \in \N^n,\\
0& \text{if } \alpha + \beta \not\in \N^n.\\
\end{cases}
\end{equation}
\begin{equation}\label{vax2}
v_\alpha *x^\beta x^\g = x^\beta v_\alpha *x^\g \;\; {\rm if}\;\;
\alpha +\g \in \N^n.
\end{equation}
There is an obvious (useful) criterion of when an element of the
algebra $\mS_n$ belongs to the ideal $F_n$. It is used in the
proof of Theorem \ref{10Feb9}.
\begin{lemma}\label{a14Feb9}
Let $a\in \mS_n$. Then $a\in F_n$ iff $a*(\sum_{i=1}^n
P_nx_i^d)=0$ for some $d\in \N$.
\end{lemma}

{\it Proof}. $(\Rightarrow )$ Trivial.

$(\Leftarrow )$ Let $C_n(d) := \{ \alpha \in \N^n\, | \,
\alpha_1\leq d, \ldots , \alpha_n\leq d\}$ and, for each element
$\alpha \in C_n(d)$,
$$ a*x^\alpha = \sum_{\beta \in \N^n}  \l_{\alpha \beta} x^\beta =
(\sum_{\beta \in \N^n} \l_{\alpha \beta} E_{\beta
\alpha})*x^\alpha$$ for some elements $\l_{\alpha\beta}\in K$, and
so $a=\sum_{\beta \in\N^n} \sum_{\alpha \in
C_n(d)}\l_{\alpha\beta}E_{\beta \alpha}\in F_n$.  $\Box $

$\noindent $

The next theorem is a criterion of when a linear map $\v \in
\End_K(P_n)$ belongs to the algebra $\mS_n$ in terms of
commutators. This result is tight when  we compare it with general
results of that sort, see Theorem \ref{a11Feb9} and Corollary
\ref{b11Feb9}. It is not obvious from the outset that the linear
maps that satisfy the commutator conditions of Theorem
\ref{10Feb9} form an algebra.

\begin{theorem}\label{10Feb9}
{\rm (A Membership Criterion)} Let $\v \in \End_K(P_n)$. Then the
following statements are equivalent.
\begin{enumerate}
\item $\v \in \mS_n$. \item $[x_1, \v]\in \gp_1, \ldots , [x_n, \v
] \in \gp_n $. \item $x_i\v = \v \cdot (x_i+p_i) +q_i$, $i=1,
\ldots , n$, for some elements $p_i, q_i\in \gp_i$.
\end{enumerate}
\end{theorem}

{\it Proof}. $(1\Rightarrow 2)$ Let $\mS_{n-1,
i}:=\bigotimes_{j\neq i}\mS_1(j)$. Recall that $[x_i,
\mS_1(i)]\subseteq F(i)$, by (\ref{comab}) for  $n=1$. Then, for
each $i=1, \ldots , n$,
$$[x_i, \mS_n]\subseteq [x_i, \mS_1(i)]\t \mS_{n-1, i}\subseteq
F(i)\t \mS_{n-1, i}=\gp_i.$$

$(2\Rightarrow 3)$ Trivial.

$(3\Rightarrow 1)$ Suppose that a map $\v$ satisfies the
conditions of statement 3. The key idea of the proof of the fact
that $\v \in \mS_n$ is to use a downward induction on a natural
number $s$ starting with $s=n$ and $\v := \v_{n+1}$ to construct
elements $a_f\in \mS_{n, f}$ ($0\neq f\in \CB_n$), elements
$q_{i,s+1}\in \gp_i$ ($i=1, \ldots , n$; $s=1, \ldots , n$), and
natural numbers $d_n\leq d_{n-1}\leq \cdots \leq d_1$ such that
the maps $\v_s: = \v - \sum_{|f|\geq s}a_f$ satisfy the following
conditions: for all $s=1, \ldots , n$, 
\begin{equation}\label{xivd1}
x_i\v_{s+1} = \v_{s+1} \cdot (x_i+p_i) +q_{i,s+1},\;\; p_i, q_{i,
s+1}\in \mS_{n-1, i}\bigotimes
\bigoplus_{k,l=0}^{d_s-1}KE_{kl}(i), \;\; i=1, \ldots , n,
\end{equation}
\begin{equation}\label{xivd2}
\v_s *(\sum_{0\leq i_1<\ldots <i_s\leq n }P_n(x_{i_1}\cdots
x_{i_s})^{d_s})=0.
\end{equation}
Note that $\v_{n+1}= \v$ and all the maps $\v_s$ satisfy the
assumptions of statement 3 since $[x_i, \mS_n]\subseteq \gp_i$,
$i=1, \ldots , n$. Suppose that we have proved this fact then, for
$s=1$, the condition (\ref{xivd2}) is
$$ (\v -\sum_{|f|\geq 1}a_f)* (\sum_{i=1}^n P_nx_i^{d_1})=0.$$
Then, by Lemma \ref{a14Feb9}, $a_0:=\v -\sum_{|f|\geq 1}a_f\in
F_n$, and so $\v = \sum_{f\in \CB_n}a_f\in \mS_n$, as required.

For $s=n$, by the assumption,   we can fix a natural number $d_n$
such that (\ref{xivd1}) holds, that is
$$
x_i\v_{n+1} = \v_{n+1} \cdot (x_i+p_i) +q_{i,n+1}; \;\; p_i, q_{i,
n+1}\in \mS_{n-1, i}\bigotimes
\bigoplus_{k,l=0}^{d_n-1}KE_{kl}(i), \;\;  i=1, \ldots , n,
$$
where $\v= \v_{n+1}$ and $q_{i, n+1}= q_i$.  We have to construct
 the element $a_f\in \mS_{n,f}= \CL_n$ where $f=(1, \ldots , 1)$ such
that (\ref{xivd2}) holds. Let $\ud_n=(d_n, \ldots , d_n)\in \N^n$.
Then
$$ \v * x^{\ud_n} = \sum_{\beta}\l_\beta x^\beta = (\sum_{\beta}
\l_\beta v_{\beta - \ud_n})*x^{\ud_n}$$ for some scalars $\l_\beta
\in K$. Let $a_f:=\sum_{\beta} \l_\beta v_{\beta - \ud_n}$. Since
$p_i*x^{\alpha + \ud_n}=0$ and $q_{i, s+1}*x^{\alpha + \ud_n}=0$,
we have
$$ \v *x^{\alpha + \ud_n} = x^\alpha \v *x^{\ud_n}\;\; {\rm for
\; all}\;\; \alpha \in \N^n.$$ Using these equalities and
(\ref{vax2}), we see that
$$ \v_n *x^{\alpha + \ud_n} = x^\alpha \v_n *x^{\ud_n}= x^\alpha (\v *x^{\ud_n} -a_f*x^{\ud_n}) =0 \;\; {\rm for
\; all}\;\; \alpha \in \N^n,$$ and so the equality (\ref{xivd2})
holds for $s=n$ and $d_n$.

Suppose that $s<n$ and we have found elements $a_f\in \mS_{n, f}$
($|f|\geq s+1$), elements $q_{i,t}\in \gp_i$ ($t=s+2, \ldots ,
n$),  and natural numbers $d_n\leq d_{n-1}\leq \cdots \leq
d_{s+1}$ that satisfy the conditions (\ref{xivd1}) and
(\ref{xivd2}) for all $s'=s+1, \ldots , n$. For the map
$\v_{s+1}$, using the assumptions of statement 3, we can fix a
sufficiently large natural number $d_s$ such that equalities
(\ref{xivd1}) hold  and that $d_s\geq d_{s+1}$. Note that the
equalities (\ref{xivd1}) hold automatically for all natural
numbers larger  than $d_s$. The precise meaning of the expression
`sufficiently large' will be given explicitly later when we find
the map $\v_s$. For a moment, any fixed value of $d_s$ such that
(\ref{xivd1}) holds and $d_s\geq d_{s+1}$ suits our purpose. For
each element $f\in \CB_n$ with $|f|=s$, the element $a_f$ is
defined as follows. The set $\{ 1, \ldots , n\}$ is a disjoint
union of its two subsets $\{ i_1, \ldots , i_s\}$ and $\{ i_{s+1},
\ldots , i_n\}$ where $f(i_1) = \cdots = f(i_s)=1$ and $f(i_{s+1})
=\cdots = f(i_n)=0$. For each vector $\nu = (\nu_{s+1},  \ldots ,
\nu_n)\in \N^{n-s}$ with all $\nu_k<d_s$, 
\begin{equation}\label{xivd3}
\v_{s+1}*((x_{i_1}\cdots
x_{i_s})^{d_s}x_{i_{s+1}}^{\nu_{s+1}}\cdots x_{i_n}^{\nu_n}) =
\sum_{\alpha\in \N^n}\l_{\alpha \nu}x^\alpha = a_f*((x_{i_1}\cdots
x_{i_s})^{d_s}x_{i_{s+1}}^{\nu_{s+1}}\cdots x_{i_n}^{\nu_n}),
\end{equation}
where $\l_{\alpha\nu}\in K$ and 
\begin{equation}\label{afdef}
a_f:=\sum_{\alpha\in
\N^n}\l_{\alpha\nu}v_{\alpha_{i_1}-d_s}(i_1)\cdots
v_{\alpha_{i_s}-d_s}(i_s)E_{\alpha_{i_{s+1}},
\nu_{s+1}}(i_{s+1})\cdots E_{\alpha_{i_n}, \nu_n}(i_n).
\end{equation}
By (\ref{xivd1}), for all elements $\alpha =(\alpha_1, \ldots ,
\alpha_s) \in \N^s$, 
\begin{equation}\label{vsa}
\v_{s+1}*(x_{i_1}^{\alpha_1}\cdots x_{i_s}^{\alpha_s}
(x_{i_1}\cdots x_{i_s})^{d_s}x_{i_{s+1}}^{\nu_{s+1}}\cdots
x_{i_n}^{\nu_n}) = x_{i_1}^{\alpha_1}\cdots x_{i_s}^{\alpha_s}
\v_{s+1}*((x_{i_1}\cdots
x_{i_s})^{d_s}x_{i_{s+1}}^{\nu_{s+1}}\cdots x_{i_n}^{\nu_n}).
\end{equation}
This equalities hold for any new $d_s$ which is not smaller than
the old $d_s$.

Define $\v_s:= \v_{s+1}-\sum_{|f|=s}a_f$ and choose a new number
$d_s$ which is not smaller than the old $d_s$ and such that
(\ref{xivd1}) holds for the map $\v_s$.  Using the equalities
(\ref{vsa}) (for all possible choices of $f$ with $|f|=s$) and for
the new choice of $d_s$ together with  (\ref{vax2}), the equality
(\ref{xivd2}) follows at once: the ideal $\sum_{0\leq i_1<\cdots <
i_{s+1}\leq n} P_n(x_{i_1}\cdots x_{i_{s+1}})^{d_s}$ is
annihilated both by the map $\v_{s+1}$ (due to (\ref{xivd2}) for
$s+1$ and $d_s\geq d_{s+1}$) and by the element $\sum_{|f|=s}a_f$,
by the choice of $d_s$, hence it is annihilated by the map $\v_s$
(each element $a_f$, where  $|f|=s$,  annihilates this ideal). In
order to prove (\ref{xivd2}) it is sufficient to show that the map
$\v_s$ annihilates the monomials of the type $u=(x_{i_1}\cdots
x_{i_s})^{d_s}x_{i_{s+1}}^{\nu_{s+1}}\cdots x_{i_n}^{\nu_n}$, but
his is obvious since
$$ \v_s *u = (\v_{s+1}-\sum_{|g|=s}a_g) *u = (\v_{s+1}-a_f)*u =0, $$
by (\ref{xivd3}) since $a_g(u) =0$ for all $g\neq f$.  $\Box $

\begin{theorem}\label{a11Feb9}
Let $A\subseteq B$ be $K$-algebras and $M$ be a faithful
$B$-module (and so $A\subseteq B\subseteq \End_K(M)$). Suppose
that $I$ is a left ideal of the algebra $B$ such that $I\subseteq
A$. Then
\begin{enumerate}
\item  the set $A':=\{ b\in B\, | \, [b,A]\subseteq I\}$ is a
subalgebra of $B$. If $[A, A]\subseteq I$ then   $A\subseteq A'$.
\item If $I$ is also an ideal of the algebra $A$, and $\{
a_s\}_{s\in S}$ is a set of $K$-algebra generators for $A$ then
$A'=\{ b\in B\, | \, [b,a_s]\in I$ for all $s\in S\}$.
\end{enumerate}
\end{theorem}

{\it Proof}. 1. The set $A'$ is a vector space over the field $K$,
to prove that $A'$ is an algebra we have to show that
$A'A'\subseteq A'$. Let $b,c\in A'$. Then
\begin{eqnarray*}
 [bc, A]&\subseteq & [b,A]c+b[c,A]\subseteq Ic+bI \\
 &\subseteq & [I,c]+cI +I\subseteq [ A, c]+I\subseteq I.
\end{eqnarray*}
If $[A, A]\subseteq I$ then, obviously,    $A\subseteq A'$.

2. Let, $A'':=\{ b\in B\, | \, [b,a_s]\in I$ for all $s\in S\}$.
Then $A'\subseteq A''$. To prove the reverse inclusion it is
enough to show that $[b, a_{s_1}\cdots a_{s_m}]\in I$ for all
products $u=a_{s_1}\cdots a_{s_m}$ of the generators $\{
a_s\}_{s\in S}$. We use induction on $m$ to prove this fact. The
case $m=1$ is obvious. So, let $m>1$ and we assume that the result
is true for all $m'<m$. Then
$$[b, a_{s_1}\cdots a_{s_m}]= [b, a_{s_1}\cdots a_{s_{m-1}}]a_{s_m}+ a_{s_1}\cdots a_{s_{m-1}} [b, a_{s_m}]  \in
Ia_{s_m}+I\subseteq I. \;\;\; \Box $$

\begin{corollary}\label{b11Feb9}
The set $\mS_1':= \{ \v \in \End_K(P_1)\, | \, [ x,\v]\in F, [y,\v
] \in F\}$ is a subalgebra of $\End_K(P_n)$ such that
$\mS_1\subseteq \mS_1'$. In fact, $\mS_1= \mS_1'$, by Theorem
\ref{10Feb9}.
\end{corollary}

{\it Proof}. This is a direct consequence of Theorem \ref{a11Feb9}
where $A=\mS_1= K\langle x,y\rangle$, $M=P_1$, $B=\End_K(P_1)$,
and $I=F$ is an ideal of $\mS_1$ such that $[\mS_1,
\mS_1]\subseteq F$. It is obvious that the ideal $F$ of the
algebra $\mS_1$ is a left ideal of the endomorphism algebra
$\End_K(P_1)$ since an element $f\in \End_k(P_1)$ belongs to $F$
iff $f*P_1x^d=0$ for some $d\in \N$. $\Box $

$\noindent $

For all integers $i,j\in \N$ (where $E_{i,-1}:=0$ and $E_{-1,
j}:=0$) 
\begin{equation}\label{cxy}
[x,y^i]=-E_{0, i-1}, \;\; [y,x^i]=E_{i-1,0},
\end{equation}
\begin{equation}\label{1cxy}
[x, E_{ij}]= E_{i+1, j}-E_{i,j-1}, \;\;
[y,E_{ij}]=E_{i-1,j}-E_{i,j+1}.
\end{equation}
For an algebra $A$ and an element $a\in A$,  let $\ad (a) := [a,
\cdot ]:b\mapsto [a,b]= ab-ba$ be the {\em inner derivation} of
the algebra $A$  determined by the element $a$. The kernel $\ker
\, \ad (a)$ of the inner derivation $\ad (a)$ is a subalgebra of
$A$.
\begin{lemma}\label{b13Feb9}
\begin{enumerate}
\item $\bigcap_{i=1}^n \ker \, \ad (x_i) = K[x_1, \ldots , x_n]$.
\item $\bigcap_{i=1}^n \ker \, \ad (y_i) = K[y_1, \ldots , y_n]$.
\end{enumerate}
\end{lemma}

{\it Proof}. 1. We use induction on $n$. Let $n=1$ and $a\in \ker
\, \ad (x_1)$. By (\ref{1SnSn}), $a=a_1+a_0$ for unique elements
$a_0\in F$ and $a_1=\sum_{i\geq 1}\l_{-i}y_1^i+p$, $p\in K[x_1]$.
Using the expressions for the commutators $[x_1,y_1^i]$ and
$[x_1,E_{ij}]$ given by (\ref{cxy}) and (\ref{1cxy}), we deduce
that $a_1=p$ and $a_0=0$, and so $a\in K[x_1]$. This proves the
equality in the case $n=1$. Let $n>1$ and we assume that the
result holds for all $n'<n$. By induction,  $\bigcap_{i=1}^{n-1}
\ker_{\mS_{n-1}} \, \ad (x_i) = P_{n-1}$. Since $\mS_n =
\mS_{n-1}\t \mS_1$, we have $\bigcap_{i=1}^{n-1} \ker_{\mS_n} \,
\ad (x_i) = P_{n-1}\t \mS_1(n)$, and finally $\bigcap_{i=1}^n \ker
\, \ad (x_i) = P_n$ since $\ker_{\mS_1(n)}\, \ad (x_n) = K[x_n]$.

2. Applying the involution $\eta$ to statement 1 we obtain
statement 2.  $\Box $

\begin{corollary}\label{a13Feb9}
$\{ \v \in \End_K(P_n)\, | \, [x_1, \v ] \in F_n ,\ldots , [x_n,
\v ] \in F_n\}=\begin{cases}
\mS_1& \text{if } n=1,\\
P_n+F_n& \text{if } n>1.\\
\end{cases}$
\end{corollary}

{\it Proof}. For $n=1$, the result follows from Theorem
\ref{10Feb9}. Let $n>1$.  Let $L$ and $R$ denote the LHS and the
RHS of the equality. Then $L\supseteq R$. Let $a\in L$, it remains
to show that $a\in R$. For each $i=1, \ldots , n$, let $\mS_{n-1,
i}:= \bigotimes_{j\neq i} \mS_1(j)$ and $F_{n-1, i}:=
\bigotimes_{j\neq i} F(j)$.

Note that $\mS_n = \mS_1 \t \mS_{n-1, 1}$ and $[x_1, \mS_1]
\subseteq F$ (see (\ref{comab}) for $n=1$). The inclusion $[ x_1,
a]\in F_n$ implies that $x_1\in K[x_1]\t \mS_{n-1, 1}+\mS_1\t
F_{n-1, 1}$. The conditions $[x_j, a]\in F_n$ for $j=2, \ldots ,
n$, imply that $a\in K[x_1]\t \mS_{n-1, 1}+F_n$ (see
(\ref{1cxy})). Then $a\in K[x_i]\t\mS_{n-1, i}+F_n$ for all $i$
(by symmetry of the indices), and $$a\in \bigcap_{i=1}^n
(K[x_i]\t\mS_{n-1, i}+F_n) = P_n+F_n. \;\;\; \Box
$$

\begin{corollary}\label{a15Feb9}
{\rm (Membership Criterion for $F_n$)}
$$\{ \v \in \End_K(P_n)\, | \, [x_i, \v ] \in F_n, [y_i, \v ] \in
F_n, i=1, \ldots , n\} =\begin{cases}
\mS_1& \text{if }n=1,\\
K+F_n& \text{if }n>1.\\
\end{cases} $$
\end{corollary}

{\it Proof}. This follows from Corollary \ref{a13Feb9} and
(\ref{cxy}).  $\Box $

$\noindent $

{\it Remarks}. 1.  The set in Corollary \ref{a15Feb9} is, in fact,
an algebra which is not obvious from the outset. This fact can be
deduced from Theorems \ref{10Feb9} and \ref{a11Feb9}: let $L$ be
 the LHS of the equality in Corollary \ref{a15Feb9}. Since $F_n\subseteq \gp_i$ for all $i$, $L\subseteq
 \mS_n$, by Theorem \ref{10Feb9}. Then $L$ is a subalgebra of
 $\mS_n$ by applying Theorem \ref{a11Feb9} in the case $A=B=
 \mS_n$ and $I=F_n$.

2. Corollaries \ref{b11Feb9} and \ref{a15Feb9} also show that in
order to have the inclusion $A\subseteq A'$ in Theorem
\ref{a11Feb9}.(1), the condition $[A,A]\subseteq I$ cannot be
dropped: for $n>1$, let $L$ be as above. By Theorem \ref{10Feb9},
$L\subseteq \mS_n$, and so $L=\{ b\in \mS_n\, | \, [b,x_i]\in F_n,
[b,y_i]\in F_n,  i=1, \ldots , n \}$, $I=F_n$ is an ideal of
$A=B=\mS_n$. Since $[\mS_n , \mS_n]\not\subseteq F_n$ and $L=
K+F_n\not\supseteq A$, we see that in Theorem \ref{a11Feb9} the
condition $[A,A]\subseteq I$ cannot be dropped and still have the
inclusion $A\subseteq A'$.


\section{The groups $\mM_n^*$ and $G_n'$}\label{CENGN}

In this section, the subgroups $\mM_n^*$ and $G_n'$ of the groups
$\mS_n^*$ and $G_n$ respectively are introduced. It is proved that
the group $\mM_n^*$ has trivial centre (Corollary \ref{a19Feb9})
and is a skew direct product of $2^n-1$ copies of the group
$\GL_\infty (K)$ (Theorem \ref{18Feb9}).  An analogue of the
polynomial Jacobian homomorphism, the so-called  global
determinant, is introduced for the group $\mM_n^*$. In Section
\ref{JACDET}, the global determinant is extended to the group
$G_n'$.

For each non-empty subset $I$ of the set of indices $\{ 1, \ldots
, n\}$, define the $K$-algebra without 1,
$$ F(I):=\bigotimes_{i\in I} F(i) = \bigoplus_{\alpha , \beta \in
\N^I} KE_{\alpha\beta}(I)\subseteq \mS_n, \;\;
E_{\alpha\beta}(I):=\prod_{i\in I}E_{\alpha_i\beta_i}(i), $$ where
 $\alpha= (\alpha_i)_{i\in I}$ and $\beta = (\beta_i)_{i\in I}$. The algebra $F(I)$ is isomorphic {\em non-canonically} to
the matrix algebra (without 1) $M_\infty (K) = \bigcup_{d\geq
1}M_d(K)$ when we fix a bijection $b: \N^m\ra \N$. Then the matrix
unit $E_{\alpha\beta}(I)$ becomes the usual matrix unit
$E_{b(\alpha ) b(\beta )}$ of the matrix algebra $M_\infty (K)$.
The function $b$ determines the finite dimensional monomial vector
space filtration $\CV_b:= \{ V_{b,i}:=\sum_{b(\alpha ) \leq i}
Kx^\alpha \}_{i\in \N}$ on $P_n$.  The algebra (without 1) $F(I)$
is  an ideal of  the following algebra with 1,
$$\mF_I:= K+F(I)\subseteq \mS_n .$$
 The algebra
$\mF_I$ contains the multiplicative monoid $\mM_I:=1+F(I)\simeq
1+M_\infty (K)$. We define the (global) determinant on $\mM_I$  as
in (\ref{gldet}): 
\begin{equation}\label{MIgldet}
\det = \det_{I, b} : \mM_I\ra K,  \;\; u\mapsto \det (u).
\end{equation}
We will see that the determinant $\det_{I, b}$ does not depend on
the bijection $b$.  The (global) determinant has usual properties
of the determinant. In particular, for all elements $u,v\in
\mM_I$,
$$\det (uv) = \det (u) \cdot \det (v).$$
The group of units $\mM_I^*$ of the monoid $\mM_I$ is
\begin{equation}\label{MIGI}
\mM_I^*=\{ u\in \mM_I\, | \, \det (u)\neq 0\}\simeq \GL_\infty
(K).
\end{equation}
It contains the normal subgroup $ S\mM_I^*=\{ u\in \mM_I\, | \,
\det (u)=1\} \simeq \SL_\infty (K)$ which is the kernel of the
group epimorphism $\det : \mM_I^* \ra K^*$. The inversion formula
for $u^{-1}$ is, basically, the Cramer's formula for the inverse
of a matrix of finite size. The group of units $\mF_I^*$ of the
algebra $\mF_I$ is $$\mF_I^*= K^*\mM_I^*\simeq K^*\times
\mM_I^*\simeq K^*\times \GL_\infty (K).$$
\begin{corollary}\label{a18Feb9}
Let $I$ be a non-empty subset of $\{ 1, \ldots , n\}$. Then
$\mM_I^* = \{ u\in \mM_I\, | \, \det (u) \neq 0\} \simeq
\GL_\infty (K)$ and $Z(\mM_I^*) = \{ 1\}$.
\end{corollary}

{\it Proof}. This follows from Theorem \ref{a13Dec8}. $\Box $

$\noindent $

{\it Definition}. Let $\mF_n:= \bigotimes_{i=1}^n \mF_{\{ i\} } =
K\bigoplus (\bigoplus_{\emptyset \neq I \subseteq \{ 1, \ldots ,
n\} } F(I))\subseteq \mS_n$ (this is a subalgebra of $\mS_n$)  and
$\mM_n:= 1+\sum_{\emptyset \neq I \subseteq \{ 1, \ldots , n\} }
F(I)$, this is a multiplicative submonoid of the algebra $\mF_n$.

$\noindent $

The group of units $\mF_n^*$ of the algebra $\mF_n$ is
$$ \mF_n^* = K^* \mM_n^* \simeq K^*\times \mM_n^*$$
where $\mM_n^*$ is the group of units of the monoid $\mM_n$. The
algebra $\mF_n$ contains all the algebras $\mF_I$, the monoid
$\mM_n$ contains all the monoids $\mM_I$, and the group $\mM_n^*$
contains all the groups $\mM_I^*$.

 Let $X_1, \ldots , X_m$ be a
non-empty subsets of a group $G$ and $X_1\cdots X_n := \{
x_1\cdots x_n\,  | \, x_i\in X_i\}$ be their {\em ordered}
product. We sometime write ${}^{set}\prod_{i=1}^nX_i$ for this
product in order to distinguish it from the direct product of
groups.  In general, $X_1\cdots X_n$ is not a subgroup of $G$. If
each element of the product  $X_1\cdots X_n$ has a unique
presentation $x_1\cdots x_n$ where $x_i\in X_i$ the we say that
the product is {\em exact} and write
$X={}^{exact}\prod_{i=1}^nX_i$.

\begin{theorem}\label{18Feb9}
$\mM_n^*\simeq \underbrace{\GL_\infty (K)\ltimes\cdots \ltimes
\GL_\infty (K)}_{2^n-1 \;\; {\rm times}}$.
\end{theorem}

{\it Proof}. The theorem follows from the fact that there is a
chain of normal subgroups of the group $\mM_n^*$: 
\begin{equation}\label{Mnch}
\mM_n^*=\mM_{n,1 }^*\supset \mM_{n,2 }^*\supset\cdots \supset
\mM_{n, i}^*\supset\cdots \supset \mM_{n,n }^*\supset\mM_{n,
n+1}^*=\{ 1\}
\end{equation}
such that, for each number $s=1, \ldots , n$, 
\begin{equation}\label{Mnch1}
\mM_{n,s}={}^{set}\prod_{|I|=s}\mM_I^*\cdot \mM_{n, s+1}^*\;\;
{\rm and}\;\; \mM_{n,s }^*/ \mM_{n,s+1 }^*\simeq \prod_{|I|=s}
\mM_I^*\simeq \GL_\infty  (K)^{n\choose s},
\end{equation}
where the first product is the product of subsets in the group
$\mM_{n,s}^*$ in arbitrary order, and the second product is  the
direct product of  groups (in particular, the  product of sets
${}^{set}\prod_{|I|=s} \mM_I^*$ has trivial intersection with the
group $\mM_{n,s+1}^*$, i.e. $\{ 1 \}$).  The groups $\mM_{n,s}^*$
are constructed below, see (\ref{MMns}). $\Box $

$\noindent $

In their construction the following two lemmas are  used
repeatedly.

\begin{lemma}\label{aa13Dec8}
Let $R$ be a ring and $I_1, \ldots , I_n$ be ideals of the ring
$R$ such that $I_iI_j=0$ for all $i\neq j$. Let $a= 1+a_1+\cdots
+a_n\in R$ where $a_1\in I_1, \ldots , a_n\in I_n$. The element
$a$ is a unit of the ring $R$ iff all the elements $1+a_i$ are
units; and, in this case,  $a^{-1}= (1+a_i)^{-1}
(1+a_2)^{-1}\cdots (1+a_n)^{-1}$.
\end{lemma}

{\it Proof}. Note that the elements $1+a_i$ commute, and  $a=
\prod_{i=1}^n (1+a_i)$. Now, the statement is obvious.  $\Box $

$\noindent $

Let $R$ be a ring, $R^*$ be its group of units, $I$ be an ideal of
$R$ such that $I\neq R$, and let  $ (1+I)^*$ be the group of units
of the multiplicative monoid $1+I$.

\begin{lemma}\label{a14Dec8}
Let $R$ and  $I$ be as above. Then
\begin{enumerate}
\item $R^*\cap (1+I)= (1+I)^*$.\item $(1+I)^*$ is a normal
subgroup of $R^*$.
\end{enumerate}
\end{lemma}

{\it Proof}. 1. The inclusion $R^*\cap (1+I)\supseteq (1+I)^*$ is
obvious. To prove the reverse inclusion, let $1+a\in R^*\cap
(1+I)$ where $a\in I$, and let $(1+a)^{-1} = 1+b$ for some $b\in
R$. The equality $1= (1+a) (1+b)$ can be written as $b= -
a(1+b)\in I$, i.e. $1+a\in (1+I)^*$. This proves the reverse
inclusion.

2. For all $a\in R^*$, $a(1+I)a^{-1}=  a(R^* \cap (1+I))a^{-1} =
aR^*a^{-1} \cap a(1+I)a^{-1} = R^* \cap (1+I) = (1+I)^*$.
Therefore, $(1+I)^*$ is a normal subgroup of $R^*$. $\Box $

$\noindent $

The set $\CF := \bigoplus_{\emptyset \neq I\subseteq \{ 1, \ldots
, n\} }F(I)$ is an ideal of the algebra $\mF_n= K+\CF$. There is
  the strictly descending chain of ideals of the algebra $\mF_n$,
$$\CF \supset \CF^2\supset \cdots \supset \CF^s\supset \cdots
\supset \CF^n = F_n , $$ where $\CF^s:=\bigoplus_{|I|\geq s}F(I)$.
The subalgebra $K+\CF^s$ of $\mF_n$ contains the multiplicative
monoid $\mM_{n,s}:=1+\CF^s$. For each number $s=1, \ldots , n$,
let 
\begin{equation}\label{MMns}
\mM^*_{n,s}:= (1+\CF^s)^*
\end{equation}
be the group of units of the monoid $\mM_{n,s}$, and so we have
the chain of normal subgroups (\ref{Mnch}) of the group $\mM_n^*$.

For each number $s=1, \ldots , n$, consider the factor algebra
$(K+\CF^s)/\CF^{s+1}=K\bigoplus \bigoplus_{|I|=s}J_I$, where
$$ J_I:= (F(I)+\CF^{s+1}) / \CF^{s+1} \simeq F(I)/ F(I)\cap \CF^{s+1}
\simeq F(I)/ 0\simeq F(I)$$ are ideals of the factor algebra such
that $J_IJ_{I'}=0$ if $I\neq I'$. By Lemma \ref{aa13Dec8},  the
group of units of the factor algebra $(K+\CF^s) / \CF^{s+1}$ is
$$ K^* \cdot \prod_{|I|=s}(1+J_I)^* \simeq K^*\times
\prod_{|I|=s}(1+J_I)^*.$$ Then the group $\mM^*_{n,s+1}$ is the
kernel of the group homomorphism 
\begin{equation}\label{Mnch2}
\mM^*_{n,s}\ra \prod_{|I|=s}(1+J_I)^*, \;\; 1+f\mapsto
1+f+\CF^{s+1}.
\end{equation}
Note that $\mM_I^*\subseteq \mM^*_{n,s}$ (where $|I|=s$),   and
the composition of the group homomorphisms
$$ \mM_I^* \ra \mM^*_{n,s}\ra  \prod_{|I'|=s}(1+J_{I'})^*\ra
1+J_{I'}$$ is an isomorphism if $I'=I$ and is the trivial
homomorphism if $I'\neq I$ (i.e. $\mM_I^*\ra 1$). Therefore, the
image of the homomorphism (\ref{Mnch2}) is isomorphic to the
direct product of  groups $\prod_{|I|=s} \mM_I^* \simeq \GL_\infty
(K)^{n\choose s}$, and (\ref{Mnch1}) follows. This completes the
proof of Theorem \ref{18Feb9}.

For each number $s=1, \ldots , n$, let
$\mM_{n,[s]}^*:={}^{set}\prod_{|I|=s}\mM_I^*$ be the product of
the sets $\mM_I^*$, $|I|=s$, in the group  $\mM_n^*$ in an {\em
arbitrary} but {\em fixed} order. By (\ref{Mnch1}), there is a
natural {\em bijections} between the sets 
\begin{equation}\label{Mnch3}
\mM_{n,[s]}^*\ra \prod_{|I|=s}\mM_I^*, \;\; u\mapsto
\prod_{|I|=s}u_i,
\end{equation}
where the RHS is the direct product of groups. So, each element
$v$ of the set $\mM_{n,[s]}^*$ is a {\em unique} product
$\prod_{|I|=s} v_I$ (in the fixed order) of elements $v_I$ of the
groups $\mM_I^*$.

\begin{corollary}\label{a20Feb9}
$\mM_n^* =\mM_{n,[1]}^*\mM_{n,[2]}^*\cdots \mM_{n,[n]}^*$  and
there is a natural bijection (determined by (\ref{Mnch3})),
$$ \mM_n^* \ra
{}^{exact}\prod_{s=1}^n{}\prod_{|I_s|=s}\mM_{I_s}^*, \;\; u\mapsto
\prod_{s=1}^n\prod_{|I_s|=s} u_{I_s},$$ where $u_{I_s}\in
\mM_{I_s}^*$. So, each element $u$ of $\mM_n^*$ is a unique
product $u=\prod_{s=1}^n{}\prod_{|I_s|=s}u_{I_s}$ where where
$u_{I_s}\in \mM_{I_s}^*$.
\end{corollary}

{\it Proof}. The result follows from (\ref{Mnch1}) and
(\ref{Mnch2}).  $\Box $

$\noindent $

For a group $G$, let $Z(G)$ denote its centre. The next corollary
shows that the group $\mM_n^*$ has trivial centre.

\begin{corollary}\label{a19Feb9}
$Z(\mM_n^*)=\{ 1\}$.
\end{corollary}

{\it Proof}. This follows from (\ref{Mnch}), (\ref{Mnch1}), and
the fact that $Z(\GL_\infty (K))=\{ 1\}$.  $\Box $

$\noindent $

The next theorem gives a characterization of the subgroup
$\CM_n:=\{ \o_u\, | \, u\in \mM_n^*\}\simeq \mM_n^*$, $\o_u\lra
u$, of $G_n$. Clearly, $\CM_n\subseteq \Inn (\mS_n)$.

\begin{theorem}\label{20Feb9}
The subgroup $\CM_n:=\{ \o_u\, | \, u\in \mM_n^*\}$ of $G_n$ is
equal to $ \{ \s \in G_n \, | \, \s (x_i) - x_i, \s (y_i) - y_i\in
\mF_n, i=1, \ldots , n\}$. Moreover, for each element $\s \in
\CM_n$,
$$ \s = \prod_{|I_1|=1}\o_{u(I_1)}\cdot
\prod_{|I_2|=2}\o_{u(I_2)}\cdots \prod_{|I_s|=s}\o_{u(I_s)}\cdots
\prod_{|I_n|=n}\o_{u(I_n)}$$ for unique elements $u(I_s)\in
\mM_{I_s}^*$ where the orders in the products are arbitrary but
fixed.
\end{theorem}

{\it Proof}. The inclusion $\{ \o_u\, | \, u\in \mM_n^*\}\subseteq
W_n:=\{ \s \in G_n \, | \, \s (x_i) - x_i, \s (y_i) - y_i\in
\mF_n, i=1, \ldots , n\}$ is obvious since
$$ \o_u (x_i) - x_i= [u,x_i]u^{-1}\in \mF_n, \;\; \o_u (y_i) - y_i= [u,y_i]u^{-1}\in
\mF_n, \;\; i=1, \ldots , n.$$ To prove the reverse inclusion it
suffices to show existence of the product for each element $\s \in
W_n$.

Uniqueness follows from Corollaries  \ref{a20Feb9} and
\ref{a19Feb9} since the RHS is equal to $\o_u$ where
$$ u =\prod_{|I_1|=1}u(I_1)\cdot  \prod_{|I_2|=2}u(I_2)\cdots
\prod_{|I_s|=s}u(I_s)\cdots \prod_{|I_n|=n}u(I_n).$$
 It follows from the explicit action of the group $S_n\ltimes
 \mT^n$ on the elements $x_i$ and $y_i$ ($i=1, \ldots , n$) and
 the equalities
 $G_n= S_n\ltimes \mT^n\ltimes \Inn (\mS_n)$ and  $\Inn (\mS_n) =
 \stCH $, that $W_n =\{ \s \in \Inn (\mS_n) \, | \, \s (x_i) - x_i, \s (y_i) - y_i\in
\mF_n, i=1, \ldots , n\}$.  Since $\Inn (\mS_n) = \stCH$ and $\s
\in W_n$, we have the inclusions (see Corollary \ref{a21Feb9}.(2))
\begin{equation}\label{sxiyi1}
\s (x_i) \in x_i+F(i) +F(i) \CF , \;\; \s (y_i) \in y_i+F(i) +F(i)
\CF , \;\; i=1, \ldots , n.
\end{equation}
It remains to prove existence of the elements $u(I_s)$. We use
induction on $n$. The case $n=1$ is obvious (Theorem
\ref{A5Feb9}). Let $n>1$ and we assume that the statement holds
for all $n'<n$. Let us find the elements $u(I_1)$, $|I_1|=1$, i.e.
the elements  $u(i)$, $i=1, \ldots , n$. Since $\s \in \Inn
(\mS_n)= \stCH$, $\s (\sum_{j\neq i}\gp_j) = \sum_{j\neq i}\gp_j$
for each number $i=1, \ldots , n$. Therefore, the automorphism
$\s$ induces an automorphism, say  $\s_i$, of the factor algebra
$$\mS_n / \sum_{j\neq i}\gp_j\simeq L_{n,i}\t
\mS_1(i),$$ where $L_{n,i}:= \bigotimes_{j\neq i} L_1(j)$, such
that $\s_i(x_j) = x_j$ for all $j\neq i$, and $\s_i(\mS_1(i)
\subseteq \mS_1(i)$, by (\ref{sxiyi1}). Then
$$\s_i(\mS_1(i))=\mS_1(i).$$ By induction, there exists an
element $u(i) \in (1+F(i))^*$ such that the inner automorphism
$\o_{u(i)}$ of the algebra $\mS_n$ induces on the factor algebra
$\mS_n/\sum_{j\neq i}\gp_j$ the automorphism $\s_i$. Let
$\o_{[1]}:=\prod_{i=1}^n \o_{u(i)}$, where the order is fixed as
in the theorem, and let $\s_{[2]}:= \o_{[1]}^{-1}\s$. Then
$$ \s_{[2]}(x_i) -x_i, \s_{[2]}(y_i)-y_i\in \bigoplus_{i\in I,
|I|\geq 2}F(I), \;\; i=1, \ldots , n.$$ Suppose that $s>1$ and we
have already found the elements $u(I)$, $|I|<s$, that satisfy the
following conditions: for all $t=2, \ldots , s$, 
\begin{equation}\label{stF}
\s_{[t]}(x_i) -x_i, \s_{[t]}(y_i)-y_i\in \bigoplus_{i\in I,
|I|\geq t}F(I), \;\; i=1, \ldots , n,
\end{equation}
where $\s_{[t]}:=\o_{[t-1]}^{-1} \cdots \o_{[1]}^{-1}\s$ and
$\o_{[r]}:= \prod_{|I_r|=r}\o_{u(I_r)}$. To finish the proof of
the theorem  by  induction on $s$ we have to find the elements
$u(I_s)$, $|I_s|=s$, such that the automorphism
$\s_{[s+1]}:=\o^{-1}_{[s]}\s_{[s]}$ satisfy (\ref{stF}) for
$t=s+1$ where $\o_{[s]}:=\prod_{|I|=s}\o_{u(I)}$, the order as in
the theorem.

{\em Case (i)}: $s<n$.  For each subset $I$ of $\{ 1, \ldots ,
n\}$, let $CI$ denote its complement. Let $|I|=s$ and $\gp_{CI}:=
\prod_{j\in CI}\gp_j$. Then  $\s_{[s]}(\gp_{CI})=\gp_{CI}$.
Therefore, the automorphism $\s_{[s]}$ induces an automorphism
$\s_{[s], I}$ of the factor algebra
$$\mS_n / \gp_{CI}\simeq L_{CI}\t \mS_I$$ where $L_{CI}:=
\bigotimes_{j\in CI}L_1(j)$ and $\mS_I:=\bigotimes_{j\in
I}\mS_1(j)$, such that $\s_{[s], I}(x_j) = x_j$ for all $j\in CI$,
and $\s_{[s], I}(\mS_I)\subseteq \mS_I$, by (\ref{stF}).
Therefore, $$\s_{[s], I}(\mS_I)= \mS_I.$$ Moreover,
$$\s_{[s], I}(x_i) -x_i, \; \s_{[s], I}(y_i)-y_i\in
F(I)=\bigotimes_{j\in I}F(j)\subseteq \mS_I,  \;\; i=1, \ldots ,
n.$$ Since $|I|=s<n$, by induction on $n$, there is an element
$u(I) \in \mM_I^*$ such that the inner automorphism $\o_{u(I)}$ of
the algebra $\mS_n$ induces the automorphism $\s_{[s],I}$. The
automorphism $\s_{[s+1]}= \o^{-1}_{[s]}\s_{[s]}$ satisfies the
condition (\ref{stF}) for $t=s+1$ where
$\o_{[s]}=\prod_{|I|=s}\o_{u(I)}$, the order as in the theorem.

{\em Case (ii)}: $s=n$. In this case, we cannot use the induction
on $n$ as we did in the previous case. Instead, we are going to
use the Membership Criterion (Corollary \ref{a15Feb9}) in the case
$n>1$. For $s=n$, the condition (\ref{stF}) states that $$p_i:=
\s_{[n]}(x_i) - x_i, \; q_i:=\s_{[n]}(y_i) - y_i\in F_n, \;\; i=1,
\ldots , n.$$ By Theorem \ref{a5Feb9}, $\s_{[n]}(a) = \v a
\v^{-1}$  (where $a\in \mS_n$) for some element $\v\in
\Aut_K(P_n)$. Then $\v x_i= (x_i+p_i) \v$ and $\v y_i= (y_i+q_i)
\v$, and so
$$ [\v, x_i]=p_i\v = \v\v^{-1}p_i\v= \v \s_{[n]}^{-1}(p_i)\in
E_n\s_{[n]}^{-1}(F_n) = E_nF_n\subseteq F_n$$ since
$\s_{[n]}^{-1}(F_n) = F_n$ (as $F_n$ is the least nonzero ideal of
the algebra $\mS_n$) and $E_nF_n\subseteq F_n$ (by Lemma\
\ref{a14Feb9}). Similarly,
$$ [\v, y_i]=q_i\v = \v\v^{-1}q_i\v= \v \s_{[n]}^{-1}(q_i)\in
E_n\s_{[n]}^{-1}(F_n) = E_nF_n\subseteq F_n. $$ By Corollary
\ref{a15Feb9}, $\v \in (K+F_n)^*= K^*\times (1+F_n)^*$, and so the
element $\v$ can be taken from the group $\mM_{\{ 1, \ldots , n\}
}= (1+F_n)^*$. Then $\s_{[n]}= \o_\v$, and the automorphism
$\s_{[n+1]}:= \o_\v^{-1}\s_{[n]}= e$ satisfies the condition
(\ref{stF}) for $t=n+1$ which states that $\s_{[n+1]}=e$. The
proof of the theorem is complete. $\Box $

$\noindent $

{\bf The group $G_n'$ and its generators}. The monoid $\mM_n$ is
stable under the action of the subgroup $S_n\ltimes \mT^n$ of
$G_n$, hence so is its group $\mM_n^*$ of units. Therefore,
$G_n':= S_n\ltimes \mT^n\ltimes \CM_n$ is a subgroup of $G_n$.

\begin{lemma}\label{x1Apr9}
$G_n'\simeq S_n\ltimes \mT^n\ltimes \underbrace{\GL_\infty
(K)\ltimes\cdots \ltimes \GL_\infty (K)}_{2^n-1 \;\; {\rm
times}}$.
\end{lemma}

{\it Proof}. $G_n'\simeq S_n\ltimes \mT^n\ltimes ( \mM_n^*/
Z(\mM_n^*))\simeq S_n\ltimes \mT^n\ltimes \mM_n^*$ (Corollary
\ref{a19Feb9}) and the statement follows from Theorem
\ref{18Feb9}. $\Box $

$\noindent $

 For each element $u\in
\mM_n^*$, let $\o_u:a\mapsto uau^{-1}$ be the inner automorphism
of $\mS_n$ determined by the element $u$. It follows from Lemma
\ref{x1Apr9} that the group $G_n'$ admits the following set of
generators (in the cases (i) and (ii) only nontrivial action of
automorphisms on the canonical generators is shown):

(i) for each pair $i\neq j$ where $i,j\in \{ 1, \ldots , n\}$,
$$ s_{ij}:x_i\mapsto x_j, \; y_i\mapsto y_j, \; x_j\mapsto x_i, \;
y_j\mapsto y_i;$$ (ii) for each $i=1, \ldots , n$ and $\l \in
K^*$,
$$ t_\l (i) : x_i\mapsto \l x_i, \; y_i\mapsto \l^{-1}y_i;$$
(iii) for each non-empty subset $I$ of $\{ 1, \ldots , n\}$,
elements  $k=(k_i)_{i\in I}, l=(l_i)_{i\in I}\in \N^I$ such that
$k\neq l$, and a scalar $\l \in K$, the inner automorphism $\o_u$
where
$$u=u(I; k,l; \l ):=1+\l\prod_{i\in
I}(x_i^{k_i}y_i^{l_i}-x_i^{k_i+1}y_i^{l_i+1}); $$ (iv) for each
non-empty subset $I$ of $\{ 1, \ldots , n\}$ and a scalar $\l \in
K\backslash \{ -1 \}$, the inner automorphism $\o_v$ where
$$v=v(I,\l ):=1+\l \prod_{i\in
I}(1- x_iy_i). $$


\section{An analogue of the Jacobian map - the global determinant}\label{JACDET}

The aim of this section is to introduce an analogue of the
polynomial Jacobian homomorphism,  the so-called global
determinant on $G_n'$ and to prove that it is a group homomorphism
from $G_n'$ to $K^*$ (Corollary \ref{e7Mar9}).

{\bf The determinant $\det$ on the group $\mM_n^*$}. By Corollary
\ref{a20Feb9}, each element $u\in \mM_n^*$ is a unique ordered
product (i.e. for fixed orders of the multiples in each set
$\mM^*_{[n], i}$)
$$ u = \prod_{s=1}^n\prod_{|I_s|=s}u_{I_s}, \;\; u_{I_s}\in
\mM^*_{I_s}, $$ and $\det_{I_s, b(I_s)}(u_{I_s})\neq 0$.

$\noindent $

{\it Definition}. The scalar $\det (u) :=
\prod_{s=1}^n\prod_{|I_s|=s}\det_{I_s, b(I_s)}(u_{I_s})\in K^*$ is
called the {\em global determinant} of the element $u$ (we will
often drop the adjective `global').

$\noindent $

We are going to prove that the determinant (map): 
\begin{equation}\label{DETM}
\det : \mM_n^* \ra K^*, \;\; u\mapsto \det (u)
\end{equation}
is well-defined (i.e. it does not depend on the orders of the
multiples in the product for $u$, and the functions $b(I_s)$),
moreover, it is a group homomorphism (Theorem  \ref{d7Mar9}).

The group $\GL_n(K)$ is the semi-direct product $U_n(K)\ltimes
E_n(K)$ of its two subgroups: $U_n(K):=\{ \l E_{11}+E-E_{11}\, |
\, \l \in K^*\}\simeq K^*$, $\l E_{11}+E-E_{11}\lra \l$, where $E$
is the $n\times n$ identity matrix, and $E_n(K)$ is the subgroup
of $\GL_n(K)$ generated by the elementary matrices $\{ E+\l
E_{ij}\, | \, \l \in K, i\neq j\}$. The group $E_n(K)$ is the
commutant $[\GL_n(K), \GL_n(K)]$ of the group $\GL_n(K)$. Apart
from the usual definition, the determinant $\det :\GL_n(K)\ra K^*$
can be defined as the group epimorphism $\det : \GL_n(K)\ra
\GL_n(K)/[\GL_n(K), \GL_n(K)]\simeq U_n(K)\simeq K^*$. Similarly,
the determinant map (\ref{DETM}) can be defined in this way (see
Theorem \ref{d7Mar9}), and using this second presentation it is
easy to prove that the determinant map (\ref{DETM}) is a group
homomorphism.

The polynomial algebra $P_n$ is equipped with the {\em cubic}
filtration $\CC := \{ \CC_m:=\sum_{\alpha\in C_m}Kx^\alpha\}_{m\in
\N}$ where $C_m:= \{ \alpha \in \N^n\, | \, {\rm all} \;
\alpha_i\leq m\}$. The filtration $\CC$ is an ascending, finite
dimensional filtration such that $P_n= \bigcup_{m\geq 0}\CC_m$ and
$\CC_m\CC_l\subseteq \CC_{m+l}$ for all $m,l\geq 0$. In the case
when $I= \{ 1, \ldots , n\}$, the next result shows that the
determinant $\det$, defined in (\ref{MIgldet}), does not depend on
the bijection $b$.

\begin{theorem}\label{6Mar9}
Let $\CV = \{ V_i\}_{i\in \N}$ be a finite dimensional vector
space filtration on $P_n$ and $a\in \mM_{\{ 1, \ldots , n\}
}=1+F_n$. Then $a(V_i) \subseteq V_i$ and $\det (a|_{V_i}) = \det
(a|_{V_j})$ for all $i,j\gg 0$. Moreover, this common value of the
determinants does not depend on the filtration $\CV$ and,
therefore, coincides with the determinant in (\ref{MIgldet}) for
$I=\{ 1, \ldots , n\}$.
\end{theorem}

{\it Proof}. Let $a\in 1+F_n$. Then $a= 1+\sum_{\alpha, \beta \in
C_d}\l_{\alpha \beta} E_{\alpha \beta}$ for some $\l_{\alpha
\beta}\in K$ and $d\in \N$. Then $a(\CC_i) \subseteq \CC_i$ for
all $i\geq d$. Note that the global determinant in
(\ref{MIgldet}), for $I=\{ 1, \ldots , n\}$, is equal to the usual
determinant $\det (a|_{\CC_i})$ for $i\geq d$;  then $\im (a-1)
\subseteq \CC_d\subseteq V_e$ for some $e\in \N$. Since $ a=
1+(a-1)$, we have $a(V_i) \subseteq V_i$ and $\det (a|_{V_i})=
\det (a|_{V_e})$ for all $i\geq e$. Note that this is true for an
arbitrary filtration $\CV$. Consider the following finite
dimensional vector space filtration $\CV':= \{ V_i':= \CC_d, i=0.
\ldots , e-1; V_j':= V_j, j\geq e\}$. Then
$$ \det (a) = \det (a|_{\CC_d}) = \det (a|_{V_{e-1}'})= \det
(a|_{V_j'}) = \det (a|_{V_j}), \;\; j\geq e.$$ This completes the
proof of the theorem. $\Box $

\begin{corollary}\label{b7Mar9}
For each non-empty subset $I$ of the set $\{ 1, \ldots , n\}$, the
determinant defined in (\ref{MIgldet}) does not depend on the
function $b$.
\end{corollary}

{\it Proof}. This is simply Theorem \ref{6Mar9} where the
polynomial algebra $P_n$ is replaced by the polynomial algebra
$P_I:= K[x_{i_1}, \ldots , x_{i_s}]$ where $I= \{ i_1, \ldots ,
i_s\}$.  $\Box $

 $\noindent $

Corollary \ref{b7Mar9} shows that the global determinant $\det$,
defined in (\ref{DETM}), does not depend on the choices of the
functions $b(I_s)$.

Each element $u\in \mM_n$ is a unique finite sum
$$ u=1+\sum_{I}\sum_{\alpha, \beta \in
\N^I}\l_{\alpha\beta}(I)E_{\alpha \beta} (I), \;\; \l_{\alpha
\beta }\in K, $$ where $I$ runs through all the non-empty subsets
of the set $\{ 1, \ldots , n\}$.

 $\noindent $

{\it Definition}. The {\em size} $s(u)$ of the element $u$ is the
maximum of all the coordinates of the vectors $\alpha$ and $\beta$
in the sum above for the element $u$ with
$\l_{\alpha\beta}(I)\neq 0$.

 $\noindent $

For all elements $u,v\in \mM_n$, $s(uv)\leq \max \{ s(u), s(v)\}$.

\begin{lemma}\label{b8Mar9}
Let $u\in \mM_m^*$ and $u=\prod_{s=1}^n\prod_{|I_s|=s}u_{I_s}$ be
its unique ordered product where $u_{I_s}\in \mM_{I_s}^*$. Then
the size $s(u)$ of the element $u$ is the maximum of the sizes
$s(u_{I_s})$ of the elements $u_{I_s}$.
\end{lemma}

{\it Proof}. Let $u_{[s]}:=\prod_{|I_s|=s}u_{I_s}$. Then
$u=u_{[1]}\cdots u_{[n]}$. The statement is obvious if $u=u_{[i]}$
for some $i$ (multiply out the elements in the product). Moreover,
by the Cramer's formula for the inverse of a matrix,
$s(u_{I_s}^{-1}) = s(u_{I_s})$ for all $I_s$ (indeed, it is
obvious that $s(u_{I_s}^{-1})\leq s(u_{I_s})$ but then
$s(u_{I_s})=s((u_{I_s}^{-1})^{-1})\leq s(u_{I_s}^{-1})$, and the
claim follows). This implies that $s(u_{[i]}^{-1}) = s(u_{[i]})$
since $u_{[i]}^{-1} =\prod_{|I_i|=i}u_{I_i}^{-1}$ (in the reverse
order to the original order) and $u_{I_i}^{-1}\in \mM_{I_i}$.
Clearly,
$$ s(u_{[i]}u_{[i+1]}\cdots u_{[n]})\geq s(u_{[i]})\;\; {\rm for \;
all}\;\; i.$$ We use  a downward induction on $i$ starting with
$i=n$ to prove that if $u=u_{[i]}\cdots u_{[n]}$ then the
statement of the lemma  holds. The statement is obvious for $i=n$,
i.e. when $u=u_{[n]}=u_{\{ 1, \ldots , n\} }$. Suppose that $i<n$,
$u=u_{[i]}\cdots u_{[n]}$ and the statement is true for all
$i'>i$. Suppose that the statement is not true for the element
$u$, we seek a contradiction. Then, $s(u_{[i]})\leq
s(u)<s(u_{[i+1]}\cdots u_{[n]})$, by induction. On the other hand,
$s(u_{[i+1]}\cdots u_{[n]})= s(u_{[i]}^{-1}u) \leq \max \{
s(u_{[i]}^{-1}), s(u)\} = \max \{ s(u_{[i]}),
s(u)\}<s(u_{[i+1]}\cdots u_{[n]})$, a contradiction. $\Box$

 \begin{corollary}\label{c8Mar9}
Let $u\in \mM_n^*$. Then $s(u^{-1})=s(u)$.
\end{corollary}

{\it Proof}. Let $u=\prod_{s=1}^n\prod_{|I_s|=s}u_{I_s}$ where
$u_{I_s}\in \mM_{I_s}^*$. Then $s(u^{-1}_{I_s})\leq s(u_{I_s})$,
hence  $s(u^{-1}) = s(\prod_{s=1}^n\prod_{|I_s|=s}u_{I_s}^{-1})\,
[ {\rm in\;  the\;  reverse\;  order}]\leq \max \{
s(u_{I_s}^{-1})\, | \, I_s\}\leq \max \{ s(u_{I_s})\, | \, I_s\}=
s(u)$, by Lemma \ref{b8Mar9}. Now, $s(u^{-1}) \leq s(u)
=s((u^{-1})^{-1})\leq s(u)$, and so $s(u^{-1}) =s(u)$. $\Box $

\begin{lemma}\label{c7Mar9}
Let $u\in \mM_I$ where $I$ is a non-empty subset of $\{ 1, \ldots
, n\}$. Then $u(\CC_i) \subseteq \CC_i$ and $u(\CC_i(I)) \subseteq
\CC_i(I) $ for all $i\geq s(u)$ (where $\CC_i(I)$ is defined in
the proof).
\end{lemma}

{\it Proof}. For $I=\{ 1, \ldots , n\}$, this is simply Theorem
\ref{6Mar9} (see the proof of Theorem \ref{6Mar9} where if $\CV =
\CC$ the elements $d$ and $e$ can be set to be equal to $s(u)$).
The case when $I\neq \{ 1, \ldots , n\}$ follows from the previous
one when we observe that $P_n= P_I\t P_{CI}$ where $P_I:=
K[x_{i_1}, \ldots , x_{i_s}]$, $I=\{ i_1, \ldots , i_s\}$, and
$CI$ is the complement of $I$. Then $\CC_i= \CC_i(I)\t \CC_i(CI)$
where $\{ \CC_i(I)\}_{i\in \N}$ and $\{ \CC_i(CI)\}_{i\in \N}$ are
the cubic filtrations for  the polynomial algebras $P_I$ and
$P_{CI}$ respectively. Note that $u|_{\CC_i}= u|_{\CC_i(I)\t
\CC_i(CI)} = u|_{\CC_i(I)}\t {\rm id}_{\CC_i(CI)}$  for all $i\geq
s(u)$. $\Box$

The group $\GL_\infty(K)$ is the semi-direct product $U(K)\ltimes
E_\infty (K)$ of its two subgroups: $U(K):=\{ \l E_{00}+1-E_{00}\,
| \, \l \in K^*\}\simeq K^*$, $\l E_{00}+1-E_{00}\lra \l$, and
$E_\infty (K)$ is the subgroup of $\GL_\infty (K)$ generated by
the elementary matrices $\{ 1+\l E_{ij}\, | \, \l \in K, i\neq
j\}$. The group $E_\infty(K)$ coincides with the  commutant
$[\GL_\infty (K), \GL_\infty (K)]$ of the group $\GL_\infty (K)$.

For each non-empty subset $I$ of $\{ 1, \ldots , n\}$, the group
$\mM_I^*$ is isomorphic to the group $\GL_\infty (K)$. Therefore,
$\mM_I^*=U_I(K)\ltimes E_I(K)$ is the semi-direct product of its
subgroups: $U_I(K):=\{ \l E_{00}(I)+1-E_{00}(I)\, | \, \l \in
K^*\}\simeq K^*$, $\l E_{00}(I)+1-E_{00}(I)\lra \l$, and $E_I (K)$
is the subgroup of $\mM_I^*(K)$ generated by the elementary
matrices $\{ 1+\l E_{\alpha \beta}(I)\, | \, \l \in K, \alpha,
\beta \in \N^I,  \alpha\neq \beta \}$. The group $E_I(K)$
coincides with the commutant $[\mM_I^*, \mM_I^*]$ of the group
$\mM_I^*$.

For $u\in U_I(K)$ and $u'\in U_{I'}(K)$, $uu'=u'u$ as follows from
$$ (\l E_{00}(I)+1-E_{00}(I))*x^\alpha= \begin{cases}
\l x^\alpha & \text{if }\forall i\in I:\alpha_i=0,\\
x^\alpha & \text{otherwise}.\\
\end{cases}
$$
So, the elements $u$ and $u'$ are diagonal matrices in the
monomial basis for $P_n$. By Corollary \ref{a20Feb9}, the subgroup
$\mU_n$ of $\mM_n$ generated by the groups $U_I(K)$ is equal to
their direct product, $\mU_n=\prod_{I\neq \emptyset } U_I(K)\simeq
K^{*(2^n-1)}$. Consider the group epimorphism \marginpar{mhk}
\begin{equation}\label{mhk}
\mu :\mU_n\ra K^*, \;\; \prod_{I\neq \emptyset }(\l_I
E_{00}(I)+1-E_{00}(I))\mapsto  \prod_{I\neq \emptyset }\l_I.
\end{equation}
For each number $s=1, \ldots , n$, let $\mU_{n,[s]}:=
\prod_{|I|=s} U_I(K)$ and $\mU_{n,s}:=\mU_{n,[s]}\times
\mU_{n,[s+1]}\times \cdots \times \mU_{n, [n]}$. By Corollary
\ref{a20Feb9}, for each $s=1, \ldots , n$, the set $E_{n,[s]}:=
\prod_{|I|=s}E_I(K)$ is an exact product of groups in arbitrary
but fixed order, and $E_{n,s}:= E_{n,[s]}E_{n,[s+1]}\cdots
E_{n,[n]}$ is the exact product of sets. We will see that the set
$E_{n,s}$ is a group.

\begin{theorem}\label{d7Mar9}
\begin{enumerate}
\item  $\mM_n^*=\mU_n\ltimes [ \mM_n^*, \mM_n^*]$ and $[ \mM_n^*,
\mM_n^*]=E_{n,1}$. \item  $\mM_{n,s}^*=\mU_{n,s}\ltimes [
\mM_{n,s}^*, \mM_{n,s}^*]$ and $[ \mM_{n,s}^*,
\mM_{n,s}^*]=E_{n,s}$ for all $s=1, \ldots , n$. \item The
determinant map $\det $ (see (\ref{DETM})) is the composition of
the  group homomorphisms (see (\ref{mhk})):
$$ \det : \mM_n^*\ra \mM_n^*/  [ \mM_n^*, \mM_n^*]\simeq
\mU_n\stackrel{\mu}{\ra}K^*.$$ In particular, $\det(uv ) = \det
(u) \det (v)$  for all $u,v\in \mM_n^*$.
\end{enumerate}
\end{theorem}

{\it Proof}. 1. Statement 1 is a part of statement 2 when $s=1$.

2. To prove statement 2 we use a downward induction on $s$
starting with $s=n$. In this case,  both statements follow at once
from the fact that $\mM_{n,n}^* =(1+F_n)^*\simeq \GL_\infty (K)=
U(K)\ltimes E_\infty (K)$ and $E_\infty (K)=[\GL_\infty (K),
\GL_\infty (K)]$ is the subgroup of $\GL_\infty (K)$ generated by
 the elementary matrices. Suppose that $s<n$ and the statements
 hold for all $s'=s+1, \ldots , n$. By the uniqueness of the
 product in Corollary \ref{a20Feb9}, $\mU_n\cap E_{n,s}=\{ 1\}$.
 It is obvious that $E_{n,s}\subseteq [\mM_n^*, \mM_n^*]$ and
 $\mM_n^* \supseteq \mU_nE_{n,s}$. Recall that the groups
 $\mM_{n,t}^*$ are normal subgroups of the group $\mM_n^*$. It
 follows that the set $E_{n,s}=E_{n,[s]}E_{n,s+1}= E_{n,[s]}[\mM_{n, s+1}^*, \mM_{n,
 s+1}^*]$ is a subgroup of $\mM_{n,s}^*$. Using elementary
 matrices and the generators for the group $\mU_{n,s}$ it is easy
 to verify that
 \marginpar{aa1}
\begin{equation}\label{aa1}
uE_{n,[s]}u^{-1}\subseteq E_{n,s}\;\; {\rm for \; all}\;\; u\in
\mU_{n,s}\;\; {\rm and\; all}\;\; s.
\end{equation}
Note that each element $u\in \mU_{n,s}$ is a diagonal matrix in
the monomial basis for $P_n$.  This implies that
$E_{n,[s]}\mU_{n,n+1}\subseteq \mU_{n,n+1}E_{n,s}$. Now,
\begin{eqnarray*}
 \mM_{n,s}^*&=& \mU_{n,[s]}E_{n,[s]}\mM_{n,s+1}^* = \mU_{n,[s]}E_{n,[s]}\mU_{n,s+1}E_{n,s+1} \\
 &\subseteq & \mU_{n,[s]}\mU_{n,s+1} E_{n,s}= \mU_{n,s}E_{n,s},
\end{eqnarray*}
and so $\mM_{n,s}^* = \mU_{n,s} E_{n,s}$. Since $E_{n,s}=
E_{n,[s]}E_{n,s+1}= E_{n,[s]}[\mM_{n,s+1}^*, \mM_{n,s+1}^*]$ and
$\mM_{n,s+1}^*$ is a normal subgroup of $\mM_n^*$, we see that
$uE_{n,s}u^{-1}\subseteq E_{n,s}$ for all elements $u\in
\mU_{n,s}$, by (\ref{aa1}), i.e. $E_{n,s}$ is a normal subgroup of
$\mM_{n,s}^*$. Hence, $\mM_{n,s}^* = \mU_{n,s}\ltimes E_{n,s}$.
Then $[\mM_{n,s}^*, \mM_{n,s}^*]\subseteq E_{n,s}$ since the group
$\mU_{n,s}$ is abelian. The opposite inclusion is obvious.
Therefore, $E_{n,s} =[\mM_{n,s}^*, \mM_{n,s}^*]$. By induction,
statement 2 holds.

3. By Corollary \ref{a20Feb9}, each element $u$ of the group
$\mM_n^*$ is the  unique product $\prod_{s=1}^n
\prod_{|I_s|=s}u_{I_s}$ where each element $u_{I_s}\in
\mM_{I_s}^*$ is a unique product $u_{I_s}(\l_{I_s})e_{I_s}$ where
$u_{I_s}(\l_{I_s}):= \l_{I_s}E_{00}(I_s)+1-E_{00}(I_s)$ and
$e_{I_s}\in E_{I_s}(K)$. Then $\det (u) =
\prod_{s=1}^n\prod_{|I_s|=s}\l_{I_s}$. By statement 2, the element
$u$ is a unique product
$\prod_{s=1}^n\prod_{|I_s|=s}u_{I_s}(\l_{I_s})\cdot e$ where $e\in
E_{n,1}$, and statement 3 follows. $\Box$

%

$\noindent $

{\bf The global determinant $\det $ on the group $G_n'$}. Recall
that $G_n'\simeq S_n\ltimes \mT^n\ltimes \mM_n^*$, it is
convenient to identify these two groups via the isomorphism. Each
element $\s$ of $G_n'$ is a unique product $\s = \tau t_\l u$
where $\tau \in S_n$, $t_\l \in \mT^n$, and $u\in \mM_n^*$.

$\noindent $

{\it Definition}. The scalar $\det (\s ) := {\rm sgn} (\tau) \cdot
\prod_{i=1}^n \l_i\cdot \det (u)\in K^*$ is called the {\em global
determinant} of the element $\s$ (we often drop the adjective
`global') where ${\rm sgn}(\tau )$ is the parity of  $\tau $.

$\noindent $

Our next goal is to prove that the determinant map
$$ \det : G_n' \ra K^*, \;\; \s \mapsto \det (\s ), $$
is a group homomorphism (Corollary \ref{e7Mar9}).

The group $S_n\ltimes \mT^n$ can be seen as a subgroup of the
general linear group $\GL (V)$ where $V= \bigoplus_{i=1}^n Kx_i
\subseteq P_n$ ($\tau (x_i) = x_{\tau (i)}$ and $t_\l (x_i) =
\l_ix_i$). The global determinant $\det (\tau t_\l ) $ of the
element $\tau t_\l\in S_n\ltimes \mT^n$ is simply the usual
determinant of the element $\tau t_\l \in \GL (V)$. So, in order
to prove Corollary \ref{e7Mar9} it suffices to show that $\det
(\tau t_\l u (\tau t_\l )^{-1}) = \det (u)$ for all $u\in \mM_n^*$
and $\tau t_\l \in S_n\ltimes \mT^n$. This follows from Theorem
\ref{d7Mar9}.(1) and the fact that the element $\tau t_\l$
respects the groups $\mU_n$ and $[\mM_n^*, \mM_n^*]$, and, for
each element $u=\prod_{I\neq \emptyset}u_I\in \mU_n$, the
conjugation $\tau t_\l u(\tau t_\l )^{-1}$  permutes the
components $u_I\in U_I(K)$.

\begin{corollary}\label{e7Mar9}
$\det (ab) = \det (a)\, \det (b)$ for all $a,b\in G_n'$.
\end{corollary}

{\bf The global determinant $\det $ on the monoids $\mM_n$ and
$S_n\ltimes \mT^n\ltimes \mM_n$}. Lemma \ref{c7Mar9} and  Theorem
\ref{d7Mar9} give an idea of how to extend the global determinant
from the group $\mM_n^*$ to the monoid $\mM_n$. Let $u\in \mM_n$
and $s(u)$ be its size. Then $u(\CC_i) \subseteq \CC_i$ for all
$i\geq s(u)$. If the map $u\in \End_K(P_n)$ is a bijection then,
by Theorem \ref{A8Mar9}, $u\in \mM_n^*$. If the map $u$ is not a
bijection then $\det (u|_{\CC_i}) = 0$ for all $i\gg 0$. Hence, if
$u,v\in \mM_n$ and $uv\in \mM_n^*$ then $u,v\in \mM_n^*$ (this
proves the first statement of Theorem \ref{8Mar9}).


$\noindent $

{\it Definition}. We can extend the (global) determinant $\det $
 to the map
 $$ \det : \mM_n\ra K, \;\; u\mapsto \begin{cases}
\det (u)& \text{if } u\in \mM_n^*,\\
0& \text{otherwise}.\\
\end{cases}$$

This common value $\det (u)$ of the determinants is called the
{\em global determinant} of the element $u\in \mM_n$ (we often
drop the adjective `global').

$\noindent $

The global determinant is a homomorphism from the monoid $\mM_n$
to the multiplicative monoid $(K, \cdot )$ (Theorem
\ref{8Mar9}.(2)), and the group $\mM_n^*$ of units of the monoid
$\mM_n$ is the set of all the elements of $\mM_n$ with nonzero
global  determinant (Corollary \ref{a7Mar9}). These results are
based on Theorem \ref{A8Mar9}. We keep the notation of Section
\ref{GAS}. The monoid $\mM_n= 1+\CF$ has the descending monoid
filtration $$\mM_n = 1+\CF \supset 1+\CF^2\supset \cdots \supset
1+\CF^n = 1+F_n.$$ For each element $u\in \mM_n$, there is a
unique number $i$ such that $u\in (1+\CF^i)\backslash
(1+\CF^{i+1})$. The number $i$ is called the {\em degree} of the
element $u$, denoted $\deg (u)$.

For each non-empty subset $I$ of $\{ 1, \ldots , n\}$, let $\CC
(I):=\{ \CC_i(I)\}_{i\in \N}$ be the cubic filtration for the
polynomial algebra $P_I:=K[x_j]_{j\in I}$.
\begin{theorem}\label{A8Mar9}
$\mM_n^*= \mM_n\cap \Aut_K(P_n)$ but $\mS_n^* \subsetneqq
\mS_n\cap \Aut_K(P_n)$.
\end{theorem}

{\it Proof}. Let $u\in \mM_n\cap \Aut_K(P_n)$. We have to show
that $u\in \mM_n^*$ since the inclusion $\mM_n^*\subseteq
\mM_n\cap \Aut_K(P_n)$ is obvious. We prove this fact by a
downward induction on the degree $i=\deg (u)$. If $i=n$, that is
$u\in (1+F_n)\cap \Aut_K(P_n) = (1+F_n)^*$, the statement is
obvious. Suppose that $i<n$, and the statement holds for all
elements $u'$ with $\deg (u') >i$.  In particular,
$(1+\CF^{i+1})\cap \Aut_K(P_n) \subseteq \mM_n^*$. Note that $u=
1+\sum_{|I|=i}a_I+\sum_{|I|>i}a_I$ for unique elements $a_I\in
F(I)$. Let $u_I:= 1+a_I$ and $u':= \prod_{|I|=i}u_I$ (in arbitrary
order). Note that $s(u_I)\leq s(u)$ for all $I$ such that $|I|=i$.
For each natural number $m>s(u)$, let $B_m(I) := \CC_m(I)\t
(\prod_{j\in CI}x_j^m\cdot P_{CI})$. By the choice of $m$,
\begin{equation}\label{uBmI}
u|_{B_m(I)}= u_I|_{B_m(I)},
\end{equation}
and so the linear map $u_I:\CC_m(I) \ra \CC_m(I)$ is an injection,
hence a bijection (since $\dim_K(\CC_m(I))<\infty$) for all
$m>s(u)$. Now,
$$ u_I\in (1+F(I))\cap \Aut_K(P_I)= (1+F(I))^*=\mM_I^* \subseteq
\mM_n^*.$$ Then $u'\in \mM_n^*$, and
$$ u(u')^{-1}\in (1+\CF^{i+1})\cap \Aut_K(P_n)\subseteq \mM_n^*,$$
therefore, $u=u(u')^{-1}\cdot u'\in \mM_n^*$.

$\mS_n^* \subsetneqq \mS_n\cap \Aut_K(P_n)$ since the element $u:=
\prod_{i=1}^n (1-y_i)$ of the algebra $\mS_n$ belongs to the set
$\Aut_K(P_n) \backslash \mS_n^*$. The element $u$ is not a unit of
the algebra $\mS_n$ since the element $u+\ga_n$ is not a unit of
the algebra $\mS_n/ \ga_n$. To show the inclusion $u\in
\Aut_K(P_n)$ we may assume that $n=1$ since $P_n=
\bigotimes_{i=1}^n K[x_i]$. The kernel of the linear map $u$ is
equal to zero since $(1-y)*p=0$ for an element $p\in K[x]$ implies
that $p=y*p= y^2*p= \cdots = y^s*p=0$ for all $s\gg 0$ ($y$ is a
locally nilpotent map). The map $u$ is surjective since for each
element $q\in K[x]$ there exists a natural number, say $t$, such
that $y^t*q=0$, and so $q=(1-y^t)*q=u(1+y+\cdots +y^{t-1})*q$.
Therefore, $u\in \Aut_K(P_n)$.  $\Box $

\begin{theorem}\label{8Mar9}
\begin{enumerate}
\item If $u,v\in \mM_n$ and $uv\in \mM_n^*$ then $u,v\in \mM_n^*$.
\item $\det
(uv) = \det (u)\, \det (v)$ for all elements $u,v\in \mM_n$.
\end{enumerate}
\end{theorem}

{\it Proof}. 
%
%
%
2. The second statement follows from the first. $\Box $

\begin{corollary}\label{a7Mar9}
\begin{enumerate}
\item $\mM_n^*= \{ u\in \mM_n\, | \, \det (u) \neq 0\}$, i.e. an
element $u\in \mM_n$ is a unit iff $\det (u)\neq 0$. \item Let
$u\in \mM_n$. Then the following statements are equivalent.
\begin{enumerate} \item The element $u$ has
left inverse in $\mS_n$ ($vu=1$ for some $v\in \mS_n$). \item The
element $u$ has right  inverse in $\mS_n$ ($uv=1$ for some $v\in
\mS_n$). \item The element $u$ is invertible in $\mS_n$.
 \item $\det (u)\neq 0$.
\end{enumerate}
\end{enumerate}
\end{corollary}

{\it Proof}. 1. Trivial.


2. Statement 2 follows from statement 1 (using the facts that
$vu=1$ implies $\det (u) \, \det (u) =1$, and $uv=1$ implies $\det
(u) \, \det (v) =1$). $\Box $

$\noindent $

We can extend the global determinant to the monoid $S_n\ltimes
\mT^n\ltimes \mM_n$ by the rule:
$$ \det : S_n\ltimes \mT^n\ltimes \mM_n \ra K, \;\; \tau t_\l
u\mapsto \det (\tau t_\l) \, \det (u),$$ where $\tau \in S_n$,
$t_\l \in \mT^n$, and $u\in \mM_n$. It follows from Corollary
\ref{a11Mar9} that this is a well-defined monoid homomorphism.

We define the {\em size} $s(a)$ of an element $a= \tau t_\l u\in
S_n\ltimes \mT^n\ltimes \mM_n $ as $s(u)$. Then $s(ab) \leq \max
\{ s(a) , s(b) \}$ for all $a,b\in S_n\ltimes \mT^n\ltimes \mM_n $
and $s(a^{-1})=s(a)$ for all $a\in S_n\ltimes \mT^n\ltimes \mM_n^*
$, by Lemma \ref{c8Mar9}.

\begin{corollary}\label{a11Mar9}
\begin{enumerate}
\item Let $a\in S_n\ltimes \mT^n\ltimes \mM_n $. Then $u(\CC_i)
\subseteq \CC_i$ 
 for
all $i,j> s(a)$.  \item $\det (ab) = \det (a)\, \det (b)$ for all
elements $a,b\in S_n\ltimes \mT^n\ltimes \mM_n $.
\end{enumerate}
\end{corollary}

\begin{corollary}\label{b11Mar9}
\begin{enumerate}
\item The group of units of the monoid $S_n\ltimes \mT^n\ltimes
\mM_n$ is $S_n\ltimes \mT^n\ltimes \mM_n^*\simeq G_n'$. \item
$S_n\ltimes \mT^n\ltimes \mM_n^*= \{ a\in S_n\ltimes \mT^n\ltimes
\mM_n\, | \, \det (a) \neq 0\}$. \item $S_n\ltimes \mT^n\ltimes
\mM_n^*= (S_n\ltimes \mT^n\ltimes \mM_n)\cap \Aut_K(P_n)$.  \item
Let $a\in S_n\ltimes \mT^n\ltimes \mM_n$. Then the following
statements are equivalent.
\begin{enumerate} \item The element $u$ has
left inverse. \item The element $u$ has right  inverse. \item The
element $u$ is invertible.
 \item $\det (u)\neq 0$.
\end{enumerate}
\end{enumerate}
\end{corollary}


\section{Stabilizers in $\Aut_{K-{\rm alg}}(\mS_n )$  of the prime or idempotent ideals of
$\mS_n$}\label{STABSN}

In this section, for each nonzero idempotent ideal $\ga $ of the
algebra $\mS_n$ its stabilizer $\St_{G_n}(\ga ) :=\{ \s \in G_n \,
| \, \s (\ga )=\ga \}$ is found (Theorem \ref{A10Feb9}). If, in
addition, the ideal $\ga$ is generic this result can be refined
even further (Corollary \ref{aA10Feb9}) where the wreath product
of groups appears. The stabilizers of all the prime ideals of the
algebra $\mS_n$ are found (Corollary \ref{b10Feb9}.(2) and
Corollary \ref{a8Mar9}). In particular, when $n>1$ the stabilizer
of each height 1 prime  of $\mS_n$ is a maximal subgroup of $G_n$
of index $n$ (Corollary \ref{b10Feb9}.(1)). It is proved that the
ideal $\ga_n$ is the only nonzero, prime, $G_n$-invariant ideal of
the algebra $\mS_n$ (Theorem \ref{B10Feb9}).

{\bf Idempotent ideals of the algebra $\mS_n$}. An ideal $\ga$  of
a ring $R$ is called an {\em idempotent} ideal (resp. a {\em
proper} ideal)  if $\ga^2=\ga$ (resp. $\ga\neq 0, R$). For an
ideal $\ga$, $\Min (\ga )$ is the set of all the minimal primes
over $\ga$. Two ideals $\ga$ and $\gb$  are called {\em
incomparable} if neither $\ga \subseteq \gb$ nor $\gb \subseteq
\ga$. The idempotent ideals of the algebra $\mS_n$ are studied in
detail in \cite{shrekalg}. Below (Theorem \ref{A24Feb9}), we
collect results on the idempotent ideals of $\mS_n$ that are used
in the proofs of this section. For the proof of Theorem
\ref{A24Feb9} and for more information on the idempotent ideals of
$\mS_n$ the interested reader is referred to \cite{shrekalg}.

\begin{theorem}\label{A24Feb9}
{\rm (\cite{shrekalg}, Theorem 7.2, Corollary 4.9, Theorem 4.13)}
\begin{enumerate}
\item Let $\ga$ be a proper, idempotent ideal of the algebra
$\mS_n$. Then $\Min (\ga )$ is a finite non-empty set each element
of which is an idempotent, prime ideal of $\mS_n$. The ideal $\ga$
is a unique product and a unique intersection of incomparable,
idempotent, prime  ideals of $\mS_n$. Moreover,
$$ a= \prod_{\gp \in \Min (\ga )}\gp =  \bigcap_{\gp \in \Min (\ga )}\gp .$$
\item Each nonzero, idempotent, prime  ideal $\gp$ of the algebra
$\mS_n$ is equal to $\gp_I:=\sum_{i\in I}\gp_i$ for some non-empty
subset of $\{ 1, \ldots , n\}$ and vice versa; and this
presentation is
 unique. \item The height of the prime ideal $\gp_I$ is $|I|$.
\end{enumerate}
\end{theorem}

\begin{corollary}\label{b10Feb9}
\begin{enumerate}
\item $\St_{G_n}(\gp_i) \simeq S_{n-1}\ltimes \mT^n\ltimes \Inn
(\mS_n)$, for $i=1, \ldots , n$.  Moreover, if $n>1$ then the
groups $\St_{G_n}(\gp_i)$ are maximal subgroups of $G_n$ (if $n=1$
then $\St_{G_1}(\gp_1)=G_1$, by Theorem \ref{B10Feb9}).  \item Let
$\gp$ be a nonzero,  idempotent, prime ideal of the algebra
$\mS_n$ and $h= \hht (\gp )$ be its height. Then $\St_{G_n}(\gp )
\simeq (S_h\times S_{n-h})\ltimes \mT^n\ltimes \Inn (\mS_n)$.\item
$\St_{G_n}(\CH_1) =\mT^n\ltimes \Inn (\mS_n)$.
\end{enumerate}
\end{corollary}

{\it Proof}. 1. Note that $\mT^n\ltimes \Inn (\mS_n) \subseteq
\St_{G_n}(\gp_i)$ and $\St_{G_n}(\gp_i) \cap S_n=\{ \tau \in S_n\,
| \, \tau (\gp_i) = \gp_i\} \simeq S_{n-1}$. Then
\begin{eqnarray*}
 \St_{G_n}(\gp_i) &=&\St_{G_n}(\gp_i)\cap G_n= \St_{G_n}(\gp_i) \cap (S_n\ltimes \mT^n\ltimes \Inn (\mS_n ))\\
 &=&(\St_{G_n}(\gp_i)\cap S_n)\ltimes \mT^n\ltimes \Inn (\mS_n )\simeq S_{n-1}\ltimes \mT^n\ltimes \Inn
(\mS_n).
\end{eqnarray*}
When $n>1$, the group $\St_{G_n}(\gp_i)$ is a  maximal subgroup of
$G_n$ since $$S_{n-1}\simeq \St_{G_n}(\gp_i)/(\mT^n\ltimes \Inn
(\mS_n))\subseteq G_n/ (\mT^n\ltimes \Inn (\mS_n))\simeq S_n$$ and
$S_{n-1}=\{ \s \in S_n \, | \, \s (i)=i \}$ is a maximal subgroup
of $S_n$.

2. By Theorem \ref{A24Feb9}.(2), $\gp=\gp_{i_1}+\cdots +
\gp_{i_h}$ for some distinct indices $i_1, \ldots , i_h\in \{ 1,
\ldots , n\}$. Let $I=\{ i_1, \ldots , i_h\}$ and $CI$ be its
complement. Since $\mT^n\ltimes \Inn (\mS_n) \subseteq
\St_{G_n}(\gp )$ and
$$\St_{G_n} (\gp ) \cap S_n = \{ \s \in S_n \, | \, \s (I) = I, \s
(CI) = CI\}\simeq S_h\times S_{n-h},$$ the result follows using
the same arguments as in the previous case.

3. Statement 3 follows from statement 1.  $\Box $

$\noindent $

Let $\Sub_n$ be the set of all subsets of $\{ 1, \ldots, n \}$.
$\Sub_n$ is a partially ordered set with respect to `$\subseteq
$'. Let $\SSub_n$ be the set of all subsets of $\Sub_n$. An
element $\{ X_1, \ldots , X_s\}$ of $\SSub_n$ is called {\em
incomparable} if for all $i\ne j$ such that $1\leq i,j\leq s$
neither $X_i\subseteq X_j$ nor $X_i\supseteq X_j$. An empty set
and one element set are called incomparable by definition. Let
$\Inc_n$ be the subset of $\SSub_n$ of all incomparable elements
of $\SSub_n$. The symmetric group   $S_n$  acts in the obvious way
on the set $\SSub_n$ ($\s \cdot \{ X_1, \ldots ,X_s\}= \{ \s
(X_1), \ldots , \s (X_s)\}$).

\begin{theorem}\label{A10Feb9}
Let $\ga $ be a proper idempotent ideal of the algebra $\mS_n$.
Then
$$ \St_{G_n}(\ga ) = \St_{S_n}(\Min (\ga ))\ltimes \mT^n\ltimes \Inn
(\mS_n)$$ where $\St_{S_n}(\Min (\ga ))=\{ \s \in S_n \, | \, \s
(\gq ) \in \Min (\ga )$ for all $\gq \in \Min (\ga )\}$. Moreover,
if $\Min (\ga ) = \{ \gq_1, \ldots , \gq_s\}$ and, for each number
$t=1, \ldots , s$, $\gq_t=\sum_{i\in I_t}\gp_i$ for some subset
$I_t$ of $\{ 1, \ldots , n\}$. Then the group $\St_{S_n}(\Min (\ga
))$ is the stabilizer in the group $S_n$  of the element $\{ I_1,
\ldots , I_s\}$ of $\SSub_n$.
\end{theorem}

{\it Remark}. Note that the group $$\St_{G_n}(\Min (\ga ))=
\St_{S_n}(\{ I_1, \ldots , I_s\})$$ (and also the group
$\St_{G_n}(\ga )$) can be effectively computed in finitely many
steps.

{\it Proof}. By Theorem \ref{A24Feb9}.(1,2),  and Corollary
\ref{b10Feb9}, $\mT^n\ltimes \Inn (\mS_n )\subseteq \St_{G_n}(\ga
)$. Note that  $\St_{G_n}(\ga )\cap S_n =\St_{S_n}(\Min (\ga ))$.
Now,
$$ \St_{G_n}(\ga ) =(\St_{G_n}(\ga )\cap S_n)\ltimes \mT^n\ltimes \Inn (\mS_n
)=\St_{S_n}(\Min (\ga ))\ltimes \mT^n\ltimes \Inn (\mS_n).
$$
By Theorem \ref{A24Feb9}.(1), $\St_{S_n}(\Min (\ga ))=
\St_{S_n}(\{ I_1, \ldots , I_s\} )$. $\Box$

$\noindent $

We are going to apply Theorem \ref{A10Feb9} to find the
stabilizers of  the generic idempotent ideals (see Corollary
\ref{aA10Feb9}) but first we recall the definition of the {\em
wreath product} $A\wr B$ of finite groups $A$ and $B$. The set
$\Fun (B,A)$ of all functions $f: B\ra A$ is a group: $(fg) (b) :=
f(b) g(b)$ for all $b\in B$ where $g\in \Fun (B,A)$. There is a
group homomorphism
$$ B\ra \Aut (\Fun (B,A)), \; b_1\mapsto (f\mapsto b_1(f):b\mapsto
f(b_1^{-1}b)).$$ Then the semidirect product $\Fun (B,A) \rtimes
B$ Is called the {\em wreath product} of the groups $A$ and $B$
denoted $A\wr B$, and so the product in $A\wr B$ is given by the
rule:
$$f_1b_1\cdot f_2b_2= f_1b_1(f_2) b_1b_2, \;\; {\rm where}\;\;
f_1, f_2\in \Fun (B,A) , \;\; b_1,b_2\in B.$$ By Theorem
\ref{A24Feb9}.(2), each nonzero, idempotent, prime ideal $\gp$ of
$\mS_n$ is a unique sum $\gp = \sum_{i\in I} \gp_i$ of height 1
prime ideals. The set $\Supp (\gp ):= \{ \gp_i\, | \, i\in I\}$ is
called the {\em support} of $\gp$.

$\noindent $

{\it Definition}. We say that a proper, idempotent ideal $\ga$ of
$\mS_n$ is {\em generic} if $\Supp (\gp ) \cap \Supp (\gq
)=\emptyset$ for all $\gp , \gq \in \Min (\ga )$ such that $\gp
\neq \gq$.

\begin{corollary}\label{aA10Feb9}
Let $\ga$ be a generic idempotent ideal of the algebra $\mS_n$,
the set $\Min (\ga )$ of minimal primes over $\ga$ is the disjoint
union of non-empty subsets $\Min_{h_1}(\ga ) \bigcup \cdots
\bigcup \Min_{h_t}(\ga )$ where $1\leq h_1<\cdots < h_t\leq n$ and
the set  $\Min_{h_i}(\ga )$ contains all the minimal primes over
$\ga$ of height $h_i$. Let $n_i:= |\Min_{h_i}(\ga )|$. Then
$$ \St_{G_n}(\ga )= (S_m\times \prod_{i=1}^t(S_{h_i}\wr
S_{n_i}))\ltimes \mT^n\ltimes \Inn (\mS_n)$$ where $m=
n-\sum_{i=1}^t n_ih_i$.
\end{corollary}

{\it Proof}. Suppose that $\Min (\ga ) = \{ \gq_1, \ldots ,
\gq_s\}$ and the sets $I_1, \ldots , I_s$ are defined in Theorem
\ref{A10Feb9}. Since the ideal $\ga$ is generic, the sets $I_1,
\ldots , I_s$ are disjoint.  By Theorem \ref{A10Feb9}, we have to
show that 
\begin{equation}\label{Smm}
\St_{S_m}(\{ I_1, \ldots , I_s\}) \simeq S_m\times
\prod_{i=1}^t(S_{h_i}\wr S_{n_i}).
\end{equation}
The ideal $\ga$ is generic, and so the set $\{ 1, \ldots , n\}$ is
the disjoint union $\bigcup_{i=0}^t M_i$ of its subsets where
$M_i:= \bigcup_{|I_j|=h_i}I_j$, $i=1, \ldots , t$, and $M_0$ is
the complement of the set $\bigcup_{i=1}^tM_i$. Let $S(M_i)$ be
the symmetric group corresponding to the set $M_i$ (i.e. the set
of all bijections $M_i\ra M_i$). Then each element $\s \in
\St_{G_n}(\{ I_1, \ldots , I_s\} )$ is a unique product $\s =
\s_0\s_1\cdots \s_t$ where $\s_i\in S(M_i)$. Moreover, $\s_0$ can
be an arbitrary element of $S(M_0) \simeq S_m$, and, for $i\neq
0$, the element $\s_i$ permutes the sets $\{ I_j\, | \,
|I_j|=h_i\}$ and simultaneously permutes the elements inside each
of the sets $I_j$, i.e. $\s_i\in S_{h_i}\wr S_{n_i}$. Now,
(\ref{Smm}) is obvious. $\Box $

\begin{corollary}\label{c10Feb9}
For each number $s=1, \ldots, n$, let
$\gb_s:=\prod_{|I|=s}(\sum_{i\in I} \gp_i)$ where $I$ runs through
all the subsets of the index set  $\{ 1, \ldots , n\}$ that
contain exactly $s$ elements.  The ideals $\gb_s$ are the only
proper, idempotent, $G_n$-invariant ideals of the algebra $\mS_n$.
\end{corollary}

{\it Proof}. By Theorem \ref{5Feb9} and Corollary
\ref{b10Feb9}.(3), the ideals  $\gb_s$ are $G_n$-invariant, and
they are proper and idempotent. The converse follows at once from
the classification of proper idempotent ideals (Theorem
\ref{A24Feb9}.(1)). $\Box $

$\noindent $

{\bf The prime ideals of the algebra $\mS_n$}. In order to prove
Theorem \ref{B10Feb9}, we recall a classification of prime ideals
for the algebra $\mS_n$ which is obtained in \cite{shrekalg}. For
a subset $\CN=\{ i_1, \ldots , i_s\}$ of the set of indices $\{ 1,
\ldots , n\}$, let $C\CN$ be its complement, $|\CN |=s$, $\mS_\CN
:= \mS_1(i_1) \t\cdots \t \mS_1(i_s)$, 
\begin{equation}\label{gmCN}
\ga_\CN :=F\t\mS_1(i_2) \t\cdots \t \mS_1(i_s)+\cdots +\mS_1(i_1)
\t\cdots \t \mS_1(i_{s-1})\t F,
\end{equation}
 $P_\CN := K[x_{i_1}, \ldots , x_{i_s}]$.
Clearly, $\mS_n= \mS_\CN \t \mS_{C\CN }$. Let $L_\CN:= K[x_{i_1},
x_{i_1}^{-1}, \ldots , x_{i_s}, x_{i_s}^{-1}]$. Then $\mS_\CN /
\ga_\CN \simeq L_\CN$. Consider the epimorphism 
\begin{equation}\label{pCN}
\pi_\CN : \mS_\CN \ra \mS_\CN/ \ga_\CN\simeq L_\CN, \;\; a\mapsto
a+ \ga_\CN.
\end{equation}
 By \cite{shrekalg}, Proposition
4.3.(2), there is the injection
$$ \spec (L_{C\CN}) \ra \spec (\mS_n) , \;\; \gq \mapsto \mS_\CN \t
\pi_{C\CN}^{-1}(\gq ).$$ The image of this injection is denoted by
$$ \spec (\mS_n, \CN ) :=\{ \mS_\CN \t
\pi_{C\CN}^{-1}(\gq ) \, | \, \gq \in \spec (L_{C\CN} )\}.$$
 Note
that $\spec (\mS_n , \emptyset )= \{ \pi^{-1}_{\{ 1, \ldots , n\}
} (\gq )\, | \, \gq \in \spec (L_n)\}\simeq \spec (L_n)$ and
$\spec (\mS_n , \{ 1, \ldots , n\} )= \{ 0\}$ since $\pi_\emptyset
:K\ra K$, $\l \mapsto \l$.

The next theorem shows that all the prime ideals of the algebra
$\mS_n$ can be obtained in this way.
\begin{theorem}\label{21Dec8}
{\rm (\cite{shrekalg}, Theorem 4.4)}
\begin{enumerate}
\item $\spec (\mS_n) = \coprod_{\CN \subseteq \{ 1, \ldots , n\} }
\spec (\mS_n , \CN )$, the disjoint union.  \item Each  prime
ideal $\gp$ of the algebra $\mS_n$ can be uniquely written as
$\mS_\CN \t \pi^{-1}_{C\CN}(\gq )$ for some subset $\CN$ of the
set $\{ 1, \ldots , n\}$ and some prime ideal $\gq$ of the algebra
$L_{C\CN }$.
\end{enumerate}
\end{theorem}

\begin{theorem}\label{B10Feb9}
The ideal  $\ga_n$ is  the only nonzero, prime, $G_n$-invariant
ideal of the algebra $\mS_n$.
\end{theorem}

{\it Proof}. By Lemma \ref{b5Feb9} (or by Corollary
\ref{b10Feb9}.(2)), the ideal $\ga_n$ is $G_n$-invariant.
Conversely, let $\gp$ be a nonzero, prime, $G_n$-invariant ideal
of the algebra $\mS_n$. By Theorem \ref{21Dec8}.(2) and the fact
that $\gp$ is also $S_n$-invariant, the ideal $\gp$ contains the
sum $\gp_1+\cdots +\gp_n = \ga_n$. Suppose that $\gp \neq \ga_n$,
we seek a contradiction. In this case, the ideal $\gp / \ga_n$ of
the algebra $\mS_n/ \ga_n = L_n$ is $\mT^n$-invariant, hence $\gp
= L_n$, a contradiction. $\Box $

$\noindent $

The classical Krull dimension of the algebra $\mS_n$ is $2n$
(\cite{shrekalg}, Theorem 4.11). For each natural number $i=0, 1,
\ldots , 2n$, let
\begin{eqnarray*}
\CH_i&:=&\{ \gp \in \Spec (\mS_n) \, | \, \hht (\gp ) =
i\},\\
\St_{G_n} (\CH_i) &:=&\{ \s \in G_n\, | \, \s (\gp ) = \gp \;\;
{\rm for \; all}\;\; \gp \in \CH_i\}.
\end{eqnarray*}
\begin{corollary}\label{d10Feb9}
$\St_{G_n} (\CH_i) =\begin{cases}
G_n& \text{if } i=0,\\
\mT^n\ltimes \Inn (\mS_n) & \text{if } i=1,\\
\Inn (\mS_n)& \text{if } i=2, \ldots , 2n.\\
\end{cases}$
\end{corollary}

{\it Proof}. The statement  is obvious for $i=0$ (since $\CH_0=\{
0\}$) and for $i=1$ (Corollary \ref{b10Feb9}.(3)). So, let $i\geq
2$. Briefly, the statement  follows from the fact that in the
algebra $L_n$ there is no proper $\mT^n$-invariant ideals (since
any such an ideal would have contained a monomial in $x_i$,
$x_i^{-1}$, $i=1, \ldots , n$; but all of them are units). Fix a
presentation $i=m+l$ where $1\leq l \leq m \leq n$. For each
subset $\CN$ of $\{ 1, \ldots , n\}$ such that $|C\CN |= m$ and,
for each prime ideal $\gq$ of $L_{C\CN }$ of height $l$,
$$\St_{G_n}(\mS_\CN \t \pi^{-1}_{C\CN}(\gq ))= S(\CN ) \ltimes
\mT^{|\CN |}(\CN ) \ltimes \St_{S(C\CN  \ltimes \mT^{|C\CN
|}(C\CN)}(\gq ) \ltimes \Inn (\mS_n)$$ where $S(\CN )$ is the
symmetric group on $\CN$ and $\mT^{|\CN |}(\CN )$ is the torus in
the group of automorphisms of the algebra $\mS_\CN$. It is obvious
that  $\Inn (\mS_n) \subseteq \St_{G_n}(\CH_i)$. For $i=2, \ldots
, 2n-1$,  $$\bigcap_{\CN , \gq }\St_{G_n}(S_\CN \t \pi^{-1}_{C\CN
}(\gq ))=\Inn (S_n),$$ and so  $\St_{G_n}(\CH_i) = \Inn (\mS_n)$.
For $i=2n$, the statement is obvious.  $\Box $

$\noindent $

Let $\gp $ be a prime ideal of the algebra $\mS_n$. When, in
addition, $\gp$ is an idempotent ideal its stabilizer is found in
Corollary \ref{b10Feb9}.(2). The next corollary, which is obtained
in the proof of Corollary \ref{d10Feb9}, gives the stabilizer of
$\gp$  when the prime ideal $\gp$ is not an idempotent ideal.

\begin{corollary}\label{a8Mar9}
Let $\gp$ be a prime ideal of the algebra $\mS_n$ which is not an
idempotent ideal, i.e. $\gp  = \mS_\CN \t \pi^{-1}_{C\CN}(\gq )$
for some subset $\CN$ of $\{ 1, \ldots , n\}$ and a nonzero prime
ideal $\gq$ of the Laurent polynomial algebra $L_{C\CN}$. Then
$\St_{G_n}(\gp ) = S(\CN ) \ltimes \mT^{|\CN |}(\CN ) \ltimes
\St_{S(C\CN  \ltimes \mT^{|C\CN |}(C\CN)}(\gq ) \ltimes\Inn
(\mS_n)$ (see the proof of Corollary \ref{d10Feb9} for details).
\end{corollary}

\begin{theorem}\label{10Mar9}
\begin{enumerate}
\item Let $n>1$ and let $\gp$ be a prime ideal of the algebra
$\mS_n$. Then the stabilizer $\St_{G_n}(\gp )$ is a maximal
subgroup of $G_n$ iff the ideal $\gp$ has height 1, and in this
case the index $[G_n:\St_{G_n}(\gp )]=n$. \item Let $n=1$ and
$\gp$ be a prime ideal of the algebra $\mS_n$. Then the stabilizer
$\St_{G_n}(\gp )$ is not a maximal subgroup of $G_n$.
\end{enumerate}
\end{theorem}

{\it Proof}. The theorem follows from Corollary \ref{b10Feb9} and
Corollary \ref{a8Mar9}. $\Box $

\begin{corollary}\label{e10Feb9}
$\St_{G_n} (\Spec (\mS_n)) = \St_{G_n} (\Max (\mS_n)) =\Inn
(\mS_n)$.
\end{corollary}

{\it Proof}. By Corollary \ref{d10Feb9},
$$ \Inn (\mS_n) \subseteq \St_{G_n} (\Spec (\mS_n))\subseteq \St_{G_n} (\Max (\mS_n))\subseteq \St_{G_n}(\CH_{2n}) =\Inn
(\mS_n), $$ and so the result. $\Box $

{\bf The algebra $\mS_n$ is $\Z^n$-graded}. The algebra $\mS_n
=\bigoplus_{\alpha \in \Z_n}\mS_{n,\alpha}$ is a $\Z^n$-graded
algebra where $\mS_{n,\alpha}:= \mS_{1,\alpha_1}\t\cdots \t
\mS_{1,\alpha_n}$, $\alpha = (\alpha_1, \ldots , \alpha_n)$,
$$ \mS_{1,i}:=
\begin{cases}
x^i\mS_{1,0}= \mS_{1,0}x^i & \text{if $i\geq 1$},\\
\mS_{1,0}& \text{if $i=0$},\\
y^{|i|}\mS_{1,0}= \mS_{1,0}y^{|i|} & \text{if $i\leq -1$,}
\end{cases}
$$
$\mS_{1,0}:= K\langle E_{00}, E_{11}, \ldots \rangle= K\oplus
KE_{00}\oplus KE_{11}\oplus\cdots$ is a commutative non-Noetherian
algebra $( KE_{00} \subset  KE_{00}\oplus KE_{11}\subset \cdots$
is an ascending chain of ideals of the algebra $\mS_{1,0}$).  For
each $i=1, \ldots , n$, and $j\in \N$, let $$v_j(i) :=
\begin{cases}
x_i^j& \text{if } j\geq 0,\\
y_i^{|j|}& \text{if } j<0,\\
\end{cases}$$ and, for $\alpha \in \Z^n$, let $v_\alpha :=
\prod_{i=1}^n v_{\alpha_i}(i)$. Then $\mS_{n, \alpha } = v_\alpha
 \mS_{n,0}= \mS_{n,0}v_\alpha$ where
 $$\mS_{n,0}:= \bigotimes_{i=1}^n \mS_{1,0}(i) = \bigotimes_{i=1}^n
 K\langle E_{00}(i), E_{11}(i), \ldots \rangle = K\bigoplus
 \bigoplus_{I}\bigoplus_{\alpha \in \N^{|I|}}KE_{\alpha
 \alpha}(I)$$
 where $I$ runs through all the non-empty subsets of $\{ 1,\ldots
 , n\}$, and $E_{\alpha\alpha}(I):= E_{\alpha_1\alpha_1}(i_1)\cdots
 E_{\alpha_s\alpha_s}(i_s)$ for $I=\{ i_1,\ldots , i_s\}$. Each
 element $a$ of the algebra $\mS_{n,0}$ is a unique finite sum
\begin{equation}\label{aEI}
a=a_0+\sum_{I}\sum_{\alpha \in \N^{|I|}}\l_{\alpha ,
I}E_{\alpha\alpha }(I)
\end{equation}
where $a_0$, $\l_{\alpha , I}\in K$. The set of elements $\{ v_\g
, v_\d (CI)E_{\alpha \beta }(I) \}$ is a $K$-basis for the algebra
$\mS_n$ where $E_{\alpha \beta}:= E_{\alpha_1\beta_1}(i_1)\cdots
 E_{\alpha_s\beta_s}(i_s)$ and, for the complement $CI=\{j_1, \ldots , j_t\}$ of the set $I$,
  $v_\d (CI):=v_{\d_1}(j_1)\cdots v_{\d_t}(j_t)$. Each nonzero element $u$ of
$\mS_n$ is a finite linear combination of the basis elements, and
each nonzero summands is called a {\em component} of $u$.

$\noindent $

{\it Definition}. The {\em volume} $\vol (u)$ of a nonzero element
$u$ of $\mS_n$ is the number of nonzero coordinates of the element
$u$ with respect to the basis $\{ v_\g , v_\d (CI) E_{\alpha \beta
}(I) \}$, or, equivalently, the number of its nonzero components.
We set $\vol (0)=0$.

$\noindent $

Note that $\vol (\s (u)) = \vol (u)$ for all $\s \in S_n\ltimes
\mT^n$.

Let $G$ be a group and $H$ be its subgroup. Then $[G:H]$ denotes
the index of $H$ in $G$.

\begin{corollary}\label{a10Mar9}
Let $\ga$ be a proper ideal of the algebra $\mS_n$. Then
$[G_n:\St_{G_n}(\ga ) ] <\infty$ iff $\ga^2=\ga$.
\end{corollary}

{\it Proof}. $(\Leftarrow )$ This implication follows from Theorem
\ref{A10Feb9}.

$(\Rightarrow )$ Suppose that $[G_n:\St_{G_n}(\ga ) ] <\infty$ for
a proper ideal $\ga$ of $\mS_n$. Note that $\mT^n = \prod_{i=1}^n
\mT^1(i)$. For each $i=1, \ldots , n$, let $T_i:= \mT^1(i) \cap
\St_{G_n}(\ga )$. Then $[\mT^1(i):T_i] \leq [G_n:\St_{G_n}(\ga ) ]
<\infty$, and so the group $T_i$ contains {\em infinitely many}
elements. Consider the subgroup $T':= T_1\times \cdots \times T_n$
of $\mT^n\cap \St_{G_n}(\ga )$. We have to show that $\ga^2=\ga$.
It suffices to show that the ideal $\ga$ is generated (as an
ideal) by elements of volume $1$. Suppose that this is not the
case for the ideal $\ga$, we seek a  contradiction. Let $v$ be the
minimum of the volumes of all the nonzero elements of the ideal
$\ga$ such that all their components do not belong to $\ga$. Fix
one such  an element $u\in \ga$ with $\vol (u) = v$. Since
$T'\subseteq \St_{G_n}(\ga )$, the element $u$ has to be of the
type $v_\beta a$ for some $\beta \in \Z^n$ and a nonzero element
$a$ of the algebra $\mS_{n,0}$. The element $a$ is a unique sum as
in (\ref{aEI}). To get a contradiction we use an induction on $n$.
Suppose that $n=1$, and so $u=v_\beta (\l + \sum_{\nu =1}^s a_\nu
E_{i_\nu i_\nu})$ for some scalars $\l$ and $a_\nu \in K^*$, $\nu
\geq 1$.

If $\l \neq 0$ then the ideal of $\mS_1$ generated by the element
$u$ is $\mS_1$. This implies that $u=v_\beta \l$ and so $\vol
(u)=1$, a contradiction.

If $\l =0$ then $uE_{i_\nu i_\nu} = a_\nu v_\beta E_{i_\nu
i_\nu}\in \ga$ for all $\nu$, a contradiction.

Suppose that $n>1$. Then, up to action of the symmetric group
$S_n$, we may assume that
$$u= v_\beta (\l +\sum_{\nu =1}^sa_\nu E_{i_\nu i_\nu}(n))$$
for some scalar $\l \in K$ and nonzero elements $a_\nu \in
\mS_{n-1}$. If $\l \neq 0$ and all $a_\nu \in K$ then the ideal of
the algebra $\mS_1(n)$ generated by the element $v_{\beta_n} (\l
+\sum_{\nu =1}^sa_\nu E_{i_\nu i_\nu}(n))\in \mS_1(n)$ is equal to
$\mS_1(n)$. Then all the summands of the element $u$ belongs to
the ideal $\ga$, a contradiction.

If $\l\neq 0$ and not all the elements $a_\nu $  belong to the
field  $K$, say $a_1\not\in K$, then the volume of the following
nonzero element of $\ga$, $uE_{i_1i_1}(n) = v_\beta (\l + a_1)
E_{i_1i_1}(n)$, is not $1$ and does not exceed $\vol (u)$.
Therefore, $a_2=\cdots = a_s=0$ and $\vol (uE_{i_1i_1}) = \vol
(u)$. Repeating the same argument several times we obtain an
element of the ideal $\ga$,
$$uE_{ii}(k) E_{jj}(k+1)\cdots
E_{i_1i_1}(n)= v_\beta (\l + b) E_{ii}(k) E_{jj}(k+1)\cdots
E_{i_1i_1}(n),$$ having volume $\vol (u)$ but $b\in F_1(k-1)$ (up
to action of the group $S_n$). Since the ideal of the algebra
$\mS_1(k-1)$ generated by its element $v_{\beta_{k-1}}(\l +b)$ is
equal to $\mS_1(k-1)$, we have a contradiction.

If $\l =0$ then all the elements $uE_{i_\nu i_\nu}(n) = v_\beta
a_\nu E_{i_\nu i_\nu}(n)$ belong to the ideal $\ga$. Therefore,
$u=v_\beta a_1 E_{i_1i_1}(n)$ for some nonzero element $a_1\in
\mS_{n-1}$ of volume $\vol (u)$. Now, repeating the same argument
as above or use induction on $n$, we come to a contradiction. The
proof of the corollary is complete.  $\Box $


\section{Endomorphisms of the algebra
$\mS_n$}\label{ENDSN}

In this section, we classify all the algebra endomorphisms of
$\mS_n$ that stabilize the elements $x_1, \ldots , x_n$ and show
that each such endomorphism is a {\em monomorphism} but {\em not}
an isomorphism provided it is not the identity map (Corollary
\ref{a7Feb9}).

Let
\begin{eqnarray*}
{\rm st}   (x_1, \ldots , x_n):= & \{ g\in E_n\, | \, g(x_1) =
x_1,
\ldots , g(x_n)=x_n\}, \\
 {\rm st}   (y_1, \ldots , y_n):= & \{ g\in E_n\,
| \, g(y_1) = y_1, \ldots , g(y_n)=y_n\}.
\end{eqnarray*}
These monoids are the stabilizers of the sets $\{ x_1, \ldots ,
x_n\}$ and $\{ y_1, \ldots , y_n\}$ in $E_n$. Note that  $$\heta
({\rm st} (x_1, \ldots , x_n))={\rm st}   (y_1, \ldots , y_n),
\;\; \heta ({\rm st}   (y_1, \ldots , y_n))= {\rm st} (x_1, \ldots
, x_n).$$ By Theorem \ref{6Feb9},
$$G_n\cap ({\rm st} (x_1, \ldots , x_n)= G_n\cap {\rm st}   (y_1, \ldots ,
y_n)=\{ e\},$$ i.e. if an algebra endomorphism of $\mS_n$ which is
not the identity map stabilizers either the set $\{ x_1, \ldots ,
x_n\}$ or $\{ y_1, \ldots , y_n\}$ then necessarily it is {\em
not} an automorphism of $\mS_n$. Our next step is to describe all
such endomorphisms and to show that all of them are {\em
monomorphisms}. Note that the algebra $\mS_n$ has plenty of ideals
(see \cite{shrekalg}) and contains the ring of infinite
dimensional matrices, so there is no problem in producing  an
algebra endomorphism which is {\em not} a monomorphism, eg
$\mS_n\ra \mS_n/(\ga_n+\sum_{i=1}^n \mS_n (x_i-1)\mS_n) \simeq
K\ra \mS_n$.

$\noindent $

In the proof of Corollary \ref{a7Feb9}, the following identities
are used. For $i=1, \ldots , n$ and $p\in K[x_1, \ldots , x_n]$,
\begin{equation}\label{yipE}
[y_i, p]=x_i^{-1}(p-p|_{x_i=0})E_{00}(i),
\end{equation}
\begin{equation}\label{yipE1}
[p, E_{00}(i)]=(p-p|_{x_i=0})E_{00}(i).
\end{equation}
In more detail, it suffices to prove the identities in the case
when $p=x_i^m$, $m\geq 1$. Then $[y_i, x_i^m]=x_i^{m-1}-x_i^my_i=
x_i^{m-1}(1-x_iy_i) = x_i^{m-1}E_{00}(i)$, and $[x_i^m ,
E_{00}(i)]= x_i^m E_{00}(i)-E_{00}(i)x_i^m = x_i^m E_{00}(i)$.

\begin{corollary}\label{a7Feb9}
\begin{enumerate}
\item The monoid ${\rm st} (x_1, \ldots , x_n)$ is an abelian
monoid each non-identity element of which is a monomorphism of the
algebra $\mS_n$ but not an automorphism. Moreover, it contains
precisely the following endomorphisms of $\mS_n$:
$$ \s_p : y_i\mapsto y_i+p_iE_{00}(i), \;\; i=1, \ldots , n,$$
where the $n$-tuple $p= (p_1, \ldots , p_n) \in K[x_1, \ldots ,
x_n]^n$ satisfies the following conditions: for each pair of
indices  $i\neq j$, 
\begin{equation}\label{pij1}
-x_j^{-1}(p_i-p_{i,j})+x_i^{-1}(p_j-p_{j,i})+p_ip_{j,i}-p_jp_{i,j}=0
\end{equation}
where $p_{i,j}:= p_i|_{x_j=0}$. \item The monoid ${\rm st} (y_1,
\ldots , y_n)$ is an abelian monoid each non-identity element of
which is a monomorphism of the algebra $\mS_n$ but not an
automorphism. Moreover, it contains precisely the following
endomorphisms of $\mS_n$:
$$ \tau_p : y_i\mapsto y_i+E_{00}(i)q_i, \;\; i=1, \ldots , n,$$
where the $n$-tuple $q= (q_1, \ldots , q_n) \in K[y_1, \ldots ,
y_n]^n$ satisfies the following conditions: for each pair of
indices  $i\neq j$, 
\begin{equation}\label{qij1}
-y_j^{-1}(q_i-q_{i,j})+y_i^{-1}(q_j-q_{j,i})+q_iq_{j,i}-q_jq_{i,j}=0
\end{equation}
where $q_{i,j}:= q_i|_{y_j=0}$.
\end{enumerate}
\end{corollary}

{\it Proof}. 1. In fact, at the beginning of the proof of Theorem
\ref{6Feb9}, we proved that each element $\s \in {\rm st} (x_1,
\ldots , x_n)$ has the form $\s = \s_p$ for {\em some} $n$-tuple
$p= (p_1, \ldots , p_n)\in K[x_1, \ldots , x_n]^n$ (there, in
proving this,  we did not use the fact the $\s $ is an
automorphism). The endomorphism $\s_p$ is well-defined iff the
elements $\s_p(y_1), \ldots , \s_p(y_n)$ commute (since
$[\s_p(y_i), \s_p(x_j)]=[\s_p(y_i), x_j]=0$ for all $i\neq j$).
Let us show that the elements $\s_p(y_1), \ldots , \s_p(y_n)$
commute iff the conditions (\ref{pij1}) hold. Moreover,  we will
prove that for each pair $i\neq j$ the condition (\ref{pij1}) is
equivalent to the condition that the elements $\s_p(y_i)$ and
$\s_p(y_j)$ commute. Indeed, using (\ref{yipE}) and (\ref{yipE1}),
we have
\begin{eqnarray*}
 0&=&  [\s_p(y_i), \s_p(y_j)]=[y_i+p_iE_{00}(i), y_j+p_jE_{00}(j)]\\
 &=&[p_i, y_j]E_{00}(i)+[y_i, p_j]E_{00}(j)+p_i[E_{00}(i), p_j]E_{00}(j)+p_j[p_i, E_{00}(j)]E_{00}(i)\\
 &=& (-x_j^{-1}(p_i-p_{i,j})+x_i^{-1}(p_j-p_{j,i})+p_ip_{j,i}-p_jp_{i,j})E_{00}(i)E_{00}(j),
\end{eqnarray*}
and so (\ref{pij1}) holds, and vice versa.

Given $\s_p, \s_{p'}\in {\rm st} (x_1, \ldots , x_n)$. Then
$$ \s_p\s_{p'}(y_i) = y_i+(p_i+p_i'-x_ip_ip_i')E_{00}(i), \;\;
i=1, \ldots , n.$$ Hence, $\s_p\s_{p'}=\s_{p'}\s_p$, and so the
monoid ${\rm st} (x_1, \ldots , x_n)$ is abelian.

 It remains to
show that each endomorphism $\s_p$ is a {\em monomorphism}, i.e.
$\ker (\s_p)=0$. Suppose that $\ker (\s_p)\neq 0$ for some $p$, we
seek a contradiction. Then $F_n\subseteq \ker (\s_p)$, since $F_n$
is the least nonzero ideal of the algebra $\mS_n$,
\cite{shrekalg}; but
$$ \s_p(E_{00}(1))= 1-x_1(y_1+p_1E_{00}(1))=
(1-x_1p_1)E_{00}(1)\neq 0, $$ a contradiction.

2. Note that $\heta ({\rm st}(x_1, \ldots , x_n))= {\rm st}(y_1,
\ldots , y_n)$ and $\heta (\s_p) = \tau_{\eta (p)}$ where $\eta
(p):= (\eta (p_1), \ldots , \eta (p_n))$ (since $\heta (\s_p)(x_i)
= \eta \s_p\eta (x_i) = \eta (y_i+p_iE_{00}(i))= x_i+E_{00}(i)\eta
(p_i))$. $\Box $

$\noindent $

For $n=1$, the conditions (\ref{pij1}) and (\ref{qij1}) are
vacuous, and so Corollary \ref{a7Feb9} takes a simpler form.
\begin{corollary}\label{b7Feb9}
\begin{enumerate}
\item ${\rm st}(x)=\{ \s_p:y\mapsto pE_{00}\, | \, p\in K[x]\}$.
\item ${\rm st}(y)=\{ \s_p:x\mapsto E_{00}q\, | \, q\in K[y]\}$.
\end{enumerate}
\end{corollary}

For each $i=1, \ldots , n$, let $G_1(i):=\Aut_{K-{\rm
alg}}(\mS_1(i))$ and $E_1(i):= \End_{K-{\rm alg}}(\mS_1(i))$.
There is a natural inclusion of groups $\prod_{i=1}^n
G_1(i)\subset G_n$. Similarly, there is a natural inclusion of
monoids $\prod_{i=1}^n E_1(i)\subset E_n$ which yields the
inclusions of submonoids:
$$ \prod_{i=1}^n {\rm st} (x_i)\subset {\rm st}(x_1, \ldots , x_n)
\;\; {\rm and } \;\;  \prod_{i=1}^n {\rm st} (y_i)\subset {\rm
st}(y_1, \ldots , y_n). $$ These inclusions are {\em not}
equalities as the following example shows.

$\noindent $

{\it Example}. Fix an {\em arbitrary} polynomial $p_i$ from the
ideal $(x_1 \cdots x_n)$ of the polynomial algebra $K[x_1, \ldots
, x_n]$, and put $p_j:= x_j^{-1}x_ip_i$ for all $j\neq i$. Then
the conditions (\ref{pij1}) hold, and so $\s_p\in E_n$ where
$p=(p_1, \ldots , p_n)$. An  element $\s_{p'}\in {\rm st}(x_1,
\ldots, x_n)$ belongs to the submonoid $\prod_{i=1}^n {\rm
st}(x_i)$ iff $p_1'\in K[x_1], \ldots , p_n'\in K[x_n]$. Now, it
is obvious that $\prod_{i=1}^n {\rm st}(x_i)\neq {\rm st}(x_1,
\ldots , x_n)$. By applying $\heta$, we see that $\prod_{i=1}^n
{\rm st}(y_i)\neq {\rm st}(y_1, \ldots , y_n)$.

Department of Pure Mathematics

University of Sheffield

Hicks Building

Sheffield S3 7RH

UK

email: v.bavula@sheffield.ac.uk

\end{document}